\documentclass[final,3p,times]{elsarticle}
\usepackage[T1]{fontenc}
\usepackage{amsmath,amssymb,mathtools,amsthm}
\usepackage{geometry}
\usepackage{microtype}
\microtypesetup{expansion=false}
\usepackage{hyperref}
\usepackage{xcolor}
\usepackage{mathrsfs}

\numberwithin{equation}{section}
\hypersetup{colorlinks=true,linkcolor=blue!45!black,citecolor=blue!45!black,
 urlcolor=blue!45!black,pdftitle={An abstract observability criterion via translation hulls, with applications to transport systems},
 pdfsubject={Kalman filtrations and real Jordan transport systems},
 bookmarksnumbered=true}
\newtheorem{theorem}{Theorem}[section]
\newtheorem{lemma}[theorem]{Lemma}
\newtheorem{proposition}[theorem]{Proposition}
\newtheorem{corollary}[theorem]{Corollary}
\theoremstyle{definition}
\newtheorem{definition}[theorem]{Definition}

\newtheorem{example}[theorem]{Example}
\theoremstyle{remark}
\newtheorem{remark}[theorem]{Remark}
\allowdisplaybreaks[2]
\begin{document}
\begin{frontmatter}

\title{An abstract observability criterion via translation hulls, with applications to transport systems}
\author[aff2,aff3]{Yacouba Simpore\corref{cor1}}
\ead{yacouba.simpore@fau.de}

\cortext[cor1]{Corresponding author.}


\affiliation[aff2]{
  organization={Universit\'e Yembila Abdoulaye Toguyeni},
  city={Fada N'Gourma},
  country={Burkina Faso}
}

\affiliation[aff3]{
  organization={Chair for Dynamics, Control, Machine Learning and Numerics,
  Alexander von Humboldt Professorship,
  Friedrich-Alexander-Universit\"at Erlangen--N\"urnberg},
  city={Erlangen},
  postcode={91058},
  country={Germany}
}
\begin{abstract}
\footnotesize
We establish an observability criterion for translation-covariant
families of operators on $L^2(\mathbb R^d)$. Under continuity of the
observed cell energies, a common summable localization bound, and
reconstruction up to a locally collectively compact error,
observability is equivalent to the absence of nonzero uniformly local
invisible states on every configuration in the full translation hull.
For one-dimensional transport systems satisfying the Kalman rank
condition, this yields localized observability in the intrinsic
Kalman--Jordan energy space for arbitrary real Jordan structure and
nonzero speeds. For unions of observed intervals with a uniform
positive lower length bound, the sharp condition is that the slowest
characteristic travels farther than the gap supremum. Observability
fails at equality even when this supremum is not attained. Duality
gives exact controllability in the dual energy space. We also obtain
sharp criteria for translating observation regions and scalar
transport with an integrable convolution kernel. The reconstruction
hypothesis is also necessary for observability. An anisotropic
heat-flux example shows that it applies beyond coercive
multiplication-plus-compact Gramian decompositions.
\end{abstract}

\end{frontmatter}
\section{Introduction}\label{sec:introduction}
Uniqueness and observability are closely related, but
they are not equivalent. The absence of nonzero
invisible solutions does not prevent sequences of
unit-energy states from becoming arbitrarily weakly
observed. On bounded domains, such sequences can often
be ruled out by a compactness--uniqueness argument.
Once an estimate has controlled high-frequency losses,
compactness produces an invisible limit of nonzero energy,
and unique continuation forces that limit to vanish
\cite{BardosLebeauRauch1992}.

On unbounded domains, this argument faces two distinct spatial
losses of compactness. Almost invisible states may escape to
infinity, where the coefficients and observation geometry can
approach limiting configurations not realized by any finite
translation. They may also spread out spatially, so that even
the largest unit-cell energy tends to zero. These mechanisms
are related to the losses of compactness analyzed by
P.-L.~Lions's concentration--compactness principle
\cite{Lions1984I}. Splitting into separated pieces, corresponding
to dichotomy in that principle, is handled by localization.
Corollary~\ref{rd:cor-ladder} provides compactly supported states
whose observation is arbitrarily small relative to their uniformly
local norm. Escape to infinity is handled by the translation hull,
which retains every limiting configuration. Spreading is handled by
cellwise normalization in $L^2_{\mathrm{uloc}}$, which keeps a
maximizing cell at unit mass. The reconstruction assumption~(R)
then supplies strong local $L^2$ compactness, ensuring that this
mass is retained in the limiting invisible profile.

\paragraph{The abstract criterion}
Let $\Gamma_\omega:\mathcal X\to\mathcal Z$, with
$\mathcal X=L^2(\mathbb R^d;\mathbb C^N)$ and
$\mathcal Z=L^2((0,T)\times\mathbb R^d;\mathbb C^m)$,
be the observation operator associated with a
configuration $\omega$. We consider the compact full
translation hull
\[
\mathcal W=
\overline{\{\theta_h\omega_0:h\in\mathbb R^d\}}
\]
of a reference configuration $\omega_0$, endowed
with a continuous translation action. The observation
operators satisfy the covariance relation
\[
\Gamma_{\theta_h\omega}\tau_h
=\tau_h\Gamma_\omega,
\]
where $\tau_hf(x)=f(x+h)$ acts on both states and
observations.

Under the assumptions of Section~\ref{abs:section},
Theorem~\ref{rev:hull} establishes
\[
\exists c>0:\quad
\|\Gamma_{\omega_0}f\|_{\mathcal Z}
\ge c\|f\|_{\mathcal X}
\quad(f\in\mathcal X)
\quad\Longleftrightarrow\quad
\ker_{L^2_{\mathrm{uloc}}}\Gamma_\omega=\{0\}
\quad(\omega\in\mathcal W).
\]
Here $L^2_{\mathrm{uloc}}$ denotes the space of
functions whose $L^2$ norms on unit cubes are
uniformly bounded, without requiring finite global
energy. The observation operators admit canonical
extensions to this space. We call the condition
on the right \emph{generalized uniqueness}.

We assume that for each bounded-support input the observed energy on
each unit cube depends continuously on $\omega$. The influence of input
cube $j$ on output cube $k$ is bounded by $\kappa(k-j)$ for a common
$\kappa\in\ell^1(\mathbb Z^d)$. The remaining assumption is a uniform
reconstruction property~(R). There are uniformly bounded operators
$B_\omega$ such that the errors $I-B_\omega\Gamma_\omega$ are
collectively compact on bounded windows for inputs of bounded uniformly
local energy. In the transport applications, a decomposition of the
Gramian $\Gamma_\omega^*\Gamma_\omega$ into a coercive multiplication
operator and a locally compact remainder provides such a reconstruction.
Conversely, (R) is necessary for observability
(Proposition~\ref{new:prop-R-necessary}). Under (H1)--(H3) on the full
hull, observability is therefore equivalent to (R) together with
generalized uniqueness (Corollary~\ref{new:cor-characterization}).

Localization transfers a failure of the global lower bound to states
whose observation is small compared with their largest unit-cell norm.
We normalize this norm to one and translate a maximizing cell to the
origin. Property~(R) gives a strongly convergent subsequence on each
bounded window, so the limit has norm one on the central cell.
Continuity of the cell energies and summability of $\kappa$ then show
that this limit is invisible for a configuration in the hull.

The requirement of uniqueness over the full translation hull is essential. As illustrated by the examples in Section~\ref{app:section}, uniqueness for the original configuration, even in the distributional sense, need not imply observability. The obstruction may instead arise from a limiting configuration at infinity.

\paragraph{Transport systems}
We now turn to the main application of the abstract
criterion: partially observed first-order transport
systems
\begin{equation}\label{eq:intro-primal}
\partial_t y+A\partial_x y=B\mathbf 1_\Omega u,
\quad (t,x)\in(0,T)\times\mathbb R,
\end{equation}
where $A\in\mathbb R^{n\times n}$ has real spectrum,
$B\in\mathbb R^{n\times m}$, and
$\Omega\subset\mathbb R$ is the control region.
The adjoint system
\begin{equation}\label{eq:intro-adjoint}
\partial_t q+A^\top\partial_x q=0,
\quad q(T)=q_T,
\end{equation}
is observed through $B^\top q$ on
$(0,T)\times\Omega$. We assume the Kalman condition
\begin{equation}\label{eq:intro-kalman}
\operatorname{rank}(B,AB,\dots,A^{n-1}B)=n,
\end{equation}
and let $K$ be the smallest integer such that
$\operatorname{rank}(B,AB,\ldots,A^{K}B)=n$.

Two features distinguish this problem. First,
the Kalman condition does not guarantee $L^2$
observability, even on the whole line. For
\[
A^\top=
\begin{pmatrix}
\lambda&1\\
0&\lambda
\end{pmatrix},
\quad B^\top=(1,0),
\]
one has, for the terminal datum $f=q_T$ and with $s=T-t$,
\[
B^\top e^{sA^\top\partial_x}f(x)
=f_1(x+\lambda s)
+s\,\partial_xf_2(x+\lambda s).
\]

Consequently, for every $T>0$, the time-integrated observation
satisfies
\[
\int_0^T\|B^\top q(T-s,\cdot)\|_2^2\,ds
\asymp_T \|f\|_{\mathcal Y^-_{K,J}}^2
=\|f_1\|_2^2+\|\partial_xf_2\|_2^2.
\]
For example, take $f_1=0$ and $f_2(x)=R^{-1/2}\phi(x/R)$ with
$\phi\in C_c^\infty(\mathbb R)$ and $\|\phi\|_2=1$. These data have
unit $L^2$ norm, while $\|\partial_x f_2\|_2=R^{-1}\|\phi'\|_2\to0$.
Whole-line observation determines, up to norm equivalence, the
intrinsic Kalman--Jordan space $\mathcal Y_{K,J}^{-}$ recalled in
Section~\ref{glob:section}, where $J$ records the Jordan structure.
Its dual $\mathcal Y_{K,J}^{+}$ is the state space for control.
A Fourier normalization identifies $\mathcal Y_{K,J}^{-}$ with $L^2$;
the uniformly local uniqueness condition applies to the normalized
datum.

Second, the observation region may be an arbitrary,
possibly non-periodic, union of intervals
\begin{equation}\label{eq:intro-gap-parameters}
\begin{aligned}
\Omega
=\bigcup_{j\in\mathbb Z}(a_j,b_j),
\quad a_j<b_j\le a_{j+1},\quad
\ell_\bullet
:=\inf_j(b_j-a_j)>0,\quad
G^\bullet
:=\sup_j(a_{j+1}-b_j)<\infty,
\end{aligned}
\end{equation}
with the case $\Omega=\mathbb R$ included by
setting $G^\bullet=0$. When
$\sigma(A)\subset\mathbb R\setminus\{0\}$,
the gap geometry and the smallest characteristic
speed determine the time scale
\begin{equation}\label{eq:intro-geometric-time}
\lambda_\bullet
:=\min_{\lambda\in\sigma(A)}|\lambda|,
\quad
T_\bullet:=\frac{G^\bullet}{\lambda_\bullet}.
\end{equation}
This is the supremum of the times that a
characteristic with the smallest speed can
spend in a gap of $\Omega$.

\paragraph{Main contributions}
Besides the abstract criterion of
Theorem~\ref{rev:hull}, our main results are
the following.
\begin{enumerate}
\item[(a)] \emph{Transport with arbitrary Jordan
structure} (Section~\ref{jo:section}).
Under \eqref{eq:intro-kalman} and
$\sigma(A)\subset\mathbb R\setminus\{0\}$,
for every $T>T_\bullet$,
Corollary~\ref{rev:boundedgaps} gives
\begin{equation}\label{eq:intro-main-observability}
c_T\|q_T\|_{\mathcal Y^-_{K,J}}^2
\le
\int_0^T\!\!\int_\Omega
|B^\top q|^2\,dx\,dt
\le
C_T\|q_T\|_{\mathcal Y^-_{K,J}}^2.
\end{equation}
Consequently, the controlled system is exactly controllable in
$\mathcal Y_{K,J}^{+}$ with $L^2$ controls.
If $G^\bullet>0$, observability fails for every $0<T\leq T_\bullet$, even when no gap attains length $G^\bullet$, as proved in Proposition~\ref{aux:necessity}. For a bounded Lipschitz weight, Theorem~\ref{jo:localized-general} characterizes observability by a positive lower bound on the time integral of the weight along every characteristic, including stationary ones.
For a single $2\times2$ Jordan block,
Proposition~\ref{rev:two-by-two} determines the optimal squared
observability cost near $T_\bullet$. With the intrinsic norm fixed in
Section~\ref{jo:cost},
$C_{\mathrm{obs}}(T,\Omega)\sim
12(1+T_\bullet^2)(T-T_\bullet)^{-3}$ as $T\downarrow T_\bullet$.

\item[(b)] \emph{Moving observation regions (Section~\ref{mov:section}).}
For $\Omega+h(t)$ with $\|h'\|_{L^\infty(0,T)}<\lambda_\bullet$,
observability holds exactly when
\[
\min_{\lambda\in\sigma(A)}|\lambda T-h(T)+h(0)|>G^\bullet.
\]
For an affine motion $h(t)=h_0+vt$, the change to a moving frame
allows arbitrary $v$ and also zero characteristic speeds. If
$G^\bullet>0$, the condition becomes
\[
T\min_{\lambda\in\sigma(A)}|\lambda-v|>G^\bullet.
\]
The relative speeds $\lambda-v$ allow a moving region to observe
stationary components.

\item[(c)] \emph{Transport with nonlocal dispersion in $\mathbb R^d$}
(Section~\ref{sec:higher-dim}).
We apply the abstract criterion to
\[
\partial_t u + V\cdot\nabla u = J*u-u
\]
with an integrable kernel $J$ and a Lipschitz observation
weight $a:\mathbb R^d\to[0,1]$. For each $T>0$, weighted $L^2$
observability holds if and only if
\[
\inf_{z\in\mathbb R^d}\int_0^T a(z+tV)\,dt>0.
\]
\item[(d)] \emph{Anisotropic heat flux (Section~\ref{sec:anisotropic-heat}).}
For the heat equation with a constant positive definite diffusion matrix,
we combine whole-space heat-flux measurements with temperature measurements
on a subset. Observability is equivalent to thickness of this subset.
The observation satisfies~(R), whereas its Gramian admits no decomposition
of the form (H4)--(H5) when the diffusion matrix has two distinct eigenvalues.
\end{enumerate}
The principal transport and nonlocal applications verify~(R) through
(H4)--(H5). The additional scope of~(R) is structural: it is preserved by
the localized output transformations of Lemma~\ref{new:lem-output}, while
(H4)--(H5) need not be, as the heat-flux example illustrates.
\paragraph{Idea of the proofs}
The transport analysis begins with a rational
Fourier multiplier $Q=DH$ identifying
$\mathcal Y^-_{K,J}$ with $L^2$. We use this order to control both frequency regimes
without assuming that the Kalman and Jordan
projectors commute.
After this normalization, the observation
operator splits into characteristic translations
and an exponentially decaying convolution
remainder. Characteristic occupation estimates
provide the coercive part of the Gramian,
while the Kalman condition and uniqueness
along characteristics, combined with a lifting
in $\mathcal S'$, yield generalized uniqueness
on the hull. For nonlocal dispersion,
generalized uniqueness follows instead from
analyticity in time in coordinates moving with the drift.

\paragraph{Related work}

\emph{Previous transport results.}
Earlier work \cite{Simpore2026} establishes localized observation for
diagonalizable systems with nonzero speeds and whole-line observation
for real Jordan blocks. Here we establish
localized observability for arbitrary real Jordan blocks on the line
and prove failure at the critical time even when the gap supremum is
not attained. The higher-dimensional transport application concerns a
scalar equation with one drift velocity and a convolution term.

\emph{Hulls and limit operators.}
Favard's work on almost-periodic differential equations
already used hulls and a separation condition for bounded
solutions of every associated homogeneous equation
\cite{Favard1928}.

In operator theory, limit operators characterize
Fredholm properties through the behavior of an
operator at infinity
\cite{RabinovichRochSilbermann2004,
LindnerSeidel2014,CLW}.
A closer connection is provided by generalized
collectively compact operator theory
\cite{ChandlerWildeZhang2002,CLW}.
In particular, \cite[Theorem~5.9]{CLW} derives a uniform
lower bound from injectivity for a translation-invariant
family of operators of the form $I-K$, under strict-topology
compactness assumptions.

Theorem~\ref{rev:hull} follows this localization and compactness
mechanism, but is formulated directly for observation maps
$\Gamma_\omega:\mathcal X\to\mathcal Z$. It uses uniqueness for the
underlying evolution equation and continuity only of the observed
cell energies. The reconstruction operators $B_\omega$ in~(R) need
only be uniformly bounded; neither their covariance nor their
continuity in $\omega$ is assumed.

The factorization $B_\omega\Gamma_\omega=I-S_\omega$ does not by itself
allow a direct application of \cite[Theorem~5.9]{CLW} to this family:
(R) imposes neither the required translation invariance nor the
strict-topology assumptions on the reconstructed operators.
Moreover, injectivity of $\Gamma_\omega$ need not imply injectivity
of a chosen reconstruction $B_\omega\Gamma_\omega$.
Under (H4)--(H5), the choice
$B_\omega=D_\omega^{-1}\Gamma_\omega^*$ gives the normalized Gramian
$D_\omega^{-1}\Gamma_\omega^*\Gamma_\omega=I+D_\omega^{-1}K_\omega$.
Those additional assumptions fail for the examples in
Sections~\ref{abs:section} and~\ref{sec:anisotropic-heat}.

For uniformly local states of infinite global energy, the identity
$\langle\Gamma_\omega^*\Gamma_\omega w,w\rangle
=\|\Gamma_\omega w\|^2$ is unavailable. Injectivity of the Gramian
therefore cannot be inferred from observation uniqueness by this
identity. Under our hypotheses, Proposition~\ref{prop:obs-vs-gram}
nevertheless proves that the two injectivity conditions, imposed on
the whole hull, are equivalent. This is a consequence of the criterion.
The contribution is the direct observation formulation with~(R),
rather than a strictly weaker injectivity condition under the same
hypotheses. Proposition~\ref{new:prop-R-necessary}, a consequence of
the lower-bound transfer in Proposition~\ref{rd:prop-rung}, also shows
that~(R) is necessary for observability. The characterization in
Corollary~\ref{new:cor-characterization} clarifies the role of this
assumption within the observation framework.

Proposition~\ref{rd:prop-rung}, which transfers lower bounds between
exponents, adapts the block argument of
\cite[Theorem~2.1]{ShinSun2009}, written there for scalar matrices,
to operator-valued blocks and to families with a common profile.

\emph{Localization and loss of compactness.}
Our approach is related to concentration--compactness
\cite{Lions1984I} and uniformly local methods
\cite{BurqJoly2016}. Limiting media also play
a role in the analysis of damped waves in
periodic and asymptotically periodic settings
\cite{OrivePazotoZuazua2001,JolyRoyer2018}.

\emph{Hyperbolic and coupled systems.}
Our results complement geometric control for
hyperbolic systems
\cite{BardosLebeauRauch1992,BastinCoron2016,
AlabauCoronOlive2017,HuOlive2024},
partial dissipation
\cite{CrinBaratDeNittiZuazua2025},
and algebraic recovery through coupling
\cite{LiardLissy2017,LissyZuazua2019,
pierre,ZUAZUA201621}.
Unlike the bounded-interval setting of
\cite{HuOlive2024}, we work on the whole line
with possibly non-periodic observation regions.
For non-diagonalizable hyperbolic evolution,
see also \cite{GarettoJahRuzhansky2018}.

\emph{Unbounded domains and moving observations.}
For the heat equation on $\mathbb R^d$,
observability of the final temperature is characterized by thickness
of the observation set
\cite{EgidiVeselic2018,
WangWangZhangZhang2019}.
Our bounded-gap regions are thick, but finite
propagation introduces a sharp positive
observation time when the set has a gap of positive length.
Moving observation regions have been studied
for waves
\cite{LeRousseauLebeauTerpolilliTrelat2017}
and for systems with stationary components
\cite{MartinRosierRouchon2013,
ChavesSilvaRosierZuazua2014,
ChavesSilvaZhangZuazua2017}.
For spectral approaches to observability,
see \cite{Wess,Miller2012Resolvent}.

\paragraph{Outline}
Sections~\ref{abs:section} and~\ref{kal:section} develop the
abstract criterion and its geometric ingredients.
Section~\ref{app:section} illustrates the necessity of the full hull. Sections~\ref{jo:section} and~\ref{mov:section} treat fixed and moving observation regions; Section~\ref{sec:higher-dim} applies the abstract criterion to transport
with nonlocal dispersion in $\mathbb R^d$. Section~\ref{sec:anisotropic-heat} illustrates the reconstruction criterion with anisotropic heat-flux measurements. The appendices collect two auxiliary lemmas and prove the necessity of the Kalman condition and the sharpness of the geometric time threshold.

\emph{Notation.}
We write $^\top$ for transposition and $^*$ for the
Hermitian or Hilbert-space adjoint, so that
$M^*=\overline M^{\,\top}$. We use
$\widehat f(\xi)=(2\pi)^{-1/2}\int_{\mathbb R}e^{-ix\xi}f(x)\,dx$ and $D_x=-i\partial_x$. For $s=T-t$,
$\widehat q(T-s,\xi)=e^{is\xi A^\top}\widehat q_T(\xi)$.

\section[An observability criterion via translation hulls in arbitrary dimension]
{An observability criterion via translation hulls in $\mathbb R^d$}
\label{abs:section}

The abstract criterion is formulated for every $d\ge1$. The transport applications of
Sections~\ref{app:section}--\ref{mov:section} use its specialization to $d=1$, retaining several
characteristic speeds and arbitrary Jordan structure; Section~\ref{sec:higher-dim} also
applies it for $d\ge2$. A translation hull records the environments seen after spatial
recentering. The translated observation sets $\Omega-x_n$ provide a
concrete example: their local compactness is explained below and
illustrates the choice of topology and the operator assumptions.

Fix $T>0$ and put
\[
 \mathcal X=L^2(\mathbb R^d;\mathbb C^N),\quad
 \mathcal Z=L^2((0,T)\times\mathbb R^d;\mathbb C^m).
\]
The operators introduced below map inputs in $\mathcal X$ to
observations in $\mathcal Z$. Scalar functions of $x$ act by
multiplication on both spaces. For $k\in\mathbb Z^d$, let
$Q_k=k+[0,1)^d$ and let $\mathbf1_k$ denote multiplication by
$\mathbf1_{Q_k}$; the corresponding observation cell is
$(0,T)\times Q_k$. For $1\le p\le\infty$, define
\begin{equation}\label{rd:amalgam}
 \|f\|_{(p)}
 =
 \bigl\|(\|\mathbf1_kf\|_{L^2})_{k\in\mathbb Z^d}\bigr\|_{\ell^p}.
\end{equation}
This norm combines the local $L^2$ sizes across spatial cells by
$\ell^p$ summation. The corresponding spaces are denoted by
$\mathcal X_p$ and $\mathcal Z_p$; they are Wiener amalgams
\cite{Feichtinger1983,FournierStewart1985,Heil2003}. In particular,
$\mathcal X_2=\mathcal X$, $\mathcal Z_2=\mathcal Z$, and
\[
 \mathcal X_\infty=L^2_{\mathrm{uloc}}(\mathbb R^d;\mathbb C^N),
 \quad
 \|f\|_{L^2_{\mathrm{uloc}}}
 =\sup_k\|f\|_{L^2(Q_k)}.
\]
Thus the uniformly local norm bounds the $L^2$ size in every cell
without requiring a finite total $L^2$ norm.

Throughout, spatial distances and lattice indices use the maximum norm
$|\cdot|$; this includes $|k-j|$, the localization profiles, and
$|x-y|$ in \eqref{abs:kernel-decay}. Vector norms are Euclidean and
matrix norms are the induced operator norms. Lipschitz constants use
the Euclidean distance:
\[
 \operatorname{Lip}(\chi)
 =\sup_{x\ne y}\frac{|\chi(x)-\chi(y)|}{|x-y|_2}.
\]
The inequality $|x|_2\le\sqrt d\,|x|$ accounts for the factors
$\sqrt d$ in Proposition~\ref{rd:prop-rung}.
\subsection{The operator assumptions}
\label{abs:assumptions}
\paragraph{Spatial translations}
We fix the convention
\begin{equation}\label{abs:translation-convention}
 (\tau_hf)(x)=f(x+h),
 \quad h,x\in\mathbb R^d,
\end{equation}
with the same translation of the spatial variable on $\mathcal Z$.
These operators are unitary on $\mathcal X$ and $\mathcal Z$ and satisfy
\[
 \tau_0=I,\quad
 \tau_h\tau_k=\tau_{h+k},\quad
 \tau_h^{-1}=\tau_{-h}.
\]
We use superscripts $\mathcal X$ and $\mathcal Z$ when it is useful
to distinguish translations of inputs and outputs. With this convention,
$\tau_{x_n}f$ expresses the profile $f$ in coordinates centered at $x_n$.
\paragraph{Translations of the environment}
An environment $\omega$ describes the spatial configuration on which
the observation operator depends. Let $\mathcal W$ be a compact
metrizable space of such configurations, equipped with maps
$\theta_h:\mathcal W\to\mathcal W$ satisfying
\begin{equation}\label{abs:environment-action}
 \theta_0\omega=\omega,\quad
\theta_h(\theta_k\omega)=\theta_{h+k}\omega
 \quad(h,k\in\mathbb R^d,\ \omega\in\mathcal W).
\end{equation}
We require $(h,\omega)\mapsto\theta_h\omega$ to be jointly continuous.
The element $\theta_h\omega$ is the environment seen after the
recentering \eqref{abs:translation-convention}. The action and the
topology of $\mathcal W$ are part of the model.
\paragraph{Example from the translated observation sets}
In one dimension, let
\[
 \Omega=\bigcup_{j\in\mathbb Z}(a_j,b_j),
 \quad a_j<b_j\leq a_{j+1},
 \quad \ell_\bullet:=\inf_j(b_j-a_j)>0.
\]
With the translation convention above, the action on observation sets is
\begin{equation}\label{abs:domain-action}
 \theta_hD=D-h,\quad
 \mathbf1_{D-h}(x)=\mathbf1_D(x+h).
\end{equation}
Indeed, $x$ belongs to the recentered set exactly when $x+h$ belongs
to the original set. Representing environments by their indicators gives
\[
 \omega_0=\mathbf1_\Omega,\quad
 \omega_n=\theta_{x_n}\omega_0=\mathbf1_{\Omega-x_n}.
\]
These translated indicators form a relatively compact family in
$L^1_{\mathrm{loc}}(\mathbb R)$. To see this, fix a bounded open interval
$I$. At most two components meeting $I$ are not contained in $I$;
all the others have length at least $\ell_\bullet$. Thus, uniformly in $n$,
\[
 \#\{\text{components of }\Omega-x_n\text{ meeting }I\}
 \leq |I|/\ell_\bullet+2.
\]
After intersection with $I$, their endpoints lie in the compact interval
$\overline I$. Passing to a subsequence fixes the number of intervals
and makes all endpoints converge. The resulting indicators converge
in $L^1(I)$, also when some limiting intervals collapse or become adjacent.
A diagonal extraction over $I=(-m,m)$ therefore gives
\[
 \omega_n\longrightarrow\omega_\infty=\mathbf1_{\Omega_\infty}
 \quad\text{in }L^1_{\mathrm{loc}}(\mathbb R).
\]
Locally, the limit is an indicator with finitely many jumps; hence
$\Omega_\infty$ may be chosen open, up to a null set.
For every bounded interval $I$, the convergence means
\[
 \|\omega_n-\omega_\infty\|_{L^1(I)}
 =|((\Omega-x_n)\triangle\Omega_\infty)\cap I|
 \longrightarrow0.
\]
Thus the configurations differ on a set of vanishing measure in each
fixed observation window. Sets are identified almost everywhere.
The orbit closure is compact in this topology, and spatial translations
act jointly continuously on it. This supplies a concrete parameter
space for the abstract criterion. Here $\theta_h$ translates the
observation geometry, while $\tau_h$ translates states and observations.

\paragraph{Operator family}
For each $\omega\in\mathcal W$, let
$\Gamma_\omega\in\mathcal L(\mathcal X,\mathcal Z)$.
In the indicator parametrization, we write
$\Gamma_D:=\Gamma_{\mathbf1_D}$.
The following assumptions describe compatibility with translations
and limits of environments, together with uniform localization and local compactness.
\begin{enumerate}
\item[(H1)] \emph{Covariance:}
\begin{equation}\label{abs:covariance}
\Gamma_{\theta_h\omega}\tau_h^{\mathcal X}
 =\tau_h^{\mathcal Z}\Gamma_\omega
 \quad(h\in\mathbb R^d,\ \omega\in\mathcal W).
\end{equation}
Recenter the input and the environment together: the resulting
observation is the original observation expressed in the new coordinates.
Equivalently,
\begin{equation}\label{abs:covariance-conjugation}
 \Gamma_{\theta_h\omega}
 =
 \tau_h^{\mathcal Z}\Gamma_\omega\tau_{-h}^{\mathcal X}.
\end{equation}
For the translated observation sets, this reads
\[
 \Gamma_{\Omega-x_n}(\tau_{x_n}^{\mathcal X}f)
 =
 \tau_{x_n}^{\mathcal Z}(\Gamma_\Omega f).
\]
\item[(H2)] \emph{Continuity of the observed energy in each cell.}
Whenever $\omega_n\to\omega$, for every $f\in\mathcal X$ of bounded
spatial support and every $k\in\mathbb Z^d$,
\begin{equation}\label{abs:strong-continuity}
 \|\mathbf1_k\Gamma_{\omega_n}f\|_{\mathcal Z}^2
 \longrightarrow \|\mathbf1_k\Gamma_\omega f\|_{\mathcal Z}^2 .
\end{equation}
In particular, strong convergence of the observations for each such
$f$ is sufficient. For translated observation sets, the condition
concerns the energy in the fixed window $(0,T)\times Q_k$.

\item[(H3)] \emph{Summable localization.} There is a common profile
$\kappa\colon\mathbb Z^d\to[0,\infty)$ such that
\begin{equation}\label{rd:localization}
 \|\mathbf1_k\Gamma_\omega\mathbf1_j\|\le\kappa(k-j),
 \qquad \sum_{n\in\mathbb Z^d}\kappa(n)<\infty
 \qquad(j,k\in\mathbb Z^d,\ \omega\in\mathcal W).
\end{equation}
The block $\mathbf1_k\Gamma_\omega\mathbf1_j$ measures how an input
in $Q_j$ contributes to the observation in $(0,T)\times Q_k$.

\item[(R)] \emph{Uniform reconstruction.}
For each $\omega\in \mathcal W$, there is a bounded linear operator
$B_\omega\colon \mathcal Z_\infty\to \mathcal X_\infty$ such that
\[
    \sup_{\omega\in \mathcal W}
    \|B_\omega\|_{\mathcal L(\mathcal Z_\infty,\mathcal X_\infty)}<\infty.
\]
Set $S_\omega=I-B_\omega\Gamma_\omega$ on $\mathcal X_\infty$, using the
uniformly local extension of $\Gamma_\omega$ given by~(H3).
For every bounded cube $Q\subset\mathbb{R}^d$, the set
\[
\bigl\{\mathbf{1}_Q S_\omega v:
\omega\in \mathcal W,\ \|v\|_{(\infty)}\leq 1\bigr\}
\]
is relatively compact in $L^2(Q;\mathbb{C}^N)$.
 The role of (R) is to prevent concentration or oscillation defects within bounded spatial windows. It ensures that sequences with vanishing observations retain strong local compactness after reconstruction.

\smallskip
\noindent\textbf{A sufficient condition for reconstruction.}
Theorem~\ref{rev:hull} uses {\rm(H1)}--{\rm(H3)} and~(R).
The following two Gramian assumptions provide a way to verify~(R)
in the transport applications.
\item[(H4)] \emph{A coercive multiplication part:}
\begin{equation}\label{abs:gramian-decomposition}
H_\omega:=\Gamma_\omega^*\Gamma_\omega=D_\omega+K_\omega,
 \quad
 (D_\omega f)(x)=d_\omega(x)f(x),
\end{equation}
where $d_\omega$ is measurable and Hermitian, with
\begin{equation}\label{abs:diagonal-ellipticity}
 0<d_\bullet I\le d_\omega(x)\le d^\bullet I
 \quad\text{a.e., uniformly in }\omega.
\end{equation}
Moreover, $d_{\omega_n}\to d_\omega$ locally in measure whenever
$\omega_n\to\omega$.
The Gramian encodes the observed energy through
$\langle H_\omega f,f\rangle_{\mathcal X}
=\|\Gamma_\omega f\|_{\mathcal Z}^2$.
Its multiplication part acts pointwise and satisfies
\[
 \langle D_\omega f,f\rangle_{\mathcal X}
 \ge d_\bullet\|f\|_{\mathcal X}^2.
\]
The remainder $K_\omega$ contains the other interactions and is
subject to the compactness requirement in (H5).

\item[(H5)] \emph{Local collective compactness:} for all bounded measurable
$I,J\subset\mathbb R^d$, the set
\begin{equation}\label{abs:collective-compactness}
 \bigcup_{\omega\in\mathcal W}
 \left\{
 \mathbf1_IK_\omega\mathbf1_Jf:
 \|f\|_{L^2(J)}\le1
 \right\}
 \quad\text{is relatively compact in }L^2(I;\mathbb C^N).
\end{equation}
Equivalently, for any sequence of environments $\omega_n$ and inputs
$f_n$ in the unit ball of $L^2(J;\mathbb C^N)$, the localized outputs
$\mathbf1_IK_{\omega_n}\mathbf1_Jf_n$ admit a strongly convergent
subsequence in $L^2(I;\mathbb C^N)$. The compactness thus holds
simultaneously over the whole family of environments.
\end{enumerate}
Both~(R) and~(H5) are \emph{collective} compactness conditions.
A family $\{G_\omega\}_{\omega\in\mathcal W}$ is \emph{collectively compact
on bounded windows} for a set $\mathcal V$ of inputs if, for every
bounded cube $I$, the set
$\bigcup_{\omega\in\mathcal W}\{\mathbf1_IG_\omega v:v\in\mathcal V\}$
is relatively compact. This is stronger than compactness of each
localized operator $\mathbf1_IG_\omega\mathbf1_J$, and the union over
$\omega$ is what Lemma~\ref{rd:lem-profile} uses, since the environment
varies along the extracted sequence. Individual compactness does not
suffice: in the scalar case, for $e_n(x)=e^{inx_1}\mathbf1_{Q_0}(x)$ and
$G_nf=\langle f,e_n\rangle e_n$, each $G_n$ has rank one, yet the
sequence $G_ne_n=e_n$ has no convergent subsequence in $L^2(Q_0)$.
Lemma~\ref{lem:cc-tools} collects the closure properties used below.

Under (H3)--(H5), Lemma~\ref{lem:gramR} gives (R) with
$B_\omega=D_\omega^{-1}\Gamma_\omega^*$.

\paragraph{A summable profile with infinite first moment}
Take a singleton hull in dimension one. Write $Uf(x)=f(x+1)$ and set
\[
 Af=f+\sum_{n\ge1}a_n U^nf,\qquad
 a_n=\frac{1}{2n(n+1)},\qquad (\Gamma f)(t,x)=Af(x).
\]
Integer translation preserves the uniformly local norm and
$\sum_{n\ge1}a_n=1/2$, so the Neumann series gives a bounded inverse
of $A$ on both $\mathcal X$ and $\mathcal X_\infty$. Thus
$Bz=A^{-1}T^{-1}\int_0^Tz(t,\cdot)\,dt$ gives (R) with $S=0$.
The family satisfies (H1)--(H3). The block of $\Gamma$ from $Q_n$
to $(0,T)\times Q_0$ has norm $\sqrt T\,a_n$ for $n\ge1$.
Hence every admissible localization profile has infinite first moment.
This example also separates (R) from (H4)--(H5). With $a_0=1$,
the block of $\Gamma^*\Gamma$ from $Q_1$ to $Q_0$ is
$T(\sum_{n\ge0}a_na_{n+1})U$ on $L^2(Q_1)$, a positive multiple
of an isometry. A multiplication operator has zero block there, so
the localized remainder cannot be compact.

\paragraph{A two-speed example for (H4)--(H5)}
Take a bounded Lipschitz weight $a:\mathbb R\to[a_0,1]$, with $a_0>0$,
and let $b$ range over its translation hull for locally uniform
convergence. In dimension one, with $N=2$ and $m=1$, consider
\[
 (\Gamma_bf)(s,x)=\sqrt{b(x)}\bigl(f_1(x+s)+f_2(x+2s)\bigr).
\]
A change of variables gives $\Gamma_b^*\Gamma_b=D_b+K_b$, where
\[
 d_b(z)=\operatorname{diag}\left(
\int_0^T b(z-s)\,ds,\ \int_0^T b(z-2s)\,ds\right),
\quad Ta_0I\leq d_b(z)\leq TI.
\]
This is the coercive multiplication part in (H4); its continuity
follows directly from locally uniform convergence of $b$.
The remainder has zero diagonal blocks and off-diagonal kernels
\[
 k_{12,b}(z,y)=\mathbf1_{\{0<y-z<T\}}b(2z-y),
 \quad k_{21,b}(z,y)=\overline{k_{12,b}(y,z)}.
\]
On bounded input and output windows, these kernels are square integrable,
so the localized remainders are Hilbert--Schmidt operators. Moreover,
they depend continuously on $b$ in Hilbert--Schmidt norm. The weight
hull is compact by Arzel\`a--Ascoli, so these localized remainders form
a norm-compact family of compact operators, which gives the collective
compactness in (H5).
The bounded profile $w=(1,-1)^\top$ nevertheless satisfies
$\Gamma_bw=0$ for every $b$: uniqueness for the full observation
remains a separate condition.

We use the notation
\[
 \|\kappa\|_1:=\sum_{n\in\mathbb Z^d}\kappa(n),
 \quad
\tau_\kappa(L):=\sum_{|n|>L}\kappa(n)
 \quad(L\in\mathbb R).
\]
Then $\tau_\kappa(L)\le\|\kappa\|_1$ and
$\tau_\kappa(L)\to0$ as $L\to\infty$.
The block estimate in (H3) gives
\[
 \|\mathbf1_k\Gamma_\omega f\|_{L^2}
 \le
 \sum_{j\in\mathbb Z^d}
 \kappa(k-j)\|\mathbf1_jf\|_{L^2}.
\]
In particular,
$\|\Gamma_\omega\|_{\mathcal L(\mathcal X,\mathcal Z)}
\le\|\kappa\|_1$, and $\Gamma_\omega$ has a canonical extension to
uniformly local inputs (Definition~\ref{def:localized} and
Lemma~\ref{rd:lem-extension}); kernels on $L^2_{\mathrm{uloc}}$ are
understood in this sense. Hypotheses (H3)--(H4) also make distant
inputs uniformly negligible in a fixed observation window
(Lemma~\ref{lem:tails}). Summability of $\kappa$ suffices for these
extension and tail estimates and for the full-hull criterion.

\subsection{The full-hull criterion}

Fix $\omega_0\in\mathcal W$ and assume, for the next theorem, that
\begin{equation}\label{abs:full-hull}
 \mathcal W
 =
 \overline{\{\theta_h\omega_0:h\in\mathbb R^d\}}.
\end{equation}
The closure is taken in the chosen topology of environments.
This is the full translation hull: it contains the reference
configuration, every translated configuration, and every limit of
such configurations.

For the one-dimensional observation sets, it becomes
\[
 \omega_0=\mathbf1_\Omega,\quad
 \mathcal W
 =
 \overline{\{\mathbf1_{\Omega-h}:h\in\mathbb R\}}
 ^{\,L^1_{\mathrm{loc}}(\mathbb R)}.
\]
Thus an environment in this hull is represented by a local
$L^1$ limit of indicators $\mathbf1_{\Omega-x_n}$. The parameters $x_n$ may remain
bounded or be unbounded. The operator family and its uniform
assumptions therefore cover the original observation geometry, its
translates, and all the limit geometries arising from spatial
recentering.

\begin{theorem}\label{rev:hull}
\label{thm:hull}
Assume (H1)--(H3) and (R), and \eqref{abs:full-hull}.
The following are equivalent:
\begin{enumerate}
\item There is $c>0$ such that
$\|\Gamma_{\omega_0}f\|_{\mathcal Z}\ge c\|f\|_{\mathcal X}$ for all
$f\in\mathcal X$.
\item For every $\omega\in\mathcal W$,
$\ker_{L^2_{\mathrm{uloc}}}\Gamma_\omega=\{0\}$.
\end{enumerate}
If the first assertion fails, there exist $\omega_\infty\in\mathcal W$ and
$w_\infty\in L^2_{\mathrm{uloc}}$ such that
\[
\Gamma_{\omega_\infty}w_\infty=0,\quad
\|w_\infty\|_{L^2(Q_0)}=\|w_\infty\|_{L^2_{\mathrm{uloc}}}=1.
\]
If it holds, there is a common $c_*>0$ such that
\[
 \|\Gamma_\omega w\|_{(p)}\ge c_*\|w\|_{(p)}
 \quad(1\le p\le\infty,\ \omega\in\mathcal W,\ w\in\mathcal X_p).
\]
In particular, the observation controls the uniformly local energy.
\end{theorem}

\begin{corollary}\label{abs:compact-family}
\label{cor:uniform}
Assume (H1)--(H3) and (R), without requiring
\eqref{abs:full-hull}. Then
\[
 \inf_{\omega\in\mathcal W}\ \inf_{f\in\mathcal X\setminus\{0\}}
 \frac{\|\Gamma_\omega f\|_{\mathcal Z}}{\|f\|_{\mathcal X}}>0
\]
if and only if $\ker_{L^2_{\mathrm{uloc}}}\Gamma_\omega=\{0\}$ for every
$\omega\in\mathcal W$. In that case a single positive lower bound
works for all $\omega\in\mathcal W$ and all $p\in[1,\infty]$.
\end{corollary}

\begin{proposition}[Observation versus Gramian uniqueness]\label{prop:obs-vs-gram}
Assume \textup{(H1)--(H3)} and \textup{(R)}. For $\omega\in\mathcal W$ let
$H_\omega=\Gamma_\omega^*\Gamma_\omega$, extended to $L^2_{\mathrm{uloc}}$ by
Lemma~\ref{rd:lem-extension}\textup{(iii)}. Then
$\ker_{L^2_{\mathrm{uloc}}}\Gamma_\omega\subset\ker_{L^2_{\mathrm{uloc}}}H_\omega$
for every $\omega\in\mathcal W$, and the following are equivalent:
\begin{enumerate}
\item[\textup{(a)}] $\ker_{L^2_{\mathrm{uloc}}}\Gamma_\omega=\{0\}$ for every $\omega\in\mathcal W$;
\item[\textup{(b)}] $\ker_{L^2_{\mathrm{uloc}}}H_\omega=\{0\}$ for every $\omega\in\mathcal W$.
\end{enumerate}
\end{proposition}

\begin{remark}\label{rem:obs-vs-gram}
The implication \textup{(a)}$\Rightarrow$\textup{(b)} uses
Corollary~\ref{abs:compact-family}, hence the profile extraction of
Lemma~\ref{rd:lem-profile}. It is not obtained from the global energy
identity, which is unavailable for $w\notin L^2$.
In applications, \textup{(a)} follows directly from uniqueness for the
evolution equation; see Lemmas~\ref{app:distribu}
and~\ref{rev:uniqueness}, and Step~3 of the proof of
Theorem~\ref{thm:nonlocal}. The Gramian condition \textup{(b)} follows
afterwards. The equivalence concerns the conditions on the entire
family; no strict inclusion of the two kernels for an individual
environment is asserted.
Proposition~\ref{new:prop-R-necessary} below shows that observability
also implies~(R). Corollary~\ref{new:cor-characterization} combines
this necessity with Theorem~\ref{rev:hull} to give a characterization
under (H1)--(H3) on the full hull.
\end{remark}

The localization tools and the proofs are given in
Sections~\ref{rd:sec-amalgams} and~\ref{rd:sec-proofs}.
\begin{remark}[Why several exponents]
For $1\le p\le q\le\infty$ and $f$ supported in $n$ cells,
\begin{equation}\label{eq:normcomp}
\|f\|_{(q)}\le\|f\|_{(p)}\le n^{1/p-1/q}\|f\|_{(q)}.
\end{equation}
A failure of the $L^2$ lower bound produces inputs with small $L^2$ observation.
Such inputs may be spread over many cells, so their largest cell norm may be much smaller
than their $L^2$ norm, and a direct norm comparison does not reach the uniformly local norm.
Proposition~\ref{rd:prop-rung} instead applies the lower bound in one exponent to cutoffs
on blocks of $R^d$ cells. The commutator error is controlled by the summable profile alone,
and the two exponents are compared only inside each block, at the fixed cost $R^{d|1/p-1/q|}$.
\end{remark}
\subsection{A concrete sufficient operator class}
\begin{proposition}
\label{abs:concrete}
Let $\mathcal W$ be the compact metrizable parameter space with
continuous $\mathbb R^d$-action introduced above. Let
$V_1,\dots,V_p\in\mathbb R^d$ be distinct, and suppose that, for
$\omega\in\mathcal W$,
\begin{equation}
\label{abs:class}
\begin{aligned}
  (\Gamma_\omega f)(t,x)
  =\sum_{\ell=1}^p c_{\ell,\omega}(t,x)f(x+tV_\ell)+\int_{\mathbb R^d}k_\omega(t,x,y)f(y)\,dy,
\end{aligned}
\end{equation}
with measurable coefficients $c_{\ell,\omega}:(0,T)\times\mathbb R^d\to\mathbb C^{m\times N}$
and kernels $k_\omega:(0,T)\times\mathbb R^d\times\mathbb R^d\to\mathbb C^{m\times N}$. Assume
that $\sup_{\ell,\omega}\|c_{\ell,\omega}\|_{L^\infty}<\infty$ and that
\begin{equation}
\label{abs:kernel-decay}
\|k_\omega(t,x,y)\|\le Ce^{-\alpha|x-y|}\quad(\alpha>0),
\end{equation}
uniformly in $\omega$. Assume the covariance relations
\begin{equation}
\label{eq:coefficient-covariance}
\begin{aligned}
  c_{\ell,\theta_h\omega}(t,x)
  =c_{\ell,\omega}(t,x+h),\quad
  k_{\theta_h\omega}(t,x,y)
  =k_\omega(t,x+h,y+h),
\end{aligned}
\end{equation}
for all $h\in\mathbb R^d$, $\omega\in\mathcal W$, $1\le\ell\le p$ and a.e.\ $(t,x,y)$. When
$\omega_n\to\omega$, assume local convergence in measure of the coefficients on
$(0,T)\times\mathbb R^d$, and local $L^2$ convergence of the kernels on $(0,T)\times\mathbb
R^d\times\mathbb R^d$. Put
\begin{equation}
\label{abs:class-d}
d_\omega(z)=\sum_{\ell=1}^p\int_0^Tc_{\ell,\omega}(t,z-tV_\ell)^*c_{\ell,\omega}(t,z-tV_\ell)\,dt
,
\end{equation}
and assume $d_\omega\ge d_\bullet I>0$ uniformly. Write
$T_{\ell,\omega}f(t,x)=c_{\ell,\omega}(t,x)f(x+tV_\ell)$, and suppose that the total
cross-velocity family
\[
R_\omega=\sum_{\ell\ne r}T_{\ell,\omega}^*T_{r,\omega}
\]
is locally collectively compact in the sense of \textup{(H5)}. Then \textup{(H1)--(H5)} hold.
The last assumption is automatic if $d=1$, and it holds in any dimension if
$c_{\ell,\omega}^*c_{r,\omega}=0$ a.e.\ for $\ell\ne r$, in particular if $p=1$. A dense orbit
is needed only to invoke Theorem~\ref{rev:hull}, not Corollary~\ref{abs:compact-family}.
\end{proposition}

\begin{proof}
Let $\rho=T\max_\ell|V_\ell|$ and
$c^\bullet=\sup_{\ell,\omega}\|c_{\ell,\omega}\|_\infty$.
The block from $Q_j$ to $(0,T)\times Q_k$ of a transport term vanishes
unless $|k-j|\le\rho+1$, and its norm is always at most
$T^{1/2}c^\bullet$. The integral term has block norm at most
$T^{1/2}Ce^\alpha e^{-\alpha|k-j|}$, by Cauchy--Schwarz on the input unit cell.
Thus (H3) holds with
\[
 \kappa(m)=T^{1/2}\bigl(p c^\bullet\mathbf1_{\{|m|\le\rho+1\}}
+Ce^\alpha e^{-\alpha|m|}\bigr).
\]
On bounded input and output windows, kernel convergence implies
operator-norm convergence of the integral terms. Exponential tails
control the remaining output. For each transport term, bounded
convergence in measure against the integrable weight
$|f(x+tV_\ell)|^2$ gives strong convergence of the observations on
bounded-support inputs and hence (H2). The same argument in
\eqref{abs:class-d} gives local $L^1$ convergence of $d_\omega$;
boundedness and the assumed lower bound prove (H4).

Let $L_\omega$ be the integral term and
$T_\omega=\sum_\ell T_{\ell,\omega}$. Then
\[
 K_\omega=R_\omega+T_\omega^*L_\omega
+L_\omega^*T_\omega+L_\omega^*L_\omega.
\]
The three terms involving $L_\omega$ have exponentially bounded
kernels, uniformly in $\omega$. For the last term, use
\[
 \int_{\mathbb R^d}e^{-\alpha|x-z|}e^{-\alpha|x-y|}\,dx
 \le e^{-\alpha|z-y|/2}\int_{\mathbb R^d}e^{-\alpha|x-z|/2}\,dx
 =C_d\alpha^{-d}e^{-\alpha|z-y|/2}.
\]
Their kernels depend continuously on $\omega$ in local $L^2$:
truncate intermediate spatial integrations, use the assumed
convergences on bounded sets, then discard uniformly exponentially
small tails. Thus, for every bounded $I,J$,
\[
 \|\mathbf1_I\bigl[(K_{\omega_n}-R_{\omega_n})-(K_\omega-R_\omega)\bigr]
       \mathbf1_J\|_{\mathrm{HS}}\longrightarrow0
 \qquad(\omega_n\to\omega).
\]
Each block is compact, and its continuous dependence on the compact
parameter space makes the family compact in operator norm.
Lemma~\ref{lem:cc-tools}\textup{(a)} gives collective compactness. Together
with the hypothesis on $R_\omega$, Lemma~\ref{lem:cc-tools}\textup{(b)}
proves (H5).

If $d=1$ and $\ell\ne r$, the change of variables
$z=x+tV_\ell$, $z'=x+tV_r$ gives the cross-term kernel
\[
 \frac{\mathbf1_{(0,T)}(t)}{|V_r-V_\ell|}
 c_{\ell,\omega}(t,z-tV_\ell)^*
 c_{r,\omega}(t,z-tV_\ell),\quad
 t=\frac{z'-z}{V_r-V_\ell}.
\]
It is bounded, has finite range, and depends continuously on $\omega$
in local $L^2$. By Lemma~\ref{lem:cc-tools}\textup{(a)}, this proves
the asserted automatic compactness.
Under channel orthogonality, each cross term vanishes in any dimension. Finally, for
$f\in\mathcal X$, \eqref{eq:coefficient-covariance} and the change of variables $y\mapsto y+h$
give
\begin{align*}
(\Gamma_{\theta_h\omega}\tau_hf)(t,x)
&=\sum_{\ell=1}^p c_{\ell,\omega}(t,x+h)f(x+h+tV_\ell)
  +\int_{\mathbb R^d}k_\omega(t,x+h,y+h)f(y+h)\,dy\\
&=(\Gamma_\omega f)(t,x+h)
 =(\tau_h\Gamma_\omega f)(t,x),
\end{align*}
which is \textup{(H1)}.
\end{proof}

\begin{remark}\label{abs:remark-multiplier}
Injectivity on global $L^2$ alone is insufficient. For example, the
scalar Fourier multiplier with symbol
\[
 h(\xi)=\frac{|\xi-\xi_*|_2^2}{1+|\xi-\xi_*|_2^2}
\]
is injective on $L^2$ but has lower norm zero, whereas the plane wave
$e^{i\xi_*\cdot x}$ is a nonzero uniformly local null profile. Here
$|\cdot|_2$ denotes the Euclidean norm.
The criterion requires uniqueness for observation maps on every
configuration of the hull. This is the form supplied by the
transport uniqueness argument below. The relation with injectivity of
$\Gamma_\omega^*\Gamma_\omega$ on uniformly local states is described in
Proposition~\ref{prop:obs-vs-gram} and Remark~\ref{rem:obs-vs-gram}.
The role of (R) can be seen even on a singleton hull. For
$\Gamma f(t)=e^\Delta f$, the Gaussian kernel satisfies (H3) and its
Fourier multiplier $e^{-|\xi|_2^2}$ is injective on tempered
distributions, hence on $L^2_{\mathrm{uloc}}$. Indeed one can divide
by this multiplier on every compact frequency window. Yet the inputs
$e^{inx_1}\varphi(x)$, for a fixed $\varphi\in C_c^\infty(Q_0^\circ)$
of $L^2$ norm one, have observation norm tending to zero by
Plancherel's theorem and dominated convergence. Assumption (R)
excludes this loss of compactness within a cell. An explicit
observability cost in terms of geometry requires a separate estimate.
\end{remark}

\begin{remark}[The obstruction for several velocities]
\label{rd:rem-velocities}
For $d\ge2$, a cross term has the form
\[
 (Kf)(z)=\int_0^T a(t,z)f(z+t(V_r-V_\ell))\,dt.
\]
If this individual operator is nonzero, choose a compactly supported
$\psi$ and a bounded output window on which $K\psi\ne0$, and take
$0\ne\eta\perp(V_r-V_\ell)$. Then
$f_n=e^{in\eta\cdot x}\psi\rightharpoonup0$, while
$Kf_n=e^{in\eta\cdot x}K\psi$. The localized output norm stays positive, so $K$ is not locally compact.

Orthogonality is therefore a useful sufficient condition, but it is not necessary for the total remainder: different cross terms may cancel. For instance, for scalar input and two output components,
\[
 (\Gamma f)(t,x)=
 \begin{pmatrix}f(x)+f(x+tv)\\ f(x+tv)-f(x+2tv)\end{pmatrix},
 \quad 0\ne v\in\mathbb R^d,
\]
has $\Gamma^*\Gamma=4TI$, although its three channel vectors are not pairwise orthogonal. This explains why Proposition~\ref{abs:concrete} allows collective compactness of the total cross-velocity family.
For the general transport systems considered later, the
one-dimensional change of variables is what verifies compactness; the theorem does not assert their multidimensional analogue.
\end{remark}

\subsection{Localization in amalgam spaces}
\label{rd:sec-amalgams}

The proof of Theorem~\ref{rev:hull} requires almost invisible inputs
whose largest cell norm can be normalized to one. We first extend the
operators to uniformly local inputs. We then transfer lower bounds
between exponents under (H3) and obtain compactly supported inputs
for profile extraction.

\paragraph{The uniformly local extension}
At the endpoint $p=\infty$, compactly supported inputs are not dense.
For instance, in the scalar case $w\equiv1$ belongs to
$\mathcal X_\infty$, but for $P_L=\mathbf1_{[-L,L]^d}$,
\[
 \|w-P_Lw\|_{(\infty)}=1,
 \quad P_Lw\longrightarrow w\quad\text{in }L^2_{\mathrm{loc}}.
\]
The extension below is therefore defined cellwise; localization will ensure that it agrees with local limits of truncations.

\begin{definition}[Localized operators and their uniformly local extension]
\label{def:localized}
Let $\mathcal E,\mathcal F\in\{\mathcal X,\mathcal Z\}$ and let
$\kappa:\mathbb Z^d\to[0,\infty)$ be summable. An operator
$G\in\mathcal L(\mathcal E,\mathcal F)$ is \emph{$\kappa$-localized} if
\[
 \|\mathbf1_kG\mathbf1_j\|\leq\kappa(k-j)
 \quad(j,k\in\mathbb Z^d).
\]
For $w\in\mathcal E_\infty$, define $Gw$ by
\begin{equation}
 \mathbf1_kGw:=\sum_{j\in\mathbb Z^d}
                 \mathbf1_kG(\mathbf1_jw)
 \quad(k\in\mathbb Z^d).
 \label{eq:cellwise-extension}
\end{equation}
This series converges absolutely in the $L^2$ space of the $k$-th
output cell, since its terms are bounded by
$\kappa(k-j)\|w\|_{(\infty)}$. Under \textup{(H3)}, it defines
$\Gamma_\omega w$ for every $w\in L^2_{\mathrm{uloc}}$.
Throughout, the notation
\[
 \ker_{L^2_{\mathrm{uloc}}}\Gamma_\omega
 :=\{w\in L^2_{\mathrm{uloc}}:\Gamma_\omega w=0\}
\]
uses this extension. For these observation maps, Lemma~\ref{lem:tails}
identifies it with the local limit of bounded-support truncations.
\end{definition}

\begin{lemma}
\label{rd:lem-extension}
\label{lem:extension}
Let $G\in\mathcal L(\mathcal E,\mathcal F)$ be $\kappa$-localized.
\begin{enumerate}
\item[\textup{(i)}]
For every $p\in[1,\infty]$, \eqref{eq:cellwise-extension} defines
a bounded operator $\mathcal E_p\to\mathcal F_p$ of norm at most
$\|\kappa\|_1$, agreeing with $G$ on $\mathcal E_p\cap\mathcal E$.
For $p<\infty$, it is the unique bounded extension from compactly
supported inputs.
\item[\textup{(ii)}]
If $w\in\mathcal E_\infty$, $L\in\mathbb R$, and
$\mathbf1_jw=0$ whenever $|j-k|\leq L$, then
\[
 \|\mathbf1_kGw\|\leq\tau_\kappa(L)\|w\|_{(\infty)}.
\]
\item[\textup{(iii)}]
If $G_1\in\mathcal L(\mathcal E,\mathcal F)$ is $\kappa_1$-localized
and $G_2\in\mathcal L(\mathcal F,\mathcal E')$ is
$\kappa_2$-localized, with $\mathcal E'\in\{\mathcal X,\mathcal Z\}$,
then $G_2G_1$ is $(\kappa_2*\kappa_1)$-localized and
\[
 (G_2G_1)w=G_2(G_1w)\qquad(w\in\mathcal E_\infty).
\]
Moreover, $G^*$ is $\widetilde\kappa$-localized, where
$\widetilde\kappa(m)=\kappa(-m)$.
\item[\textup{(iv)}]
Under \textup{(H3)--(H4)}, the operators
$\Gamma_\omega^*$, $H_\omega$, $D_\omega$, and $K_\omega$
have the respective localization profiles
\[
 \widetilde\kappa,\quad
 \widetilde\kappa*\kappa,\quad
 d^\bullet\delta_0,\quad
 \kappa_K:=\widetilde\kappa*\kappa+d^\bullet\delta_0,
\]
uniformly in $\omega$. For every $w\in L^2_{\mathrm{uloc}}$,
\[
 H_\omega w=\Gamma_\omega^*(\Gamma_\omega w)
=d_\omega w+K_\omega w.
\]
\end{enumerate}
\end{lemma}

\begin{proof}
\textbf{(i)} The cellwise estimate
\[
 \|\mathbf1_kGw\|
 \leq\sum_j\kappa(k-j)\|\mathbf1_jw\|
\]
and Young's inequality on $\mathbb Z^d$ give the bound for every $p$.
If $w\in\mathcal E_p\cap\mathcal E$, its cell decomposition converges
in $\mathcal E$, so applying $G$ and restricting to the $k$-th output
cell gives \eqref{eq:cellwise-extension}. For $p<\infty$, uniqueness
follows from the density of compactly supported inputs in $\mathcal E_p$.

\textbf{(ii)} Only indices with $|j-k|>L$ contribute to
\eqref{eq:cellwise-extension}.

\textbf{(iii)} On $L^2$ inputs, inserting the intermediate cell
decomposition gives
\[
 \|\mathbf1_kG_2G_1\mathbf1_j\|
 \leq\sum_l\kappa_2(k-l)\kappa_1(l-j)
 =(\kappa_2*\kappa_1)(k-j).
\]
For $w\in\mathcal E_\infty$, the double series is absolutely convergent:
\[
 \sum_{l,j}\|\mathbf1_kG_2\mathbf1_lG_1\mathbf1_jw\|
 \leq\|\kappa_2\|_1\|\kappa_1\|_1\|w\|_{(\infty)}.
\]
Summing first in $j$ gives $\mathbf1_kG_2(G_1w)$, while summing
first in $l$ gives $\mathbf1_k(G_2G_1)w$. The adjoint profile follows
from $\|\mathbf1_jG^*\mathbf1_k\|=\|\mathbf1_kG\mathbf1_j\|$.

\textbf{(iv)} Apply (iii) to $\Gamma_\omega^*\Gamma_\omega$.
Multiplication by $d_\omega$ acts within each cell with norm at most
$d^\bullet$, and $K_\omega=H_\omega-D_\omega$. These observations
give the stated profiles. Both Gramian identities hold cellwise and
therefore also on $L^2_{\mathrm{uloc}}$.
\end{proof}

The next lemma expresses the spatial meaning of localization: inputs
far from a fixed output window have a uniformly small effect there.
It will allow local convergence arguments even when the inputs have
infinite global energy.

\begin{lemma}
\label{lem:tails}
Assume \textup{(H3)--(H4)} and put $P_L=\mathbf1_{[-L,L]^d}$.
For every bounded measurable $I\subset\mathbb R^d$, there is a function
$\varepsilon_I$, depending only on $I,d,\kappa,d^\bullet$, with
$\varepsilon_I(L)\to0$ as $L\to\infty$, such that
\begin{equation}
 \begin{aligned}
 \|\mathbf1_I\Gamma_\omega(1-P_L)w\|_{\mathcal Z}
  +\|\mathbf1_IK_\omega(1-P_L)w\|_{\mathcal X}
  \leq\varepsilon_I(L)\|w\|_{(\infty)}\quad
  (\omega\in\mathcal W,\ L\geq1,\ w\in L^2_{\mathrm{uloc}}).
 \end{aligned}
 \label{abs:uniform-tails}
\end{equation}
In particular, $\mathbf1_I\Gamma_\omega(P_Lw)\to\mathbf1_I\Gamma_\omega w$
and $\mathbf1_IK_\omega(P_Lw)\to\mathbf1_IK_\omega w$ in $L^2$.
The estimate for $\Gamma_\omega$ alone requires only the summable block
bound in \eqref{rd:localization}.
\end{lemma}

\begin{proof}
By Lemma~\ref{rd:lem-extension} (iv), the two operator families are
localized with profiles $\kappa$ and $\kappa_K$. Choose an integer
$R_I\geq0$ such that $I\subset\bigcup_{|k|\leq R_I}Q_k$.
If $Q_j$ meets the complement of $[-L,L]^d$, then $|j|>L-1$.
Thus, for $|k|\leq R_I$, the input $(1-P_L)w$ vanishes in every
cell with $|j-k|\leq L-R_I-1$.
For either operator $G_\omega\in \{\Gamma_\omega, K_\omega\}$, with its corresponding profile
$\rho\in\{\kappa,\kappa_K\}$, Lemma~\ref{rd:lem-extension} (ii) gives
\[
 \|\mathbf1_kG_\omega(1-P_L)w\|
 \leq\tau_\rho(L-R_I-1)\|w\|_{(\infty)}.
\]
Summing the squares over the $(2R_I+1)^d$ output cells proves
\eqref{abs:uniform-tails} with
\[
 \varepsilon_I(L)=(2R_I+1)^{d/2}
 \bigl(\tau_\kappa(L-R_I-1)+\tau_{\kappa_K}(L-R_I-1)\bigr).
\]
Both tails tend to zero by summability. The convergence of truncations
follows by linearity, since $P_Lw\in\mathcal X$.
\end{proof}

\paragraph{Transfer between exponents}
Let $\mathcal E,\mathcal F\in\{\mathcal X,\mathcal Z\}$. For a localized
operator $G\in\mathcal L(\mathcal E,\mathcal F)$ put
\[
 \nu_p(G):=\inf_{f\in\mathcal E_p\setminus\{0\}}
              \frac{\|Gf\|_{(p)}}{\|f\|_{(p)}} .
\]
The following argument groups cells into larger blocks, as in the
localization proof of \cite[Theorem~2.1]{ShinSun2009} for scalar matrices.

\begin{proposition}\label{rd:prop-rung}
Let $\mathcal G\subset\mathcal L(\mathcal E,\mathcal F)$ have a common
summable localization profile $\kappa$. If $\nu_p(G)\ge a>0$ for all
$G\in\mathcal G$ at one exponent $p\in[1,\infty]$, then
\[
 \inf_{G\in\mathcal G,\;1\le q\le\infty}\nu_q(G)>0 .
\]
\end{proposition}
\begin{proof}
Fix $G\in\mathcal G$ and $q\in[1,\infty]$, and write $G_{kj}=\mathbf 1_kG\mathbf 1_j$, so that
$\|G_{kj}\|\le\kappa(k-j)$. Throughout, $1/\infty=0$ and $s_+=\max\{s,0\}$.

\emph{Step 1: block norms.} For an integer $R\ge1$, the blocks
$D_l=Rl+\{0,\dots,R-1\}^d$, $l\in\mathbb Z^d$, partition $\mathbb Z^d$.
For an input or output $v$ put
\[
F_l(v)=\big\|(\|\mathbf 1_jv\|_{L^2})_{j\in D_l}\big\|_{\ell^p},\qquad
M(v)=\|(F_l(v))_l\|_{\ell^q}.
\]
Since $\#D_l=R^d$, Hölder's inequality on each block gives
\begin{equation}\label{eq:Mq}
R^{-d(1/q-1/p)_+}\|v\|_{(q)}\le M(v)\le R^{d(1/p-1/q)_+}\|v\|_{(q)}.
\end{equation}
In particular $M(v)<\infty$ and $M(Gv)<\infty$ for $v\in\mathcal E_q$.

\emph{Step 2: cutoffs.} Choose a Lipschitz function $\chi:\mathbb R^d\to[0,1]$, equal to one
on $[0,1]^d$ and supported in $[-1,2]^d$, and let $\chi_l$ multiply the cell $j$, of an input
or an output, by $\chi(j/R-l)$. Then $\chi_l=1$ on $D_l$, and $\chi_l$ vanishes outside
$\bigcup_{b\in B}D_{l+b}$, where $B=\{-1,0,1,2\}^d$. Hence
\begin{equation}\label{eq:chiM}
M(v)\le\big\|(\|\chi_lv\|_{(p)})_l\big\|_{\ell^q}\le 4^d M(v).
\end{equation}
Each $\chi_lv$ has bounded support, so $\chi_lv\in\mathcal E\cap\mathcal E_p$,
and the extensions of $G$ agree on it by Lemma~\ref{lem:extension}\,(i).

\emph{Step 3: commutators.} For $n\in\mathbb Z^d$ let $G^{(n)}$ be the $n$th diagonal,
$\mathbf 1_kG^{(n)}v=G_{k,k-n}\mathbf 1_{k-n}v$. By Lemma~\ref{lem:extension},
$Gv=\sum_nG^{(n)}v$ with absolute convergence in every cell, and
\[
\mathbf 1_k[G^{(n)},\chi_l]v
=\big(\chi(\tfrac{k-n}{R}-l)-\chi(\tfrac kR-l)\big)\,G_{k,k-n}\mathbf 1_{k-n}v .
\]
The scalar factor is at most $c_\chi\min\{1,|n|/R\}$, with
$c_\chi=\max\{1,\sqrt d\,\mathrm{Lip}(\chi)\}$. It vanishes unless $k$ or $k-n$ lies in
$R(l+[-1,2]^d)$. Equivalently, the input index $j=k-n$ lies in
$\bigcup_{b\in B'}D_{l+b}$ or in $\bigcup_{b\in B'}D_{l-\lfloor n/R\rfloor+b}$,
where $B'=\{-2,\dots,2\}^d$ and the floor is taken componentwise.
Since $k\mapsto k-n$ is a bijection,
\[
\|[G^{(n)},\chi_l]v\|_{(p)}\le c_\chi\min\{1,\tfrac{|n|}{R}\}\,\kappa(n)
\sum_{b\in B'}\big(F_{l+b}(v)+F_{l-\lfloor n/R\rfloor+b}(v)\big).
\]
For every fixed shift $c$, $\|(F_{l+c}(v))_l\|_{\ell^q}=M(v)$. Summing over $n$ by
Minkowski's inequality gives
\begin{equation}\label{eq:comm}
\big\|(\|[G,\chi_l]v\|_{(p)})_l\big\|_{\ell^q}\le C_\chi\,\eta_R\,M(v),\qquad
\eta_R=\sum_{n\in\mathbb Z^d}\kappa(n)\min\{1,|n|/R\},
\end{equation}
with $C_\chi=2\cdot5^dc_\chi$. Since $\kappa\in\ell^1$, $\eta_R\to0$ by dominated convergence.
Only the summability of $\kappa$ is used here. All sums converge absolutely,
also when $p$ or $q$ is infinite.

\emph{Step 4: conclusion.} The lower bound at the exponent $p$, applied to $\chi_lv$, gives
\[
a\|\chi_lv\|_{(p)}\le\|G(\chi_lv)\|_{(p)}\le\|\chi_lGv\|_{(p)}+\|[G,\chi_l]v\|_{(p)}.
\]
Take the $\ell^q$ norm in $l$, and use \eqref{eq:chiM} for $v$ and for $Gv$, together
with \eqref{eq:comm}:
\[
aM(v)\le 4^dM(Gv)+C_\chi\eta_RM(v).
\]
Fix an integer $R$, depending only on $a$, $\kappa$, $d$ and $\chi$, such that
$C_\chi\eta_R\le a/2$. Since $M(v)<\infty$, we get $M(v)\le 2\cdot4^d a^{-1}M(Gv)$.
Applying \eqref{eq:Mq} to $v$ and to $Gv$ then gives
\begin{equation}\label{eq:transfer}
\|Gv\|_{(q)}\ge\frac{a}{2\cdot4^d}\,R^{-d|1/p-1/q|}\,\|v\|_{(q)}\qquad(v\in\mathcal E_q).
\end{equation}
Since $|1/p-1/q|\le1$ and $R$ depends neither on $G$ nor on $q$, we obtain
$\nu_q(G)\ge a/(2\cdot4^dR^d)$ for all $G\in\mathcal G$ and $q\in[1,\infty]$.
\end{proof}

\begin{corollary}\label{rd:cor-ladder}
Let $\mathcal G\subset\mathcal L(\mathcal E,\mathcal F)$ have a common
summable localization profile $\kappa$, and let $p,q\in[1,\infty]$.
If $\inf_{G\in\mathcal G}\nu_p(G)=0$, then for every $\varepsilon>0$
there are $G\in\mathcal G$ and a nonzero compactly supported
$v\in\mathcal E$ such that
\[
 \|Gv\|_{(q)}\leq\varepsilon\|v\|_{(q)} .
\]
\end{corollary}

\begin{proof}
By Proposition~\ref{rd:prop-rung} applied with $q$ as the starting
exponent, $\inf_{G\in\mathcal G}\nu_q(G)=0$. For $q<\infty$, the
claim follows by density of compactly supported inputs in $\mathcal
E_q$ and boundedness of the extension. For $q=\infty$, take
$\|u\|_{(\infty)}=1$ and $\|Gu\|_{(\infty)}<\varepsilon/4$ for some
$G\in\mathcal G$. Choose $j_0$ with $\|\mathbf1_{j_0}u\|>1/2$.
Let $\psi$ be compactly supported and Lipschitz, with $0\le\psi\le1$
and $\psi(0)=1$, and let $\psi_R$ multiply cell $j$ by
$\psi((j-j_0)/R)$. The difference of its multipliers on cells $j,k$
is at most $C\min(1,|k-j|/R)$. Summing the block bounds gives
$\|[G,\psi_R]u\|_{(\infty)}\le C\eta_R$, where $\eta_R\to0$ is as
in the proof above. For $R$ with $C\eta_R<\varepsilon/4$ the compactly
supported $v=\psi_Ru$ has $\|v\|_{(\infty)}>1/2$ and
$\|Gv\|_{(\infty)}<\varepsilon/2$, as required.
\end{proof}

\paragraph{Collective compactness}
The following closure properties are used in Lemma~\ref{lem:gramR},
Proposition~\ref{abs:concrete} and Theorem~\ref{thm:nonlocal}.

\begin{lemma}\label{lem:cc-tools}
Let $\mathcal W$ be compact, let $I,J\subset\mathbb R^d$ be bounded, and
let $\{G_\omega\}_{\omega\in\mathcal W}$ be localized operators.
\begin{enumerate}
\item[\textup{(a)}] If $\omega\mapsto\mathbf1_IG_\omega\mathbf1_J$ is
continuous in operator norm with compact values, then
$\{\mathbf1_IG_\omega\mathbf1_J\}_{\omega}$ is collectively compact on the
unit ball of $L^2(J)$.
\item[\textup{(b)}] Finite sums of collectively compact families are
collectively compact. If each family $\{G^{(m)}_\omega\}$ is collectively
compact on the unit ball of $L^2(J)$ and
$\sup_\omega\|\mathbf1_I(G_\omega-G^{(m)}_\omega)\mathbf1_J\|\to0$, then
$\{\mathbf1_IG_\omega\mathbf1_J\}_\omega$ is collectively compact on the
same ball.
\item[\textup{(c)}] Let $\{G_\omega\}$ be collectively compact on a set
$\mathcal V$, and let $m_\omega$ be matrix-valued multipliers with
$\sup_\omega\|m_\omega\|_{L^\infty}<\infty$ and $m_{\omega_n}\to m_\omega$
locally in measure whenever $\omega_n\to\omega$. Then
$\{m_\omega G_\omega\}$ is collectively compact on $\mathcal V$.
\item[\textup{(d)}] If $\{G_\omega\}$ has a common summable profile and is
collectively compact on the unit ball of $L^2(J)$ for every bounded $J$,
then it is collectively compact on the unit ball of $\mathcal X_\infty$.
\end{enumerate}
\end{lemma}

\begin{proof}
(a) The set $\{\mathbf1_IG_\omega\mathbf1_J\}_\omega$ is a norm-compact set
of compact operators. A finite $\varepsilon$-net of it shows that the
union of the images of the unit ball is totally bounded.
(b) The same $\varepsilon$-net argument applies.
(c) Given $(\omega_n,v_n)$ with $v_n\in\mathcal V$, extract a subsequence
with $\omega_n\to\omega$ and $\mathbf1_IG_{\omega_n}v_n\to g$ in $L^2(I)$.
Then $m_{\omega_n}\mathbf1_IG_{\omega_n}v_n\to m_\omega g$, by the uniform
bound and dominated convergence along a further subsequence.
(d) Put $P_L=\mathbf1_{[-L,L]^d}$. By
Lemma~\ref{rd:lem-extension}\textup{(ii)}, as in the proof of
Lemma~\ref{lem:tails}, the common summable profile gives
\[
 \sup_{\omega\in\mathcal W}\sup_{\|v\|_{(\infty)}\le1}
 \|\mathbf1_I G_\omega(I-P_L)v\|_{L^2(I)}
 \le\eta_I(L)\longrightarrow0.
\]
For fixed $L$, the inputs $P_Lv$ lie in a common bounded ball of
$L^2([-L,L]^d)$, whose images are relatively compact by assumption.
The uniform tail estimate extends this compactness to the original
inputs.
\end{proof}

\subsection{Extraction of profiles and proofs of the criterion}
\label{rd:sec-proofs}
\begin{lemma}\label{lem:gramR}
Assume (H3)--(H5). Then (R) holds with $B_\omega = D_\omega^{-1}\Gamma_\omega^*$
and $S_\omega = -D_\omega^{-1}K_\omega$ on $\mathcal X_\infty$.
\end{lemma}

\begin{proof}
By Lemma~\ref{rd:lem-extension} (iii), the reflected localization profile bounds $\Gamma_\omega^*$
uniformly from $\mathcal Z_\infty$ to $\mathcal X_\infty$, and (H4) gives a uniform bound for
$D_\omega^{-1}$. By Lemma~\ref{rd:lem-extension} (iv),
$\Gamma_\omega^*(\Gamma_\omega w) = d_\omega w + K_\omega w$ for every
$w \in L^2_{\mathrm{uloc}}$; hence $I - B_\omega\Gamma_\omega = -D_\omega^{-1}K_\omega$
on $\mathcal X_\infty$.

Fix a bounded cube $Q$ and a sequence $(\omega_n, v_n)$ with
$\|v_n\|_{(\infty)} \le 1$. Since $\mathcal W$ is compact, we may assume that
$\omega_n \to \omega$. Lemma~\ref{lem:tails} allows us to truncate the input to a larger
cube with an error tending uniformly to zero. For each fixed truncation, (H5)
gives relative compactness of the restricted outputs of $K_{\omega_n}$; thus
the same relative compactness holds before truncation
(Lemma~\ref{lem:cc-tools}\textup{(d)}). By (H4), the inverse
multipliers $d_{\omega_n}^{-1}$ converge locally in measure; their uniform
bound gives strong convergence on $L^2(Q;\mathbb{C}^N)$. Applying them to a
convergent subsequence of the restricted outputs preserves convergence
(Lemma~\ref{lem:cc-tools}\textup{(c)}).
This proves (R).
\end{proof}
\begin{lemma}
\label{rd:lem-profile}
\label{lem:profile}
Assume \textup{(H2), (H3), and (R)}. Let $\omega_n\to\omega_\infty$ in $\mathcal W$, and let $w_n\in
L^2_{\mathrm{uloc}}$ satisfy
\begin{equation}
\label{rd:normalized}
\|w_n\|_{(\infty)}=\|w_n\|_{L^2(Q_0)}=1,\quad\|\Gamma_{\omega_n}w_n\|_{(\infty)}\to0 .
\end{equation}
Then a subsequence of $(w_n)$ converges strongly in $L^2_{\mathrm{loc}}$ to some $w_\infty\in
L^2_{\mathrm{uloc}}$ with $\Gamma_{\omega_\infty}w_\infty=0$ and
$\|w_\infty\|_{L^2(Q_0)}=\|w_\infty\|_{(\infty)}=1$.
\end{lemma}

\begin{proof}
By (R),
\[
 w_n=B_{\omega_n}\Gamma_{\omega_n}w_n+S_{\omega_n}w_n .
\]
The first term tends to zero in $\mathcal X_\infty$.
By (R) the second has a strongly convergent subsequence on every
bounded cube. A diagonal extraction gives a single subsequence
with $w_n\to w_\infty$ strongly in $L^2_{\rm loc}$.
Cell norms pass to the limit: $\|\mathbf 1_kw_\infty\|=\lim_n\|\mathbf 1_kw_n\|\le1$ for
every $k$, with equality for $k=0$. Hence $w_\infty\in L^2_{\mathrm{uloc}}$ and
$\|w_\infty\|_{(\infty)}=\|w_\infty\|_{L^2(Q_0)}=1$.

Fix a cell $k$ and $P_L=\mathbf1_{[-L,L]^d}$. By (H3) and
Lemma~\ref{rd:lem-extension} (ii),
\[
 \|\mathbf1_k\Gamma_\omega(1-P_L)u\|
 \le e_k(L)\|u\|_{(\infty)},\qquad
 e_k(L):=\tau_\kappa(L-|k|-1)\longrightarrow0 .
\]
For fixed $L$, $P_Lw_n\to P_Lw_\infty$ in $\mathcal X$ and the
operator norms are uniformly bounded by $\|\kappa\|_1$. Applying
(H2) to the fixed input $P_Lw_\infty$ and taking square roots gives
\[
 \|\mathbf1_k\Gamma_{\omega_\infty}P_Lw_\infty\|
 =\lim_n\|\mathbf1_k\Gamma_{\omega_n}P_Lw_n\|\le e_k(L),
\]
where the inequality follows from \eqref{rd:normalized} and the tail
bound. Applying the same bound to $w_\infty$ gives
$\|\mathbf1_k\Gamma_{\omega_\infty}w_\infty\|\le 2e_k(L)$.
Let $L\to\infty$ and then vary $k$ to obtain
$\Gamma_{\omega_\infty}w_\infty=0$.

\end{proof}

\begin{proof}[Proof of Theorem~\ref{rev:hull}]
\textbf{Observability implies injectivity and uniform lower bounds.}
This is the only part of the proof that uses \eqref{abs:full-hull}.
Assume the first assertion, with constant $c>0$.
Covariance and the unitarity of translations on $\mathcal X$ and
$\mathcal Z$ give the same lower bound on every translate of
$\omega_0$. Given $\omega\in\mathcal W$, choose $h_n$ such that
$\theta_{h_n}\omega_0\to\omega$. Fix a bounded-support input
$f\in\mathcal X$. By (H3), the output cell norms are dominated by
the fixed $\ell^2$ sequence $\kappa*(\|\mathbf1_jf\|)_j$.
By (H2) their squares converge cell by cell. Dominated convergence
for the sum of these squares gives convergence of the global
observation norms, so
\[
  \|\Gamma_\omega f\|_{\mathcal Z}\ge c\|f\|_{\mathcal X}.
\]
By density of bounded-support inputs and boundedness of
$\Gamma_\omega$, this holds for every $f\in\mathcal X$.
Consequently,
\[
  \inf_{\omega\in\mathcal W}\nu_2(\Gamma_\omega)\ge c>0.
\]
Proposition~\ref{rd:prop-rung}, applied to the whole family, gives
\[
  c_*:=\inf_{\omega\in\mathcal W,\;1\le p\le\infty}
       \nu_p(\Gamma_\omega)>0 .
\]
The case $p=\infty$ excludes every nonzero uniformly local kernel
element and proves the second assertion.

\textbf{Failure of observability produces a nonzero invisible profile.}
Suppose the first assertion fails, so that
$\nu_2(\Gamma_{\omega_0})=0$. Apply
Corollary~\ref{rd:cor-ladder} to the singleton family
$\{\Gamma_{\omega_0}\}$, from $p=2$ to $q=\infty$.
For any sequence $\varepsilon_n\downarrow0$, it gives compactly
supported $f_n\ne0$ satisfying
\[
  \|\Gamma_{\omega_0}f_n\|_{(\infty)}
  \le\varepsilon_n\|f_n\|_{(\infty)}.
\]
Normalize $\|f_n\|_{(\infty)}=1$. Since $f_n$ meets finitely many
cells, there is $k_n\in\mathbb Z^d$ such that
$\|f_n\|_{L^2(Q_{k_n})}=1$. Set
\[
  w_n=\tau_{k_n}f_n,
  \quad\omega_n=\theta_{k_n}\omega_0.
\]
Integer translations preserve the amalgam norms. Covariance thus gives
\[
\|w_n\|_{(\infty)}=\|w_n\|_{L^2(Q_0)}=1,\quad
\|\Gamma_{\omega_n}w_n\|_{(\infty)}\le\varepsilon_n.
\]
By compactness, after passing to a subsequence,
$\omega_n\to\omega_\infty\in\mathcal W$.
Lemma~\ref{rd:lem-profile} yields
\[
\Gamma_{\omega_\infty}w_\infty=0,
 \quad
\|w_\infty\|_{L^2(Q_0)}=\|w_\infty\|_{(\infty)}=1.
\]
This proves the asserted profile statement and contradicts the
second assertion whenever injectivity holds throughout the hull.
\end{proof}

\begin{proof}[Proof of Corollary~\ref{abs:compact-family}]
If the family has a uniform $L^2$ lower bound,
Proposition~\ref{rd:prop-rung} gives a common lower bound for all
environments and all $p\in[1,\infty]$. In particular, its $p=\infty$ lower bound excludes
nonzero uniformly local kernel elements.

Conversely, suppose
$\inf_{\omega\in\mathcal W}\nu_2(\Gamma_\omega)=0$.
Applying Corollary~\ref{rd:cor-ladder} to the entire family, from
$p=2$ to $q=\infty$, gives parameters $\varpi_n\in\mathcal W$
and compactly supported $f_n\ne0$ such that
\[
  \|f_n\|_{(\infty)}=1,
  \quad
  \|\Gamma_{\varpi_n}f_n\|_{(\infty)}\le\varepsilon_n,
  \quad\varepsilon_n\downarrow0.
\]
Choose maximizing cells $k_n$ and set
\[
w_n=\tau_{k_n}f_n,
\quad\omega_n=\theta_{k_n}\varpi_n.
\]
As above, these satisfy \eqref{rd:normalized}.
Compactness gives a convergent subsequence of $(\omega_n)$, and
Lemma~\ref{rd:lem-profile} produces a nonzero uniformly local kernel
element for its limit. This contradicts injectivity on every
configuration. The argument uses compactness and invariance of
$\mathcal W$, without requiring a dense orbit.
\end{proof}

\begin{proof}[Proof of Proposition~\ref{prop:obs-vs-gram}]
If $\Gamma_\omega w=0$, then $H_\omega w=\Gamma_\omega^*(\Gamma_\omega w)=0$ by
Lemma~\ref{rd:lem-extension}\textup{(iii)}. This proves the inclusion and
\textup{(b)}$\Rightarrow$\textup{(a)}.
Assume \textup{(a)}. Corollary~\ref{abs:compact-family} yields $c>0$ such that
$\|\Gamma_\omega f\|_{\mathcal Z}\ge c\|f\|_{\mathcal X}$ for all
$\omega\in\mathcal W$ and $f\in\mathcal X$. Since
\[
c^2\|f\|_{\mathcal X}^2\le\|\Gamma_\omega f\|_{\mathcal Z}^2
=\langle H_\omega f,f\rangle_{\mathcal X}
\le\|H_\omega f\|_{\mathcal X}\|f\|_{\mathcal X} ,
\]
we get $\nu_2(H_\omega)\ge c^2$ for all $\omega$. By
Lemma~\ref{rd:lem-extension}\textup{(iii)}, the operators $H_\omega$ share
the summable profile $\widetilde\kappa*\kappa$, so
Proposition~\ref{rd:prop-rung}, with $\mathcal E=\mathcal F=\mathcal X$,
gives $\inf_{\omega\in\mathcal W}\nu_\infty(H_\omega)>0$. This excludes
nonzero elements of $\ker_{L^2_{\mathrm{uloc}}}H_\omega$.
\end{proof}

\begin{proposition}[Necessity of (R)]\label{new:prop-R-necessary}
Assume (H3), and suppose that
\begin{equation}\label{new:eq-uniform-lower}
  c:=\inf_{\omega\in\mathcal W}\ \inf_{f\in \mathcal X\setminus\{0\}}
  \frac{\|\Gamma_\omega f\|_{\mathcal Z}}{\|f\|_{\mathcal X}}>0 .
\end{equation}
Then each Gramian $H_\omega=\Gamma_\omega^*\Gamma_\omega$, extended to
$\mathcal X_\infty$ by Lemma~\ref{rd:lem-extension}, is invertible on $\mathcal X_\infty$, and
$\sup_{\omega}\|H_\omega^{-1}\|_{\mathcal L(\mathcal X_\infty)}<\infty$. The operators
\[
  B_\omega:=H_\omega^{-1}\Gamma_\omega^*\colon \mathcal Z_\infty\to \mathcal X_\infty
\]
are uniformly bounded and satisfy $B_\omega\Gamma_\omega=I$ on $\mathcal X_\infty$.
In particular, (R) holds with $S_\omega=0$. If (H1)--(H2) and \eqref{abs:full-hull} hold,
the same conclusion follows as soon as the observability inequality holds
for $\omega_0$.
\end{proposition}

\begin{proof}
By (H3) and Lemma~\ref{rd:lem-extension} \textup{(iii)}, the adjoints $\Gamma_\omega^*$ are
$\widetilde\kappa$-localized, and the Gramians $H_\omega$ share the summable
profile $\widetilde\kappa*\kappa$. By \eqref{new:eq-uniform-lower},
\[
  c^2\|f\|_{\mathcal X}^2\le\langle H_\omega f,f\rangle\le\|H_\omega f\|_{\mathcal X}\|f\|_{\mathcal X} ,
\]
so $\nu_2(H_\omega)\ge c^2$ for every $\omega$. Proposition~\ref{rd:prop-rung} then gives
\[
  c_1:=\inf_{\omega\in\mathcal W,\ 1\le p\le\infty}\nu_p(H_\omega)>0 .
\]

Fix $\omega$ and write $H=H_\omega$. The cellwise identification
$\mathcal X_1=\ell^1(\{L^2(Q_k;\mathbb C^N)\}_{k\in\mathbb Z^d})$
gives $\mathcal X_1^*\cong\mathcal X_\infty$ under the summed $L^2$ pairing;
see also \cite{FournierStewart1985} for amalgam-space duality.
Let $f\in\mathcal X_1$ and $\psi\in\mathcal X_\infty$. The
terms of the double series
$\sum_{j,k}\langle\mathbf 1_kH\mathbf 1_jf,\mathbf 1_k\psi\rangle$ are bounded
by $(\widetilde\kappa*\kappa)(k-j)\,\|\mathbf 1_jf\|\,\|\mathbf 1_k\psi\|$, so the
series converges absolutely. Summing it in both orders and using $H^*=H$ gives
\begin{equation}\label{new:eq-duality}
  \langle Hf,\psi\rangle=\langle f,H\psi\rangle
  \qquad(f\in \mathcal X_1,\ \psi\in \mathcal X_\infty),
\end{equation}
where $H\psi$ denotes the uniformly local extension. On $\mathcal X_1$, the operator $H$
is bounded below by $c_1$, so its range is closed. If $\psi\in \mathcal X_\infty\cong
\mathcal X_1^*$ annihilates $H(\mathcal X_1)$, then \eqref{new:eq-duality} gives $H\psi=0$, and
$\nu_\infty(H)\ge c_1$ forces $\psi=0$. Hence $H$ is invertible on $\mathcal X_1$,
with $\|H^{-1}\|_{\mathcal L(\mathcal X_1)}\le c_1^{-1}$. By \eqref{new:eq-duality}, the
extension of $H$ to $\mathcal X_\infty$ is the Banach adjoint of $H|_{\mathcal X_1}$. It is
therefore invertible on $\mathcal X_\infty$, with the same bound.

Consequently
$\|B_\omega\|_{\mathcal L(\mathcal Z_\infty,\mathcal X_\infty)}\le c_1^{-1}\|\kappa\|_{1}$. For
$w\in \mathcal X_\infty$, Lemma~\ref{rd:lem-extension} \textup{(iii)} gives
$\Gamma_\omega^*(\Gamma_\omega w)=H_\omega w$, and hence $B_\omega\Gamma_\omega w=w$.

Under (H1)--(H2) and \eqref{abs:full-hull}, the first part of the proof of Theorem~\ref{rev:hull}
transfers the lower bound from $\omega_0$ to every $\omega\in\mathcal W$, with
the same constant.
\end{proof}

\begin{corollary}[Characterization of observability]\label{new:cor-characterization}
Assume (H1)--(H3) and \eqref{abs:full-hull}. The observability inequality holds for
$\omega_0$ if and only if
\begin{enumerate}
  \item[(R)] the uniform reconstruction property holds, and
  \item[(U)] $\ker_{L^2_{\mathrm{uloc}}}\Gamma_\omega=\{0\}$ for every
             $\omega\in\mathcal W$.
\end{enumerate}
Under (H1)--(H3), neither condition implies the other. The Gaussian observer of
Remark~\ref{abs:remark-multiplier} satisfies (U) but not (R). The two-speed example of Section~\ref{abs:assumptions}
satisfies (H4)--(H5), hence (R), but not (U).
\end{corollary}

\begin{proof}
Sufficiency is Theorem~\ref{rev:hull}. For necessity,
Proposition~\ref{new:prop-R-necessary} first gives (R); Theorem~\ref{rev:hull}
then gives (U). The Gaussian observer is
injective on $L^2_{\rm uloc}$ but not observable, so by Theorem~\ref{rev:hull} it cannot
satisfy (R).
\end{proof}

\begin{remark}[Local versus global]\label{new:rem-local-global}
Proposition~\ref{new:prop-R-necessary} uses the observability inequality to
construct an exact left inverse. Condition (R) asks for much less: a uniformly
bounded reconstruction modulo errors that are collectively compact on bounded
windows. Once uniformly bounded reconstruction operators have been constructed,
the compactness of the errors can be checked window by window,
for instance through (H4)--(H5). Theorem~\ref{rev:hull} converts this local information, together with
uniqueness on the hull, into the global estimate.
\end{remark}

\begin{lemma}[Output transformations]\label{new:lem-output}
Assume that $\{\Gamma_\omega\}_{\omega\in\mathcal W}$ satisfies (H1), (H3) and
(R). Let $\mathcal Z'=L^2((0,T)\times\mathbb R^d;\mathbb C^{m'})$, and let
$M_\omega\in\mathcal L(\mathcal Z,\mathcal Z')$, $\omega\in\mathcal W$, satisfy:
\begin{enumerate}
  \item[(a)] covariance: $M_{\theta_h\omega}\tau_h^{\mathcal Z}=\tau_h^{\mathcal Z'}M_\omega$ for
             $h\in\mathbb R^d$;
  \item[(b)] localization: $\|\mathbf 1_kM_\omega\mathbf 1_j\|\le\nu(k-j)$ for
             some $\nu\in\ell^1(\mathbb Z^d)$ independent of $\omega$;
  \item[(c)] uniform left invertibility: there are operators
             $L_\omega\in\mathcal L(\mathcal Z'_\infty,\mathcal Z_\infty)$ with
             $\sup_\omega\|L_\omega\|<\infty$ and $L_\omega M_\omega z=z$ for
             all $z\in \mathcal Z_\infty$.
\end{enumerate}
Then $\Gamma'_\omega:=M_\omega\Gamma_\omega$ satisfies (H1), (H3) with profile
$\nu*\kappa$, and (R) with
\[
  B'_\omega=B_\omega L_\omega,\qquad S'_\omega=S_\omega .
\]
Moreover, $\ker_{L^2_{\rm uloc}}\Gamma'_\omega=\ker_{L^2_{\rm uloc}}\Gamma_\omega$
for every $\omega\in\mathcal W$.
\end{lemma}

\begin{proof}
Covariance follows from (a) and (H1). By Lemma~\ref{rd:lem-extension}\textup{(iii)}, with the obvious
change of output space, $\Gamma'_\omega$ is $(\nu*\kappa)$-localized, and
$\Gamma'_\omega w=M_\omega(\Gamma_\omega w)$ for $w\in \mathcal X_\infty$. Hence
\[
  B'_\omega\Gamma'_\omega w=B_\omega L_\omega M_\omega\Gamma_\omega w
  =B_\omega\Gamma_\omega w ,
\]
so $I-B'_\omega\Gamma'_\omega=S_\omega$ on $\mathcal X_\infty$, and
$\sup_\omega\|B'_\omega\|\le\sup_\omega\|B_\omega\|\,\sup_\omega\|L_\omega\|$.
If $\Gamma'_\omega w=0$, then $\Gamma_\omega w=L_\omega M_\omega\Gamma_\omega w=0$.
The converse inclusion is immediate.
\end{proof}

\begin{remark}[Continuity of the cell energies]\label{new:rem-H2}
Assumption (H2) does not pass automatically from $\Gamma_\omega$ to
$M_\omega\Gamma_\omega$, because it controls only cell energies. Two sufficient
conditions cover the examples of this paper.
\begin{enumerate}
  \item[(i)] $M_\omega$ is multiplication by a matrix field $m_\omega(t,x)$ with
        $\sup_\omega\|m_\omega\|_{L^\infty}<\infty$ and $m_{\omega_n}\to m_\omega$
        locally in measure, and $\Gamma_{\omega_n}f\to\Gamma_\omega f$ in
        $L^2_{\rm loc}$ for every bounded-support $f$.
  \item[(ii)] The output splits into channels,
        $\Gamma_\omega=(\Gamma^{s}_\omega,\Gamma^{e}_\omega)$, where $\Gamma^{s}$
        has the strong convergence in (i), and
        $M_\omega=\mathrm{diag}(m_\omega,I)$ with $m_\omega$ as in (i).
        The cell energy of $\Gamma^{e}$ then converges, since it is the
        difference between the total energy, controlled by (H2), and the
        energy of $\Gamma^{s}$.
\end{enumerate}
\end{remark}

\begin{remark}[(R) versus (H4)--(H5)]\label{new:rem-R-vs-H45}
The Gramian of $M_\omega\Gamma_\omega$ is
$\Gamma_\omega^*M_\omega^*M_\omega\Gamma_\omega$. The decomposition required by
(H4)--(H5) need not be preserved by the output transformations of
Lemma~\ref{new:lem-output}, whereas (R) is preserved. Both examples of this
paper that separate (R) from (H4)--(H5) have this form:
\begin{itemize}
  \item In the singleton example of Section~\ref{abs:assumptions}, $\Gamma=A\Gamma^0$ with
        $\Gamma^0f(t,x)=f(x)$, whose Gramian is $T\,I$, and $L=A^{-1}$ is given
        by the Neumann series.
  \item In Section~\ref{new:sec-heat-flux}, the measured flux is $M\Gamma^0$,
        where $\Gamma^0$ is the energy-normalized observation and
        $M=\mathrm{diag}(-\mathsf A^{1/2},1)$.
\end{itemize}
\end{remark}
\begin{remark}[Several velocities in higher dimension]\label{new:rem-multi-d}
Let $d\ge2$, $V_1\ne V_2$, and consider the whole-space observation
$y=q_1(x+tV_1)+q_2(x+tV_2)$. Its Fourier Gramian has eigenvalues
$T\pm|g(\vartheta)|$, where $g(\vartheta)=\int_0^Te^{is\vartheta}\,ds$ and
$\vartheta=(V_1-V_2)\cdot\xi$. Its smallest eigenvalue is therefore comparable
to $\min\{1,\vartheta^2\}$. The Kalman-type degeneracy thus occurs on the whole
hyperplane $\{(V_1-V_2)\cdot\xi=0\}$, including arbitrarily high frequencies,
and Remark~\ref{rd:rem-velocities} shows that (H5) fails. By
Proposition~\ref{new:prop-R-necessary}, (R) is necessary for observability in any normalized
space in which (H1)--(H3) hold. The open problem is to verify it locally,
without assuming observability.
\end{remark}
\section{Bounded-gap observation sets and compact weight families}
\label{kal:section}
We record the geometric and compactness ingredients used in the localized observability argument. These depend only on the characteristic speeds and the observation geometry.
Let $\Lambda\subset\mathbb R\setminus\{0\}$ be a finite nonempty set, and define
\[
\lambda_\bullet=\min_{\lambda\in\Lambda}|\lambda|,\quad\lambda^\bullet=\max_{\lambda\in\Lambda}|\lambda|.
\]
For $\ell>0$ and $G\geq0$, let $\mathcal G_{\ell,G}$ be the class
of open sets
\begin{equation}\label{kal:bands}
 \Omega=\bigcup_{n\in\mathbb Z}(a_n,b_n),
 \quad a_n<b_n\leq a_{n+1}, \quad b_n-a_n\geq\ell,\quad a_{n+1}-b_n\leq G,
\end{equation}
with endpoints tending to $\pm\infty$ as appropriate.
We also include $\Omega=\mathbb R$ in this class.
For $T>0$ and a bounded measurable weight $a:\mathbb R\to[0,1]$,
set
\begin{equation}\label{kal:occupation}
 m_{\lambda,a}(z)
 :=\int_0^T a(z-\lambda s)\,ds.
\end{equation}

\begin{lemma}[Geometric intersection bound]
\label{rev:intersection}
Let $O\in\mathcal G_{L_0,G_0}$, where $L_0>0$ and $G_0\geq0$.
Every interval $J$ of length $L>G_0$ satisfies
\begin{equation}\label{kal:elementarygeometry}
 |J\cap O|\geq\min\{L_0,L-G_0\}.
\end{equation}
\end{lemma}

\begin{proof}
The assertion is immediate if $O=\mathbb R$.
If $J$ contains a whole observed interval, then
$|J\cap O|\geq L_0$. Otherwise, it meets at most two observed intervals, and its
unobserved part is contained, up to endpoints, in a single gap. Hence $|J\setminus O|\leq G_0$, giving
$|J\cap O|\geq L-G_0$.
\end{proof}

\begin{lemma}
\label{kal:compactweights}
Assume $T>G/\lambda_\bullet$
and choose
\[
 0<\delta<\min\left\{\frac{\ell}{4},
\frac{\lambda_\bullet T-G}{4}\right\}.
\]
For $\Omega\in\mathcal G_{\ell,G}$, define
\begin{equation}\label{kal:minorant}
 a_\Omega(x)=
 \begin{cases}
 \displaystyle
 \min\left\{1,\delta^{-1}
\operatorname{dist}(x,\mathbb R\setminus\Omega)\right\},
& \Omega\neq\mathbb R,\\
 1, & \Omega=\mathbb R.
 \end{cases}
\end{equation}
Then
 $0\leq a_\Omega\leq\mathbf1_\Omega,\quad
 \operatorname{Lip}(a_\Omega)\leq\delta^{-1},$
and, for every $\lambda\in\Lambda$ and $z\in\mathbb R$,
\begin{equation}
 m_{\lambda,a_\Omega}(z)
 \geq
 \frac{1}{|\lambda|}
 \min\{\ell-2\delta,\ |\lambda|T-G-2\delta\}
 \geq m_0,\quad \text{where} \quad
 m_0:=
 \frac{1}{\lambda^\bullet}
 \min\{\ell-2\delta,\ \lambda_\bullet T-G-2\delta\}>0.
\end{equation}
Moreover, the family
\begin{equation}\label{kal:compactfamily}
 \mathfrak A=
 \left\{
 a:\mathbb R\to[0,1]:
 \begin{array}{l}
  \operatorname{Lip}(a)\leq\delta^{-1},\\
  m_{\lambda,a}(z)\geq m_0
  \quad\text{for all }\lambda\in\Lambda,\ z\in\mathbb R
 \end{array}
 \right\}
\end{equation}
is compact for locally uniform convergence, invariant under spatial translations, and contains every $a_\Omega$.
\end{lemma}
\begin{proof}
The range and Lipschitz bounds follow from the definition. For $\Omega\neq\mathbb R$, the weight $a_\Omega$ equals one on
\[
 O_\delta:=\bigcup_{n\in\mathbb Z}(a_n+\delta,b_n-\delta).
\]
The intervals of $O_\delta$ have lengths at least $\ell-2\delta$, and their gaps have lengths at most $G+2\delta$.
The characteristic interval
$J_{\lambda,z}:=\{z-\lambda s:0<s<T\}$
has length $|\lambda|T>G+2\delta$.
Therefore, by \eqref{kal:elementarygeometry},
\[
 \begin{aligned}
 m_{\lambda,a_\Omega}(z)
 \geq\frac{|J_{\lambda,z}\cap O_\delta|}{|\lambda|}\geq\frac{1}{|\lambda|}
 \min\{\ell-2\delta,\ |\lambda|T-G-2\delta\}\geq m_0.
 \end{aligned}
\]
For $\Omega=\mathbb R$, the occupation equals $T$, and the same bounds hold.
The functions in $\mathfrak A$ are uniformly bounded and share a common Lipschitz constant. Arzel\`a--Ascoli and a diagonal argument give relative compactness for locally uniform convergence.
The family is closed: the range and Lipschitz bounds pass to the limit, as do the occupation bounds, since $a_j\to a$ locally uniformly implies
\[
 m_{\lambda,a_j}\longrightarrow m_{\lambda,a}
 \quad\text{locally uniformly in }z.
\]
Thus $\mathfrak A$ is compact.
Finally, for $a_h(x):=a(x+h)$,
$m_{\lambda,a_h}(z)=m_{\lambda,a}(z+h),$
which proves translation invariance.
\end{proof}
\section{Why the translation hull is necessary}
\label{app:section}

This section shows, on two explicit examples, that uniqueness for the original configuration
does not suffice.
\subsection{The Gramian of a variable observer}
\label{app:gramian}
Let $\mathbb V=\operatorname{diag}(\lambda_1I_{d_1},\ldots,
\lambda_pI_{d_p})$, with distinct real $\lambda_j$ and
$n=\sum_jd_j$. For $f=(f_1,\ldots,f_p)$, set
$(U(t)f)_j(x)=f_j(x-\lambda_jt)$.
Given $B^\top=(C_1,\ldots,C_p)\in
C_b^\infty(\mathbb R;\mathbb C^{m\times n})$ and $T>0$,
define the bounded observation operator
\[
 \mathcal O_Bf(t,x)=B^\top(x)U(t)f(x).
\]
Its translation hull
$\mathscr H(B^\top )=
\overline{\{B^\top(\cdot+h):h\in\mathbb R\}}^{\,C^1_{\mathrm{loc}}}$
is compact and retains the uniform derivative bounds.

The diagonal characteristic contribution to
$G_B=\mathcal O_B^*\mathcal O_B$ is the multiplication
operator $D_B=\operatorname{diag}(D_1,\ldots,D_p)$, where
\[
 D_j(z)=\int_0^T
 C_j(z+\lambda_jt)^* C_j(z+\lambda_jt)\,dt.
\]
Assume uniform characteristic visibility:
\begin{equation}
 D_j(z)\ge dI_{d_j}
 \quad(z\in\mathbb R,\ 1\le j\le p),
 \quad d>0.
 \label{app:diagonal-visibility}
\end{equation}
This condition is preserved throughout $\mathscr H(B^\top)$.

\subsection{The same rank and visibility tests, different limits at infinity}
\label{app:examples}

From now on, take
\[
 \mathbb V=\operatorname{diag}(I_2,2I_2),\quad
 \mathsf J=\begin{pmatrix}0&-1\\1&0\end{pmatrix},\quad
 R(x)=e^{x\mathsf J},\quad
 C_r(x)=\bigl(I_2,-r(x)R(x)\bigr).
\]
Let $0\le\chi\le1$ be smooth and $2\pi$-periodic, with $\omega_\chi=\{\chi>0\}\neq\emptyset$ a finite union of open bands per period,
and let $G$ be the maximal gap length. The case $\chi\equiv1$ is allowed,
with $G=0$. We observe
\[
 B^\top(x)q(t,x)=\chi(x)\bigl(q_1(t,x)-r(x)R(x)q_2(t,x)\bigr),
 \quad 0<t<T.
\]
For $r\in C_b^\infty(\mathbb R)$ bounded above and below by positive constants, the
unlocalized pointwise observability matrix is
\[
 \begin{pmatrix}C_r(x)\\ C_r(x)\mathbb V\end{pmatrix}
 =\begin{pmatrix}I_2&-r(x)R(x)\\I_2&-2r(x)R(x)\end{pmatrix}.
\]
It is uniformly invertible: subtracting the first block row from the second
leaves the invertible block $-r(x)R(x)$.  For the localized observer, the
diagonal Gramian blocks are
\begin{align*}
 D_1(z)=\int_0^T\chi(z+t)^2\,dt\,I_2,\quad
 D_2(z)=\int_0^T\chi(z+2t)^2r(z+2t)^2\,dt\,I_2.
\end{align*}
When $T>G$, every characteristic interval meets $\omega_\chi$.
The integrals of $\chi^2$ are therefore positive continuous periodic
functions, with positive minima; the lower bound on $r$ gives
\eqref{app:diagonal-visibility}. If $\chi\equiv1$ and $r\ge c$, then
$D_B\ge T\operatorname{diag}(I_2,c^2I_2)$.

\begin{lemma}
\label{app:distribu}
Suppose $r>0$ is smooth and $h=r'/r$ is nonconstant on every nonempty open
interval.  If a distributional transport solution satisfies
$\chi(x)(q_1-rRq_2)=0$ on $(0,T)\times\mathbb R$, with $T>G$, then its
initial datum is zero.  No integrability assumption is needed.
\end{lemma}

\begin{proof}
On each component $I$ of $\omega_\chi$, $q_1=rRq_2$.
Apply $\partial_t+\partial_x$ and use
$(\partial_t+2\partial_x)q_2=0$ to obtain
$\partial_xq_2=(hI_2+\mathsf J)q_2$ in distributions.
Thus $\partial_x(r^{-1}R(-x)q_2)=0$, so
$q_2(t,x)=r(x)R(x)v_I(t)$ for a vector distribution $v_I$.
Its transport equation and its $x$ derivative give
\[
 v_I'+2\mathsf Jv_I+2h(x)v_I=0,\quad h'(x)v_I=0.
\]
Since $h'$ is not identically zero on $I$, testing in $x$ gives
$v_I=0$ and hence $q_1=q_2=0$ on $(0,T)\times \omega_\chi$.
Write $q_j(t,x)=f_j(x-\lambda_jt)$, with $\lambda_1=1$, $\lambda_2=2$.
Each interval $z+\lambda_j(0,T)$ has length greater than $G$ and meets $\omega_\chi$;
Lemma~\ref{aux:pullback}, applied with speed $-\lambda_j$, therefore
implies $f_j=0$.
\end{proof}

\begin{theorem}
\label{app:periodic-observability}
Let $r(x)=2+\sin x$.  For the observer $\chi C_r$ above and $T>0$,
there is a constant $C_T<\infty$ such that
\[
 \|q(0)\|_2^2\le C_T
 \int_0^T\!\int_{\mathbb R}\chi(x)^2
|q_1(t,x)-(2+\sin x)R(x)q_2(t,x)|^2\,dx\,dt
\]
if and only if $T>G$. In particular, with $\chi\equiv1$, observability holds for every $T>0$.
\end{theorem}

\begin{proof}
Suppose first that $T>G$. Visibility was verified above. By periodicity every hull element is $ \chi(x+\theta)
 \bigl(I_2,-(2+\sin(x+\theta))R(x+\theta)\bigr).$
The logarithmic derivative $\cos x/(2+\sin x)$ is real analytic and
nonconstant, and so are all its translates. Since $R'=\mathsf JR$ is preserved
by translation and the gaps retain the bound $G$,
Lemma~\ref{app:distribu} gives distributional uniqueness on every hull
element.
Proposition~\ref{abs:concrete} applies to this hull with
characteristic speeds $-1$ and $-2$, zero integral correction,
and the coercivity bound \eqref{app:diagonal-visibility}.
Together with uniqueness on every hull element,
Theorem~\ref{rev:hull} yields the observability estimate.

If $0<T\le G$, let $z$ be the left endpoint of a gap of maximal
length. For initial data $(f_1,0)$ the observed energy is
$\int m(y)|f_1(y)|^2\,dy$, where $m(y)=\int_0^T\chi(y+t)^2\,dt$.
Since $m$ is continuous and $m(z)=0$, unit packets supported in
shrinking neighborhoods of $z$ have observation tending to zero.
\end{proof}

\begin{theorem}
\label{app:tanh-nonobservability}
Let $r(x)=2+\tanh x$.  The observer $\chi C_r$ has distributional uniqueness
for every $T>G$, and satisfies the same uniform algebraic and characteristic
visibility conditions as above.  Nevertheless, for every fixed $T>0$ there
are smooth compactly supported initial data $f_n$ with $\|f_n\|_2=1$ and
$\|\mathcal O_{\chi C_r}f_n\|_2\to0$.
\end{theorem}

\begin{proof}
The logarithmic derivative $h(x)=\operatorname{sech}^2x/(2+\tanh x)$
is real analytic and nonconstant, so Lemma~\ref{app:distribu} gives
uniqueness for $T>G$. The bounds $1<r<3$ give the algebraic and
characteristic visibility for $T>G$ already established.

For the failure of observability, choose a real
$\phi\in C_c^\infty(-1,1)$ with $\|\phi\|_2=1$, a unit vector
$v\in\mathbb R^2$, and set
\begin{align*}
 a_n(x)=n^{-1/2}\phi\!\left(\frac{x-n^2}{n}\right),\quad
 \alpha=\frac3{\sqrt{10}},\quad
 f_{1,n}(x)=\alpha a_n(x)R(2x)v,\quad
 f_{2,n}(x)=\frac\alpha3a_n(x)R(x)v.
\end{align*}
Orthogonality of $R$ gives
$\|f_n\|_2^2=\alpha^2(1+1/9)=1$.  The two transports satisfy
\[
 q_{1,n}-rRq_{2,n}
 =\alpha R(2x-2t)v
\left[a_n(x-t)-\frac{r(x)}3a_n(x-2t)\right].
\]
Now $\|a_n'\|_2=n^{-1}\|\phi'\|_2$, and
\[
 \|a_n(\,\cdot-t)-a_n(\,\cdot-2t)\|_2
 \le t\|a_n'\|_2.
\]
Whenever $a_n(x-2t)\ne0$ and $0<t<T$, one has $x\ge n^2-n$.  Therefore,
with
$\varepsilon_n=\sup_{x\ge n^2-n}|1-r(x)/3|\to0$ and $0\le\chi\le1$,
\[
 \|\mathcal O_{\chi C_r}f_n\|_2^2
 \le2\alpha^2\left(
       \frac{T^3}{3n^2}\|\phi'\|_2^2+T\varepsilon_n^2\right)
 \longrightarrow0.
\]
This contradicts any uniform observability inequality.
\end{proof}

\begin{remark}[What the translated limits detect]
\label{app:interpretation}
Along sequences $h_\nu\to\pm\infty$ with
$h_\nu\bmod2\pi\to\theta$, the observers in
Theorem~\ref{app:tanh-nonobservability} converge to
\[
 C_{\pm,\theta}(x)=\chi(x+\theta)
\bigl(I_2,-c_\pm R(x+\theta)\bigr),\quad c_+=3,\quad c_-=1.
\]
These limits admit the nonzero bounded invisible solutions
\[
 q_1(t,x)=R(2x-2t+\theta)v,\quad
 q_2(t,x)=c_\pm^{-1}R(x-2t)v.
\]
These are $L^2_{\mathrm{uloc}}$ profiles, not $L^2$ states.
The normalized packets approximate the limiting profile at positive infinity and have
frequencies concentrating near $\{\pm1,\pm2\}$: this obstruction is
not the zero-frequency degeneracy of a constant Kalman observer.
Thus distributional uniqueness for the original medium does not replace
uniqueness on the full hull. The two examples share exactly the same
localization $\chi$ and speeds; only the second has an invisible
limiting configuration.
\end{remark}

\section{Localized observability with arbitrary real Jordan structure}
\label{jo:section}
\subsection{The intrinsic observation space and Jordan coordinates}
\label{glob:section}

To apply the abstract criterion to transport systems with Jordan blocks,
we first identify the norm measured by whole-line observation and its
dual control space. We recall the characterization from
\cite{Simpore2026}, with a short proof, and then pass to Jordan coordinates.
Zero characteristic speeds are allowed in the whole-line construction.

Assume that $A\in\mathbb R^{n\times n}$ has real spectrum and that
$B\in\mathbb R^{n\times m}$ satisfies
\[
 \operatorname{rank}(B,AB,\ldots,A^{n-1}B)=n.
\]
Set $E_{-1}=\{0\}$, $E_k=\operatorname{Ran}(B,\ldots,A^kB)$,
and $K=\min\{k:E_k=\mathbb R^n\}$. Let $\Pi_k^A$ be the orthogonal
projector onto $E_k\cap E_{k-1}^{\perp}$. For
$\lambda\in\Lambda:=\sigma(A)$, let $P_\lambda^A$ be the generalized
spectral projector of $A^\top$ and $m_\lambda$ the largest Jordan-block
size at $\lambda$. Define
\[
 A_{\lambda,j}
 =\frac1{j!}B^\top(A^\top-\lambda I)^jP_\lambda^A,
 \quad 0\le j<m_\lambda,
\]
and
\begin{equation}
\begin{gathered}
 W_K(\xi)=\sum_{k=0}^K|\xi|^{2k}\Pi_k^A,
 \quad
 W_J(\xi)=\sum_{\lambda\in\Lambda}
\sum_{j=0}^{m_\lambda-1}|\xi|^{2j}A_{\lambda,j}^*A_{\lambda,j},
 \\
 W(\xi)=
 \begin{cases}
  W_K(\xi),&0<|\xi|\le1,\\
  W_J(\xi),&|\xi|>1.
 \end{cases}
\end{gathered}
\label{glob:weight}
\end{equation}
The Kalman filtration describes the degeneracy at zero, whereas the
nilpotent parts determine the growth at high frequencies. Blocks with
the same eigenvalue contribute to a single output and are therefore
summed before taking its squared norm. The value of $W$ at zero is immaterial.

\begin{theorem}[Whole-line observation]
\label{glob:global-theorem}
For every $T>0$, the Gramian
\begin{equation}
 G_T(\xi)=\int_0^T e^{-is\xi A}BB^\top e^{is\xi A^\top}\,ds
 \label{glob:gramian}
\end{equation}
satisfies
\begin{equation}
 c_TW(\xi)\le G_T(\xi)\le C_TW(\xi),
 \quad \xi\ne0,
 \label{glob:matrix-equivalence}
\end{equation}
for some $c_T,C_T>0$. In particular, $W(\xi)$ is positive definite
for $\xi\ne0$.
\end{theorem}

\begin{proof}
The Kalman condition gives $G_T(\xi)>0$ for $\xi\ne0$: if
$B^\top e^{is\xi A^\top}z$ vanishes on $[0,T]$, its derivatives at
$s=0$ imply $z=0$.

Near zero, set $D_\xi=\sum_{k=0}^K(i\xi)^k\Pi_k^A$.
The identities $B^\top(A^\top)^j\Pi_k^A=0$ for $j<k$ give
\[
 B^\top e^{is\xi A^\top}D_\xi^{-1}
 =\sum_{k=0}^K\frac{s^k}{k!}B^\top(A^\top)^k\Pi_k^A
   +O_T(|\xi|),
\]
uniformly for $s\in[0,T]$. The limiting polynomial map is injective:
each coefficient must vanish, and $B^\top(A^\top)^k$ is injective on
$E_k\cap E_{k-1}^{\perp}$. Its Gramian is positive definite, so
$G_T(\xi)\asymp_T D_\xi^*D_\xi=W_K(\xi)$ near zero.

At high frequencies, the exact expansion is
\[
 B^\top e^{is\xi A^\top}z
 =\sum_{\lambda\in\Lambda}e^{is\xi\lambda}
   \sum_{j=0}^{m_\lambda-1}s^j(i\xi)^jA_{\lambda,j}z.
\]
The Gram matrix of the functions $s^je^{is\xi\lambda}$ has fixed
positive definite diagonal blocks
$\bigl(\int_0^T s^{j+\ell}\,ds\bigr)_{j,\ell}$.
Its off-diagonal blocks are $O_T(|\xi|^{-1})$ by integration by parts.
Applying this matrix estimate to the coefficients
$(i\xi)^jA_{\lambda,j}z$ yields $G_T(\xi)\asymp_T W_J(\xi)$ for
large $|\xi|$. The expansion and the Kalman condition also show that
$W_J(\xi)>0$ for every $\xi\ne0$. Positivity and compactness on the
remaining frequency annuli complete the proof.
\end{proof}

\begin{definition}
\label{glob:space}
The intrinsic observation space $\mathcal Y^-_{K,J}$ is the completion
of $\mathcal S(\mathbb R;\mathbb C^n)$ for
\begin{equation}
 \|f\|_{\mathcal Y^-_{K,J}}^2
 =\int_{\mathbb R}\widehat f(\xi)^*W(\xi)\widehat f(\xi)\,d\xi.
 \label{glob:norm}
\end{equation}
It is realized by measurable Fourier fields with
$W^{1/2}\widehat f\in L^2$, modulo equality almost everywhere.
Its dual $\mathcal Y^+_{K,J}$ has norm
\begin{equation}
 \|y\|_{\mathcal Y^+_{K,J}}^2
 =\int_{\mathbb R}\widehat y(\xi)^*W(\xi)^{-1}\widehat y(\xi)\,d\xi
 \label{glob:dual-norm}
\end{equation}
and pairing
$\langle y,f\rangle=\int_{\mathbb R}\widehat y(\xi)^*
\widehat f(\xi)\,d\xi$. Here $J$ records the Jordan structure;
we write $\mathcal Y^\pm_{K,J}(A,B)$ when the pair needs to be specified.
\end{definition}

Any fixed positive cutoff in \eqref{glob:weight} gives an equivalent
norm. Elements of $\mathcal Y^-_{K,J}$ need not be componentwise
$L^2$ functions, and no canonical distributional realization is assumed.
Evolution and observation will be defined by continuity from Schwartz data.

\paragraph{Observation and evolution}
For Schwartz data, put
$U(s)f=\mathcal F^{-1}(e^{is\xi A^\top}\widehat f)$.
Plancherel and \eqref{glob:matrix-equivalence} give
\[
 c_T\|f\|_{\mathcal Y^-_{K,J}}^2
 \le\int_0^T\|B^\top U(s)f\|_{L^2}^2\,ds
 \le C_T\|f\|_{\mathcal Y^-_{K,J}}^2.
\]
Applying the lower bound with $T=1$ to $U(t)f$ uses the observation
on $[t,t+1]\subset[-R,R+1]$ for $|t|\le R$. On the negative half of
this interval use the upper bound at $-\xi$ and the identity
$W(-\xi)=W(\xi)$. Together the two halves give
\[
 \sup_{|t|\le R}\|U(t)f\|_{\mathcal Y^-_{K,J}}^2
 \le \frac{2C_{R+1}}{c_1}\|f\|_{\mathcal Y^-_{K,J}}^2.
\]
Density gives a strongly continuous group on $\mathcal Y^-_{K,J}$.
Its transpose $S(t)=e^{-tA\partial_x}$ is a strongly continuous group
on $\mathcal Y^+_{K,J}$.

The whole-line upper bound also implies that every norm satisfying
$\|f\|_\sharp^2\lesssim
\int_0^T\|(B^\top U(s)f)|_\Omega\|_{L^2(\Omega)}^2\,ds$
on Schwartz data obeys
\begin{equation}
 \|f\|_\sharp\lesssim\|f\|_{\mathcal Y^-_{K,J}}.
 \label{rev:domination}
\end{equation}
Thus \eqref{glob:norm} is the strongest norm controlled by whole-line
observation, up to equivalence. For a positive definite Hermitian
multiplier $\Psi$, the corresponding whole-line lower observation bound
holds exactly when $\Psi(\xi)\lesssim W(\xi)$ almost everywhere.

\paragraph{Localized observation and control}
For every measurable $\Omega$, the map
$C_\Omega f=(B^\top f)|_\Omega$, initially defined on Schwartz data,
extends boundedly from $\mathcal Y^-_{K,J}$ into
$\mathcal U:=L^2(\Omega;\mathbb C^m)$. Indeed,
$B^\top=B^\top\Pi_0^A=\sum_\lambda A_{\lambda,0}$ implies
\[
 |B^\top z|^2\le c_B\,z^*W(\xi)z,
 \quad c_B=\max\{\|B^\top\|^2,\#\Lambda\}.
\]
Thus spatial restriction is applied to the $L^2$ output $B^\top f$.
Let $B_\Omega=C_\Omega'$ in the preceding duality and set
\[
 \Phi_tu=\int_0^tS(t-s)B_\Omega u(s)\,ds,
 \quad y(t)=S(t)y_0+\Phi_tu.
\]
Since $\Phi_T'f=C_\Omega U(T-\cdot)f$, Hilbert space duality
\cite{Lions1988,TucsnakWeiss2009} gives
$\operatorname{Ran}\Phi_T=\mathcal Y^+_{K,J}$ if and only if
\begin{equation}
 \int_0^T\|C_\Omega U(T-t)f\|_{\mathcal U}^2\,dt
 \ge c\|f\|_{\mathcal Y^-_{K,J}}^2,
 \qquad f\in\mathcal Y^-_{K,J},
 \label{glob:HUM}
\end{equation}
for some $c>0$. In that case, for $y_0,y_d\in\mathcal Y^+_{K,J}$,
the minimum-norm control steering $y_0$ to $y_d$ satisfies
\begin{equation}
 \|u_{\min}\|_{L^2(0,T;\mathcal U)}^2
 \le c^{-1}\|y_d-S(T)y_0\|_{\mathcal Y^+_{K,J}}^2.
 \label{glob:correct-cost}
\end{equation}

\paragraph{Jordan coordinates}
Choose a fixed real change of coordinates $P$ such that
\[
 \begin{aligned}
 F=P^{-1}A^\top P
   =\bigoplus_{\lambda\in\Lambda}(\lambda I+N_\lambda),
 C=B^\top P,\quad
 N_\lambda^{r_\lambda}=0,\quad
  r_\lambda=m_\lambda.
 \end{aligned}
\]
All Jordan blocks at the same speed are grouped in the coordinate
subspace $E_\lambda$. Its orthogonal projector is
$P_\lambda=P^{-1}P_\lambda^A P$; put $C_\lambda=C|_{E_\lambda}$
and extend maps on $E_\lambda$ by zero when needed.
The Kalman condition and its minimal index are unchanged:
\begin{equation}
 \bigcap_{j=0}^K\ker(CF^j)=\{0\}.
 \label{jo:kalman}
\end{equation}
On each $E_\lambda$, this implies
\begin{equation}
 \bigcap_{j=0}^{r_\lambda-1}\ker(C_\lambda N_\lambda^j)=\{0\},
 \label{jo:block-kalman}
\end{equation}
since a vector in this intersection is killed by every $CF^j$.

The orthogonal Kalman projectors are now computed in these coordinates.
Set
\[
 \mathcal N_0=\mathbb C^n,
 \quad \mathcal N_k=\bigcap_{j=0}^{k-1}\ker(CF^j)
 \quad(1\le k\le K+1),
\]
and let $\Pi_k$ project onto
$\mathcal N_k\cap\mathcal N_{k+1}^{\perp}$. Then
\begin{equation}
 CF^j\Pi_k=0\quad(j<k),
 \quad CF^k|_{\operatorname{Ran}\Pi_k}\text{ is injective}.
 \label{jo:low-identities}
\end{equation}
In general, $\Pi_k\ne P^{-1}\Pi_k^A P$, because a nonorthogonal change
of coordinates does not preserve orthogonal projections.
Thus $P_\lambda^A$ and $P_\lambda$ are spectral projectors, whereas
$\Pi_k^A$ and $\Pi_k$ are orthogonal projectors for the respective
Kalman filtrations.

The identity $B^\top e^{is\xi A^\top}P=Ce^{is\xi F}$ and
Theorem~\ref{glob:global-theorem} give
\begin{equation}
 \|Pf\|_{\mathcal Y^-_{K,J}(A,B)}
 \asymp\|f\|_{\mathcal Y^-_{K,J}(F^\top,C^\top)}.
 \label{jo:space}
\end{equation}
The block-coordinate weight is obtained from \eqref{glob:weight} by
replacing $\Pi_k^A$ with $\Pi_k$ and $A_{\lambda,j}$ with
$C_\lambda N_\lambda^jP_\lambda/j!$. The fixed factorials may be
omitted without changing norm equivalence. We use the notation
$\mathcal Y^-_{K,J}$ in either coordinate system, with
\eqref{jo:space} understood.

In the construction below assume $0\notin\Lambda$ and fix $T>0$,
with $s=T-t$. The case of a zero speed in
Theorem~\ref{jo:localized-general} follows directly from whole-line observation. The remaining task is to place observation in
the $L^2$ framework of the abstract criterion. We construct a rational
normalization that compensates for both the Kalman degeneracy at zero
and the growth caused by the nilpotent parts at high frequencies.
Although its inverse need not have finite propagation, the normalized
observation kernel has exponential decay, which supplies the required
localization estimates.
\subsection{High-frequency layers and a rational normalization}
The high-frequency decomposition used here is an increasing filtration of each generalized
eigenspace. It differs from the low-frequency Kalman filtration, and the two families of
projectors need not commute.
\begin{lemma}
\label{jo:high-layers}
For $0\leq j\leq r_\lambda$, put
\[
 \mathcal F_{\lambda,j}
=\bigcap_{\ell=j}^{r_\lambda-1}\ker(C_\lambda N_\lambda^\ell),\quad
 \mathcal F_{\lambda,r_\lambda}=E_\lambda,
\]
and let $H_{\lambda,j}$ be the orthogonal projector onto
\[
 E_{\lambda,j}=\mathcal F_{\lambda,j+1}\cap\mathcal F_{\lambda,j}^\perp,\quad 0\le j<r_\lambda.
\]

These projectors give an orthogonal decomposition of $\mathbb C^n$.  Moreover,
\begin{equation}
 C_\lambda N_\lambda^\ell H_{\lambda,j}=0\quad(\ell>j),
 \quad C_\lambda N_\lambda^j|_{E_{\lambda,j}}
\text{ is injective}.
 \label{jo:high-identities}
\end{equation}
There are constants $c,C>0$ such that, for $|\xi|\geq1$,
\begin{equation}
c\sum_{\lambda,j}|\xi|^{2j}|H_{\lambda,j}z|^2\leq\sum_{\lambda,\ell}|\xi|^{2\ell}
|C_\lambda N_\lambda^\ell P_\lambda z|^2
 \leq C\sum_{\lambda,j}|\xi|^{2j}|H_{\lambda,j}z|^2.
 \label{jo:high-equivalence}
\end{equation}
\end{lemma}
\begin{proof}
The subspaces $\mathcal F_{\lambda,j}$ increase with $j$, and
$\mathcal F_{\lambda,0}=\{0\}$ by \eqref{jo:block-kalman}.  This proves the orthogonal decomposition.  The vanishing in \eqref{jo:high-identities} follows
from $E_{\lambda,j}\subset\mathcal F_{\lambda,j+1}$.  If a vector in
$E_{\lambda,j}$ is also killed by $C_\lambda N_\lambda^j$, it belongs to $\mathcal F_{\lambda,j}$ and must be zero.
For completeness, set
$L_R=\sum_{\lambda,j}R^{-j}H_{\lambda,j}$, $R\geq1$, and consider the stacked map whose $(\lambda,\ell)$ component is
\[
 R^\ell C_\lambda N_\lambda^\ell P_\lambda L_Rz
=\sum_{j\geq\ell}R^{\ell-j}
C_\lambda N_\lambda^\ell H_{\lambda,j}z.
\]
As $R\to\infty$, this map converges in matrix norm to the map with components
$C_\lambda N_\lambda^\ell H_{\lambda,\ell}z$.  The latter is injective, by
\eqref{jo:high-identities} and the orthogonal decomposition.  Thus the stacked
map is bounded above and below uniformly for large $R$.  On a remaining compact interval $1\leq R\leq R_0$, the same conclusion follows by continuity and
\eqref{jo:block-kalman}.  Substituting $z=L_R^{-1}w$ proves
\eqref{jo:high-equivalence}.
\end{proof}
Choose an integer $M\ge1$ such that $2M+1>K$, and define
\begin{equation}
 \varphi(\xi)=\frac{\xi^{2M+1}}{(1+\xi^2)^M},\quad
 H(\xi)=\sum_{\lambda,j}(1+i\varphi(\xi))^jH_{\lambda,j},\quad
 D(\xi)=\sum_{k=0}^K\Bigl(\frac{i\xi}{1+i\xi}\Bigr)^k\Pi_k.
 \label{jo:normalizer-def}
\end{equation}
Our normalization is
\begin{equation}
 Q(\xi)=D(\xi)H(\xi).
 \label{jo:Q}
\end{equation}
The order of the factors in \eqref{jo:Q} matters: the two systems of
projectors need not commute. Each factor is harmless where the other
one acts: $H=I+O(|\xi|^{2M+1})$ near zero, whereas $D$ and $D^{-1}$ are
bounded for $|\xi|\ge1$. In the proof of Lemma~\ref{rev:normalizer} this
gives $DHD^{-1}=I+O(|\xi|^{2M+1-K})$ near zero and $|DHz|\asymp|Hz|$ at
infinity. With the reverse order, the same argument would require
$H(D-I)H^{-1}\to0$ at infinity. When the Kalman projectors mix the
high-frequency layers, this term is only $O(|\xi|^{r-2})$, with
$r=\max_\lambda r_\lambda$, and it need not tend to zero. For instance,
let $\lambda\ne\mu$ be nonzero reals and take
\[
 F=\begin{pmatrix}\lambda&1&0\\0&\lambda&0\\0&0&\mu\end{pmatrix},
 \qquad C=(1,0,1).
\]
Then $K=2$, $\operatorname{Ran}\Pi_2$ is spanned by
$(1,\mu-\lambda,-1)^\top$, and
$H(\xi)=\operatorname{diag}(1,1+i\varphi(\xi),1)$. Thus $H(\xi)$ and
$\Pi_2$ do not commute for $\xi\ne0$, and the $(2,1)$ entry of
$H(D-I)H^{-1}$ tends to $-(\mu-\lambda)/(2+(\mu-\lambda)^2)\ne0$.
This example shows why the small-perturbation argument cannot simply
be reversed; noncommutation alone does not disprove norm equivalence
for $HD$.

\begin{lemma}
\label{rev:normalizer}
The matrix $Q(\xi)$ is invertible for real $\xi\ne0$, is smooth on the real
axis, and has rational entries.  For Schwartz data,
\begin{equation}
 \|Q(D_x)f\|_{L^2}\asymp\|f\|_{\mathcal Y_{K,J}^{-}}.
 \label{jo:Q-equivalence}
\end{equation}
Consequently $Q(D_x)$ extends to an isomorphism from
$\mathcal Y_{K,J}^{-}$ onto $L^2(\mathbb R;\mathbb C^n)$.
\end{lemma}

\begin{proof}
For real $\xi$, the number $1+i\varphi(\xi)$ does not vanish.  All the projectors in $H$ are orthogonal, so $H$ and $H^{-1}$ have rational entries
without real poles.  Also
\[
 H(\xi)=I+O(|\xi|^{2M+1})\quad(\xi\to0),
 \quad \varphi(\xi)=\xi+O(|\xi|^{-1})\quad(|\xi|\to\infty).
\]
Near zero, $\|D(\xi)\|$ is bounded and
$\|D(\xi)^{-1}\|\leq C|\xi|^{-K}$.  Hence
\[
 D(\xi)H(\xi)D(\xi)^{-1}
       =I+O(|\xi|^{2M+1-K}).
\]
It follows that $|Q(\xi)z| \asymp |D(\xi)z|$ near zero,
with constants uniform in $z$, and the latter norm is uniformly equivalent to the low-frequency part of the $\mathcal Y^-_{K,J}$ norm.  On $|\xi|\geq1$, both $D$ and $D^{-1}$ are uniformly
bounded.  Orthogonality of the high-frequency layers gives
\[
 |H(\xi)z|^2
   =\sum_{\lambda,j}(1+\varphi(\xi)^2)^j|H_{\lambda,j}z|^2
\asymp\sum_{\lambda,j}|\xi|^{2j}|H_{\lambda,j}z|^2.
\]
Lemma~\ref{jo:high-layers} gives the required high-frequency equivalence.
Continuity and invertibility handle any compact annulus between the two
estimates.  This proves \eqref{jo:Q-equivalence} by Plancherel.  Finally, the
range of $Q(D_x)$ on Schwartz functions is dense in $L^2$: it contains every
function whose Fourier transform is smooth, compactly supported, and supported
away from zero.  Completion proves the last assertion.
\end{proof}
\subsection{The normalized observation kernel}

For a Schwartz datum $f$, set $u=Q(D_x)f$. In the normalized variable $u$,
the observation $Ce^{sF\partial_x}f$ is represented by the multiplier
\begin{equation}
 M_s(\xi)=Ce^{is\xi F}Q(\xi)^{-1}
=Ce^{is\xi F}H(\xi)^{-1}D(\xi)^{-1},
 \quad 0\leq s\leq T,\quad \xi\ne0.
 \label{jo:Ms}
\end{equation}
We will show that $M_s(D_x)$ is a finite sum of characteristic translations and a convolution operator with an exponentially decaying kernel. This
provides the representation needed to apply Proposition~\ref{abs:concrete}
to the weighted observations in the next subsection.

To fix the normalization of convolution kernels, write
\[
 (\mathcal F_0k)(\xi)
 :=\int_{\mathbb R}e^{-ix\xi}k(x)\,dx
 =\sqrt{2\pi}\,\widehat k(\xi).
\]
Thus $k*u=(\mathcal F_0k)(D_x)u$, and the symbol $e^{is\lambda\xi}$
acts as the translation $u(x)\mapsto u(x+\lambda s)$.
Two elementary identities illustrate the kernel contributions below:
\[
 \begin{aligned}
  \mathcal F_0\bigl(\mathbf1_{(-s,0)}\bigr)(\xi)
    =\frac{e^{is\xi}-1}{i\xi},\quad
  \mathcal F_0\bigl(e^{-x}\mathbf1_{(0,\infty)}(x)\bigr)(\xi)
    =\frac{1}{1+i\xi}.
 \end{aligned}
\]
In the first identity, the quotient extends to the value $s$ at zero:
the cancellation produces a kernel supported on a bounded interval.
In the second, the nonreal pole produces an exponentially decaying kernel.

\begin{lemma}
\label{rev:kernel}
The singularity of $M_s$ at zero is removable, uniformly for $s\in[0,T]$,
and
\begin{equation}
 M_s(0)=\sum_{k=0}^{K}\frac{s^k}{k!}CF^k\Pi_k.
 \label{jo:Ms-zero}
\end{equation}
Define
\begin{equation}
 V_\lambda(s)=\sum_{j=0}^{r_\lambda-1}
       \frac{s^j}{j!}C_\lambda N_\lambda^jH_{\lambda,j}.
 \label{jo:principal-V}
\end{equation}
There are measurable matrix kernels $k_s:\mathbb R\to\mathbb C^{m\times n}$
and constants $c,C_T>0$ such that, for every Schwartz input $u$,
\begin{equation}
 \begin{aligned}
  M_s(D_x)u(x)
   =\sum_{\lambda\in\Lambda}V_\lambda(s)u(x+\lambda s)
     +\int_{\mathbb R}k_s(x-y)u(y)\,dy,\quad
  |k_s(x)|\leq C_Te^{-c|x|}.
 \end{aligned}
 \label{jo:kernel-representation}
\end{equation}
The constants are independent of $s\in[0,T]$. This representation defines
uniformly bounded operators on both $L^2$ and $L^2_{\mathrm{uloc}}$.
\end{lemma}

\begin{proof}
\textit{Step 1: removal of the apparent singularity at zero.}
Put $p=2M+1>K$. By \eqref{jo:low-identities},
\[
 Ce^{is\xi F}\Pi_k
  =(i\xi)^k\left(\frac{s^k}{k!}CF^k\Pi_k
          +\xi R_{s,k}(\xi)\right),
\]
where $R_{s,k}$ is analytic near zero, with uniform bounds for
$s\in[0,T]$. Since
\[
 D(\xi)^{-1}
  =\sum_{k=0}^{K}\left(\frac{1+i\xi}{i\xi}\right)^k\Pi_k,
\]
these factors $(i\xi)^k$ cancel the poles of $D^{-1}$.
Hence $Ce^{is\xi F}D(\xi)^{-1}$ is analytic at zero and its value
there is the right-hand side of \eqref{jo:Ms-zero}.

The additional factor $H^{-1}$ does not change this value. Indeed,
$H(\xi)^{-1}-I=\xi^pR_H(\xi)$ with $R_H$ analytic near zero, whereas
$D^{-1}$ has a pole of order at most $K$. Therefore
\[
 Ce^{is\xi F}\bigl(H(\xi)^{-1}-I\bigr)D(\xi)^{-1}
   =\xi^{p-K}\widetilde R_s(\xi),
\]
with $\widetilde R_s$ analytic near zero. Since $p>K$, this term
vanishes at zero. Thus the full symbol $M_s$ extends analytically
across zero, with all local bounds uniform for $s\in[0,T]$.

\medskip
\noindent\textit{Step 2: extraction of the characteristic translations.}
Nilpotency and \eqref{jo:high-identities} give
\begin{equation}
 \begin{aligned}
 M_s(\xi)
  =\sum_{\lambda\in\Lambda}e^{is\lambda\xi}R_{\lambda,s}(\xi),\quad
 R_{\lambda,s}(\xi)
  =\left[\sum_{j=0}^{r_\lambda-1}\sum_{\ell=0}^{j}
     \frac{s^\ell(i\xi)^\ell}{\ell!}
     (1+i\varphi(\xi))^{-j}
     C_\lambda N_\lambda^\ell H_{\lambda,j}\right]D(\xi)^{-1}.
 \end{aligned}
 \label{jo:exact-rational}
\end{equation}
Since $\varphi(\xi)=\xi+O(|\xi|^{-1})$ and
$D(\xi)^{-1}=I+O(|\xi|^{-1})$, only the terms with $\ell=j$
contribute to the limit at infinity. Consequently,
\[
 R_{\lambda,s}(\xi)=V_\lambda(s)+O(|\xi|^{-1})
 \quad (|\xi|\to\infty),
\]
uniformly in $s$. The terms $e^{is\lambda\xi}V_\lambda(s)$
give the translations in \eqref{jo:kernel-representation}.
The remainders $R_{\lambda,s}-V_\lambda(s)$ are rational and vanish
at infinity. Their poles lie in a fixed finite set consisting
of zero and the nonreal poles of $H^{-1}$; the pole at zero has order
at most $K$. Their partial-fraction coefficients depend continuously
on $s$ and are uniformly bounded on $[0,T]$.

\medskip
\noindent\textit{Step 3: cancellation of the poles at zero.}
For $K\geq1$, group these contributions across all characteristic speeds:
\[
 Z_s(\xi)=\sum_{\lambda\in\Lambda}e^{is\lambda\xi}
             \sum_{d=1}^{K}\frac{Z_{\lambda,d}(s)}{(i\xi)^d}.
\]
All other terms are analytic at zero, so Step~1 implies that $Z_s$
is analytic there as well. This cancellation must be used before
identifying its convolution kernel.

Define
\[
 z_s(x)=\sum_{\lambda\in\Lambda}\sum_{d=1}^{K}
       Z_{\lambda,d}(s)\frac{(x+\lambda s)_+^{d-1}}{(d-1)!},
\]
where $x_+^0=\mathbf1_{(0,\infty)}(x)$.
This function vanishes to the left of all points $-\lambda s$.
To their right, it is a polynomial whose coefficient of $x^m$,
$0\leq m\leq K-1$, is
\[
 \frac1{m!}\sum_{\lambda\in\Lambda}\sum_{d=m+1}^{K}
   Z_{\lambda,d}(s)\frac{(\lambda s)^{d-1-m}}{(d-1-m)!}.
\]
This is $1/m!$ times the coefficient of $(i\xi)^{-m-1}$ in the
Laurent expansion of $Z_s$, hence it vanishes. Thus $z_s$ is
supported in the convex hull of $\{-\lambda s:\lambda\in\Lambda\}$,
and both its support and its size are uniformly bounded for $s\in[0,T]$.

To verify its Fourier symbol, differentiate in distributions:
\[
 \begin{aligned}
  \partial_x^Kz_s
    =\sum_{\lambda\in\Lambda}\sum_{d=1}^{K}
      Z_{\lambda,d}(s)\,\partial_x^{K-d}\delta_{-\lambda s},\quad
  (i\xi)^K(\mathcal F_0z_s)(\xi)
=(i\xi)^KZ_s(\xi).
 \end{aligned}
\]
Hence $\mathcal F_0z_s=Z_s$ for $\xi\ne0$ and, by continuity,
also at zero. This identifies the contribution of the poles at zero
with a uniformly bounded, compactly supported kernel.
For $K=0$, no poles at zero occur and we set $Z_s=z_s=0$.

\medskip
\noindent\textit{Step 4: exponential decay of the remaining kernels.}
For a nonreal pole $\zeta$ and an integer $\nu\geq1$, the kernel
with symbol $(\xi-\zeta)^{-\nu}$ is
\[
 \kappa_{\zeta,\nu}(x)=
 \begin{cases}
  i^\nu\mathbf1_{(0,\infty)}(x)
       \dfrac{x^{\nu-1}}{(\nu-1)!}e^{i\zeta x},
       &\operatorname{Im}\zeta>0,\\[4pt]
  (-i)^\nu\mathbf1_{(-\infty,0)}(x)
       \dfrac{(-x)^{\nu-1}}{(\nu-1)!}e^{i\zeta x},
       &\operatorname{Im}\zeta<0.
 \end{cases}
\]
The identity $\mathcal F_0\kappa_{\zeta,\nu}=(\xi-\zeta)^{-\nu}$
follows from the Laplace integral. Each kernel is a polynomial times
a decaying exponential on one half-line. Multiplication by
$e^{is\lambda\xi}$ replaces it by
$\kappa_{\zeta,\nu}(x+s\lambda)$.
There are finitely many poles, their distances from the real axis
are positive, and their orders are fixed. The coefficients and the
shifts are uniformly bounded for $s\in[0,T]$. Decreasing the
exponential rate absorbs all polynomial factors. Adding these kernels
to $z_s$ gives $k_s$ and proves \eqref{jo:kernel-representation},
including the uniform exponential bound.

Finally, $\sup_{s\in[0,T]}\|k_s\|_{L^1}<\infty$, so Young's
inequality gives the $L^2$ bound. Translations are uniformly bounded
on $L^2_{\mathrm{uloc}}$, since a translated unit interval meets at
most two unit cells. Minkowski's inequality then bounds convolution
by a constant times $\|k_s\|_{L^1}$ on this space. Together with the
uniform bounds on $V_\lambda(s)$, this proves the last assertion.
\end{proof}

\subsection{Characteristic occupation and the Gramian}
Let
\begin{equation}
 0\le a\le1,\quad a\in W^{1,\infty}(\mathbb R),\quad
 \inf_{\lambda\in\Lambda}\inf_{y\in\mathbb R}
 \int_0^T a(y-\lambda s)\,ds\ge\mu>0.
 \label{jo:weight-assumption}
\end{equation}
The occupation integrals are continuous, so essential infima give the same condition. The translation hull $\mathfrak H(a)$, taken in the locally uniform topology, is compact and translation invariant; all its elements retain the same occupation and Lipschitz bounds.
For $b\in\mathfrak H(a)$, set
\[
 (\Gamma_bu)(s,x)=\sqrt{b(x)}M_s(D_x)u(x),
 \quad G_b=\Gamma_b^*\Gamma_b.
\]
By Lemma~\ref{rev:kernel}, these observations are uniformly bounded on $L^2$ and also act locally on $L^2_{\mathrm{uloc}}$.

\begin{lemma}
\label{jo:principal}
The family $\Gamma_b$, $b\in\mathfrak H(a)$, satisfies (H1)--(H5). Its Gramian decomposes as
\begin{equation}
 G_b=\mathsf D_b+\mathsf K_b,\quad
 \mathsf D_b(y)=\sum_{\lambda\in\Lambda}\int_0^T
 b(y-\lambda s)V_\lambda(s)^* V_\lambda(s)\,ds.
 \label{jo:gramian-decomposition}
\end{equation}
The matrices $\mathsf D_b$ are uniformly Lipschitz and satisfy
\begin{equation}
 d_0I\le\mathsf D_b(y)\le d_1I
 \label{jo:D-coercivity}
\end{equation}
for constants $d_0,d_1>0$ independent of $b$ and $y$.
The remainder has a matrix kernel with the uniform bound
\begin{equation}
 |\mathsf K_b(x,y)|\le C_Te^{-c|x-y|}.
 \label{jo:gramian-kernel-bound}
\end{equation}
\end{lemma}
\begin{proof}
We first establish coercivity of the matrix in
\eqref{jo:gramian-decomposition}. Each $V_\lambda$ acts on the
domain block $E_\lambda$, so $\mathsf D_b$ is block diagonal.
For $z\in E_\lambda$, the polynomial formula
\[
 V_\lambda(s)z=\sum_{j=0}^{r_\lambda-1}
 \frac{s^j}{j!}C_\lambda N_\lambda^jH_{\lambda,j}z
\]
and Lemma~\ref{aux:polynomial-weight}, applied to
$b(y-\lambda s)$, give
\[
 \int_0^T b(y-\lambda s)|V_\lambda(s)z|^2\,ds
 \ge c_0\sum_j|C_\lambda N_\lambda^jH_{\lambda,j}z|^2
 \ge c_1|z|^2,
\]
where the last inequality follows from \eqref{jo:high-identities}.
The constants are uniform in $b$ and $y$, and there are finitely
many speeds. This proves the lower bound in
\eqref{jo:D-coercivity}; the upper and Lipschitz bounds follow
directly from the formula for $\mathsf D_b$.

Now apply Proposition~\ref{abs:concrete} with
\[
 c_{\lambda,b}(s,x)=\sqrt{b(x)}V_\lambda(s),\quad
 k_b(s,x,y)=\sqrt{b(x)}k_s(x-y).
\]
Lemma~\ref{rev:kernel} supplies the representation and uniform bounds.
Locally uniform convergence $b_n\to b$ gives coefficient continuity
and local $L^2$ kernel continuity. Translation covariance follows from
convolution, and the principal matrix is precisely $\mathsf D_b$.
Proposition~\ref{abs:concrete} therefore gives (H1)--(H5). In particular,
its proof yields
\[
 \|\mathbf1_I(\mathsf K_{b_n}-\mathsf K_b)\mathbf1_J\|_{\mathrm{HS}}
 \longrightarrow0
 \quad\text{for all bounded }I,J.
\]
Compactness of $\mathfrak H(a)$ then gives collective compactness
uniformly over the hull. The bound~\eqref{jo:gramian-kernel-bound} supplies localization;
it is not, by itself, a proof of collective compactness.
For \eqref{jo:gramian-kernel-bound}, the cross kernel for $\lambda\ne\nu$
is bounded by $\|V_\lambda\|_\infty\|V_\nu\|_\infty/|\lambda-\nu|$
and vanishes when $|x-y|>T|\lambda-\nu|$. The remaining kernels contain
the exponentially decaying remainder from Lemma~\ref{rev:kernel}.
\end{proof}

\subsection{Generalized uniqueness on the hull}
The uniqueness argument is stated for normalized data in $L^2_{\rm uloc}$. It therefore excludes both translated concentration and spreading sequences. No assertion about a nonzero global $L^2$ profile is needed.

Let $\mathcal O_M(\mathbb R)$ denote the smooth functions whose derivatives
grow at most polynomially. For a matrix symbol $m$ with entries in
$\mathcal O_M$, $m(D_x)$ maps $\mathcal S$ to $\mathcal S$ and
$\mathcal S'$ to $\mathcal S'$, and $(m_1m_2)(D_x)=m_1(D_x)m_2(D_x)$ on
$\mathcal S'$. The next lemma makes explicit how the normalized datum
$u\in L^2_{\mathrm{uloc}}$ is related to a solution of the original
transport system. The lift need not be unique, but its observation is
independent of the choice. Since $Q^{-1}$ may have a pole at zero,
the lift cannot be defined by multiplying an arbitrary tempered
distribution by $Q^{-1}$.

\begin{lemma}[Distributional lift and observation identity]
\label{lem:lift}
\begin{enumerate}
\item[\textup{(i)}] The symbols $D$, $H$, $H^{-1}$, $Q$,
$Ce^{is\xi F}$, $M_s$ and $\mathcal F_0k_s$ have entries in
$\mathcal O_M$, with bounds uniform in $s\in[0,T]$. For $K\ge1$, the
symbols $D^{-1}$ and $Q^{-1}$ have a pole of order $K$ at $\xi=0$ and are
smooth elsewhere.
\item[\textup{(ii)}] For every $u\in L^2_{\mathrm{uloc}}(\mathbb R;\mathbb C^n)$
there is $f\in\mathcal S'(\mathbb R;\mathbb C^n)$ with $Q(D_x)f=u$.
\item[\textup{(iii)}] The kernel of $Q(D_x)$ on $\mathcal S'$ is
\[
 \mathcal N_Q=\Bigl\{\textstyle\sum_{k=1}^{K}g_k:\ g_k\text{ is a polynomial
 of degree at most }k-1\text{ with values in }\operatorname{Ran}\Pi_k\Bigr\},
\]
and $Ce^{sF\partial_x}g=0$ for every $g\in\mathcal N_Q$ and every
$s\in\mathbb R$. In particular, the lift in \textup{(ii)} is unique
exactly when $K=0$.
\item[\textup{(iv)}] For every lift $f$ as in \textup{(ii)} and every
$s\in[0,T]$,
\[
 Ce^{sF\partial_x}f=M_s(D_x)u\quad\text{in }\mathcal S'.
\]
The right-hand side is given by the absolutely convergent formula
\eqref{jo:kernel-representation}, and on each cell it coincides with the
uniformly local extension of Definition~\ref{def:localized}. The map
$s\mapsto Ce^{sF\partial_x}f$ defines a distribution on
$(0,T)\times\mathbb R$, equal to the locally square integrable function
$y(s,x)=M_s(D_x)u(x)$. In particular $\Gamma_bu=\sqrt b\,y$, with
$\Gamma_b$ understood as in Theorem~\ref{rev:hull}. If
$u=Q(D_x)P^{-1}q_T$ with $q_T\in\mathcal S$, then
$y(s,\cdot)=B^\top q(T-s,\cdot)$.
\end{enumerate}
\end{lemma}

\begin{proof}
(i) The entries of $D$ and $H^{-1}$ are rational functions without real
poles that are bounded at infinity, and those of $H$ have polynomial
growth. The entries of $Ce^{is\xi F}$ are exponentials $e^{is\lambda\xi}$
times polynomials in $s\xi$. By Lemma~\ref{rev:kernel}, $M_s$ is analytic
at zero and equals $\sum_\lambda e^{is\lambda\xi}R_{\lambda,s}(\xi)$, where
each $R_{\lambda,s}$ is rational and bounded at infinity; all derivatives
of $M_s$ are therefore bounded on $\{|\xi|\ge1\}$. Since $k_s$ decays
exponentially, $x^jk_s\in L^1$ for every $j$, so $\mathcal F_0k_s$ is smooth
with bounded derivatives. The last assertion follows from the formula for
$D^{-1}$ in the proof of Lemma~\ref{rev:kernel} and from $H(0)=I$.

(ii) For $k\ge1$, define
\[
 w_k(x)=\frac1{(k-1)!}\int_0^x(x-t)^{k-1}\Pi_ku(t)\,dt,
\]
with the oriented integral for $x<0$, and set
\[
v_k=(1+\partial_x)^k w_k,
\quad v_0=\Pi_0u,
\quad v=\sum_{k=0}^K v_k.
\]
The uniformly local bound gives
$|w_k(x)|\le C_k(1+|x|)^k\|u\|_{L^2_{\mathrm{uloc}}}$, so all
these distributions are tempered. Since
$\partial_x^kw_k=\Pi_ku$ and $(1+\partial_x)^{-k}$ is a regular
multiplier on $\mathcal S'$, we have
\[
 D(D_x)v_k
 =\partial_x^k(1+\partial_x)^{-k}v_k=\Pi_ku
 \quad(k\ge1),\qquad D(D_x)v_0=\Pi_0u.
\]
Thus $D(D_x)v=u$.
By (i), $H^{-1}$ acts on $\mathcal S'$, and
$f=H(D_x)^{-1}v$ satisfies $Q(D_x)f=u$.

(iii) Let $Q(D_x)g=0$ and set $h=H(D_x)g\in\mathcal S'$.
Projecting $D(D_x)h=0$ onto each Kalman layer and applying
$(1+\partial_x)^k$ gives
\[
 \Pi_0h=0,\qquad \partial_x^k\Pi_kh=0\quad(1\le k\le K).
\]
Hence $\Pi_kh$ is a polynomial of degree at most $k-1$.
The symbol $H^{-1}-I$ vanishes to order $2M+1>K$ at zero, so its
multiplier annihilates these polynomials. Thus $g=H(D_x)^{-1}h=h$.
Conversely, every such polynomial sum lies in $\ker Q(D_x)$, proving
$\ker Q(D_x)=\mathcal N_Q$. Finally,
\[
 Ce^{sF\partial_x}g_k=\sum_{j\ge0}\frac{s^j}{j!}CF^j\Pi_k\,\partial_x^jg_k.
\]
The terms with $j<k$ vanish by \eqref{jo:low-identities}, and those with
$j\ge k$ vanish because $\deg g_k\le k-1$.

(iv) The removable singularity established in
Lemma~\ref{rev:kernel} is important here. The identity
$M_s(\xi)Q(\xi)=Ce^{is\xi F}$ holds for $\xi\ne0$, and both sides are
continuous, so it holds on the whole real axis. By (i), it is an identity
of multipliers on $\mathcal S'$, including distributions whose Fourier
transform is supported at $\xi=0$. Applying it to $f$ gives
$Ce^{sF\partial_x}f=M_s(D_x)u$. Two lifts differ by an element of
$\mathcal N_Q$, which is invisible by (iii), so the observed distribution
does not depend on the choice of lift.

For $u\in L^2_{\rm uloc}$ the right-hand side of
\eqref{jo:kernel-representation} is an absolutely convergent integral.
By Fubini's theorem on each cell it coincides with the uniformly local
extension of Definition~\ref{def:localized}. It also coincides with the
action of the multiplier $M_s$ on $\mathcal S'$: since $k_s$ decays
exponentially, it is a rapidly decreasing distribution, and
$k_s*u=(\mathcal F_0k_s)(D_x)u$ in $\mathcal S'$. The map
$s\mapsto e^{sF\partial_x}f$ is smooth from $\mathbb R$ into $\mathcal S'$,
since $\partial_s$ acts through the symbol $i\xi F$, so it defines a
distribution on $(0,T)\times\mathbb R$. The kernel representation makes
$y$ jointly measurable and locally square integrable in space-time, and
the two distributions coincide. Hence $\Gamma_bu=\sqrt b\,y$:
the weight multiplies the locally square integrable observation,
not the possibly distributional lift. The last
assertion follows from $C=B^\top P$ and
$e^{sF\partial_x}P^{-1}=P^{-1}e^{sA^\top\partial_x}$.
\end{proof}

The analyticity of $M_s$ at the origin is used in (iv) on distributions
whose Fourier transform is supported at $\xi=0$, such as constants,
which belong to $L^2_{\mathrm{uloc}}$. This is where the choice $2M+1>K$
enters: the correction $H^{-1}-I$ vanishes to order $2M+1$, which exceeds
the order $K$ of the pole of $D^{-1}$. For a constant $u_0$ one has
$M_s(D_x)u_0=M_s(0)u_0$, which is used at the end of the proof of
Lemma~\ref{rev:uniqueness}.

\begin{lemma}
\label{rev:uniqueness}
For every $b\in\mathfrak H(a)$,
\begin{equation}
 u\in L^2_{\rm uloc},\quad \Gamma_bu=0
\quad \Longrightarrow \quad u=0.
 \label{jo:generalized-invisibility}
\end{equation}
\end{lemma}

\begin{proof}
Let $f\in\mathcal S'(\mathbb R;\mathbb C^n)$ be a lift of $u$ given by
Lemma~\ref{lem:lift}\textup{(ii)}, and put
\[
q(s)=e^{sF\partial_x}f,
\quad y(s)=Cq(s).
\]
By Lemma~\ref{lem:lift}\textup{(iv)}, $y(s)=M_s(D_x)u$, the function $y$
is locally square integrable in space-time, and $\Gamma_bu=\sqrt b\,y$.
The choice of the lift does not affect $y$; we do not claim that $f$
itself is uniquely determined by $u$.

Let $O_b=\{x:b(x)>0\}$ and
$U_b=(0,T)\times O_b$.
Since $\Gamma_bu=\sqrt b\,y=0$, restriction to compact
subsets where $b$ is bounded below by a positive number
gives $y=0$ in $\mathcal D'(U_b)$.

The equation is $\partial_sq=F\partial_xq$.  For $0\leq j\leq K$,
\[
 CF^j\partial_x^Kq=\partial_s^j\partial_x^{K-j}y=0
\quad\text{in }\mathcal D'(U_b).
\]
The finite-dimensional injectivity \eqref{jo:kalman} therefore gives
\begin{equation}
 \partial_x^Kq=0\quad\text{in }\mathcal D'(U_b).
 \label{jo:local-polynomial}
\end{equation}
For the generalized eigenspace component,
\[
q_\lambda(s,x)=e^{sN_\lambda\partial_x}f_\lambda(x+\lambda s),
\quad f_\lambda=P_\lambda f.
\]
The inverse nilpotent exponential is a finite differential operator.  Hence
\begin{align*}
 \partial_x^Kf_\lambda(x+\lambda s)&=\sum_{j=0}^{r_\lambda-1}
\frac{(-s)^j}{j!}N_\lambda^j
\partial_x^{K+j}q_\lambda(s,x)=0\text{ in }\mathcal D'(U_b),
\end{align*}
where the second equality follows by spatial differentiation of
\eqref{jo:local-polynomial}.

The map $(s,x)\mapsto x+\lambda s$ sends $U_b$ onto the whole real line. Indeed, \eqref{jo:weight-assumption} guarantees, for every $y$, a time $s\in(0,T)$ with $b(y-\lambda s)>0$.  Lemma~\ref{aux:pullback}, applied componentwise, therefore gives
$\partial_x^K f_\lambda=0$ on $\mathbb R$.  For $K=0$ this already says $f=0$. Otherwise each $f_\lambda$, and therefore $f$, is a polynomial distribution of degree at most $K-1$.

The multiplier $Q$ is smooth at zero, so $u=Q(D_x)f$ is also a polynomial. A polynomial in $L^2_{\rm uloc}(\mathbb R)$ is constant; write $u=u_0$.
Lemma~\ref{rev:kernel} now gives
\[
 M_s(D_x)u_0=M_s(0)u_0
=\sum_{k=0}^{K}\frac{s^k}{k!}CF^k\Pi_ku_0.
\]
The weight $b$ is not identically zero.  Since $\Gamma_bu_0=0$, this polynomial
in $s$ is zero on $(0,T)$.  Every coefficient $CF^k\Pi_ku_0$ is zero.
Injectivity in \eqref{jo:low-identities} gives $\Pi_ku_0=0$ for every $k$,
and hence $u_0=0$.  This proves \eqref{jo:generalized-invisibility}.
\end{proof}
\subsection{Localized observability and bounded-gap sets}
\begin{theorem}
\label{jo:localized-general}
Assume the Kalman condition \eqref{jo:kalman} and $\sigma(A)\subset\mathbb R$.
Let $T>0$ and let $a:\mathbb R\to[0,1]$ be Lipschitz. There exist
$c_T,C_T>0$ such that every terminal datum $q_T\in\mathcal Y_{K,J}^{-}$
and its adjoint solution satisfy
\begin{equation}
 c_T\|q_T\|_{\mathcal Y_{K,J}^{-}}^2
 \leq\int_0^T\int_{\mathbb R}a(x)|B^\top q(t,x)|^2\,dx\,dt
 \leq C_T\|q_T\|_{\mathcal Y_{K,J}^{-}}^2
 \label{jo:localized-observability}
\end{equation}
if and only if
\[
 \min_{\lambda\in\sigma(A)}\inf_{y\in\mathbb R} m_{\lambda,a}(y)>0,
 \qquad m_{\lambda,a}(y)=\int_0^T a(y-\lambda s)\,ds .
\]
If $0\in\sigma(A)$, this condition is equivalent to $\inf_{\mathbb R}a>0$.
The upper bound holds for every such weight. For fixed $A,B,T$, the
two bounds are uniform on the hull $\mathfrak H(a)$ whenever the
condition holds. More generally they are uniform on any
translation-invariant family of $[0,1]$-valued weights that is compact
for locally uniform convergence and has a common Lipschitz bound and a
common positive lower bound on the occupation integrals.
\end{theorem}

\begin{proof}
The upper bound follows from $a\le1$, Plancherel's theorem and
Theorem~\ref{glob:global-theorem}. For necessity, fix a nonzero
eigenvector $A^\top v=\lambda v$. Kalman's condition gives
$B^\top v\ne0$. For $q_T=v\phi$, $\phi\in C_c^\infty(\mathbb R)$,
the observed energy is
$|B^\top v|^2\int m_{\lambda,a}(y)|\phi(y)|^2\,dy$.
The whole-line upper bound gives
\[
 T|B^\top v|^2\|\phi\|_2^2\le C_T^{\rm glob}\|v\phi\|_{\mathcal Y^-_{K,J}}^2.
\]
Thus the lower bound in \eqref{jo:localized-observability} implies
$m_{\lambda,a}\ge c_T T/C_T^{\rm glob}$ almost everywhere and, by
continuity, everywhere.

For sufficiency, let $\mu>0$ be a common lower bound for the occupation
integrals. If $0\in\sigma(A)$, then $m_{0,a}=Ta$ and $a\ge\mu/T$.
Conversely, $\inf a>0$ yields $m_{\lambda,a}\ge T\inf a$ at every
speed. In this case the lower estimate follows from
Theorem~\ref{glob:global-theorem} with the factor $\mu/T$.

Now suppose $0\notin\sigma(A)$. For Schwartz data put
$u=Q(D_x)P^{-1}q_T$. Lemma~\ref{rev:normalizer} and \eqref{jo:space}
give $\|u\|_2\asymp\|q_T\|_{\mathcal Y^-_{K,J}}$. For
$b\in\mathfrak H(a)$, Lemma~\ref{lem:lift}\textup{(iv)} shows that the
normalized observation satisfies
\[
 (\Gamma_bu)(s,x)=\sqrt{b(x)}\,B^\top q(T-s,x).
\]
Lemma~\ref{jo:principal} verifies the hypotheses on the observation
family, and Lemma~\ref{rev:uniqueness} proves uniqueness in
$L^2_{\rm uloc}$ for every hull element. Theorem~\ref{rev:hull}
therefore gives $\|\Gamma_bu\|_{\mathcal Z}\ge c_*\|u\|_2$ with $c_*>0$
independent of $b$. This proves the lower bound on Schwartz data.
The upper bound makes the observation continuous for the intrinsic
norm, so both estimates extend by density to $\mathcal Y^-_{K,J}$.
For a compact family as in the statement the same two lemmas apply
throughout the family, and Corollary~\ref{abs:compact-family} gives
a common lower bound.
\end{proof}

\begin{corollary}\label{rev:boundedgaps}
Assume \eqref{jo:kalman} and $\sigma(A)\subset\mathbb R\setminus\{0\}$.
Let $\Omega$ be as in \eqref{eq:intro-gap-parameters} and define
$T_\bullet$ by \eqref{eq:intro-geometric-time}. For $T>0$, there is a
constant such that
\begin{equation}
 \|q_T\|^2_{\mathcal Y^-_{K,J}}\le C_{T,\ell_\bullet,G^\bullet,A,B}
 \int_0^T\!\!\int_\Omega|B^\top q(t,x)|^2\,dx\,dt\label{jo:BG-observability}
\end{equation}
for every $q_T\in\mathcal Y^-_{K,J}$ if and only if $T>T_\bullet$.
When $\Omega=\mathbb R$ the constant depends only on $T,A,B$.
\end{corollary}

\begin{proof}
For $\Omega=\mathbb R$ apply Theorem~\ref{glob:global-theorem}.
Otherwise suppose $T>T_\bullet$ and apply Lemma~\ref{kal:compactweights} with the
characteristic-speed set $\Lambda$ and the parameters
$\ell_\bullet,G^\bullet,T$. It provides a common $\delta>0$,
a constant $m_0>0$, and weights
$a_\Omega\leq\mathbf1_\Omega$ such that
\[
\operatorname{Lip}(a_\Omega)\leq\delta^{-1},
\quad\int_0^T a_\Omega(y-\lambda s)\,ds\geq m_0.
\]
These weights belong to the compact
translation-invariant family $\mathfrak A$
defined in \eqref{kal:compactfamily}.
The uniform version of
Theorem~\ref{jo:localized-general} therefore applies.
Since $a_\Omega\leq\mathbf1_\Omega$, its weighted observation is bounded above by the observation on $\Omega$. This proves
\eqref{jo:BG-observability}. For $0<T\le T_\bullet$, failure follows
from Proposition~\ref{aux:necessity}, including equality when no gap
attains $G^\bullet$.
\end{proof}
\begin{corollary}
Under the hypotheses of Corollary~\ref{rev:boundedgaps}, system \eqref{eq:intro-primal} is
exactly controllable in $\mathcal Y^+_{K,J}$ at time $T>0$ if and only if
$T>T_\bullet$. In that case, for $y_0,y_d\in\mathcal Y^+_{K,J}$,
the minimum-norm control satisfies
\[
 \|u_{\min}\|^2_{L^2((0,T)\times\Omega)}\le C_{T,\ell_\bullet,G^\bullet,A,B}
 \|y_d-e^{-TA\partial_x}y_0\|^2_{\mathcal Y^+_{K,J}}.
\]
The observation space is optimal in the norm-domination sense of
\eqref{rev:domination}.
\end{corollary}

\begin{proof}
By the HUM equivalence \eqref{glob:HUM}, Corollary~\ref{rev:boundedgaps}
gives exactly the stated range of control times.
The minimum-control estimate \eqref{glob:correct-cost} yields the cost
bound, and whole-line observation gives the norm domination \eqref{rev:domination}.
\end{proof}
\subsection[Sharp observation cost for a 2 by 2 Jordan block]{Sharp observation cost for a $2\times2$ Jordan block}
\label{jo:cost}
As above, $F=A^\top$, $C=B^\top$, and $s=T-t$.
Let
\[
F=\begin{pmatrix}\lambda&1\\0&\lambda\end{pmatrix},
\quad C=(1,0),\quad \lambda\in\mathbb R\setminus\{0\}.
\]
Here $K=1$, both weights are $\operatorname{diag}(1,|\xi|^2)$, and
\[
\|f\|_{\mathcal Y^-_{K,J}}^2
=\|f_1\|_2^2+\|\partial_xf_2\|_2^2,
\quad
Ce^{sF\partial_x}f(x)=f_1(x+\lambda s)+s\partial_xf_2(x+\lambda s).
\]
The completion imposes no additional $L^2$ condition on $f_2$.
Write $\mathcal I_{T,\Omega}(f)
=\int_0^T\int_\Omega|Ce^{sF\partial_x}f|^2\,dx\,ds$.

\begin{proposition}
\label{rev:two-by-two}
Let $\Omega\in\mathcal G_{\ell,G}$. For $T>G/|\lambda|$ set
$\mu_T=\min\{\ell,|\lambda|T-G\}/|\lambda|$. Then
\begin{equation}
\frac{\mu_T^3}{12(1+T^2)}
 \|f\|_{\mathcal Y^-_{K,J}}^2
\le \mathcal I_{T,\Omega}(f)
\le \left(T+\frac{T^3}{3}\right)
 \|f\|_{\mathcal Y^-_{K,J}}^2
\label{rev:two-by-two-estimate}
\end{equation}
for every $f\in\mathcal Y^-_{K,J}$.

Suppose now that the supremum $G^\bullet$ of the actual gap lengths is
positive, and let $\ell_\bullet$ be the infimum of the observed interval
lengths. Set $T_\bullet=G^\bullet/|\lambda|$ and $\delta=T-T_\bullet$.
Write $C_{\mathrm{obs}}(T,\Omega)$ for the least constant in
$\|f\|_{\mathcal Y^-_{K,J}}^2\le C_{\mathrm{obs}}\mathcal I_{T,\Omega}(f)$.
For $0<\delta<\ell_\bullet/|\lambda|$ this constant is exactly
\[
 C_{\mathrm{obs}}(T,\Omega)=\frac{1}{\lambda_{\min}(M_\delta)},
 \qquad M_d=\int_{T-d}^T\begin{pmatrix}1&s\\s&s^2\end{pmatrix}\,ds
 \quad(0<d\le T).
\]
In particular,
\[
 \lim_{T\downarrow T_\bullet}(T-T_\bullet)^3 C_{\mathrm{obs}}(T,\Omega)
 =12(1+T_\bullet^2).
\]
\end{proposition}

\begin{proof}
For Schwartz data set $v(y)=(f_1(y),\partial_xf_2(y))^\top$ and
$w_y(s)=\mathbf1_\Omega(y-\lambda s)$. Changing variables gives
\[
\mathcal I_{T,\Omega}(f)=\int_{\mathbb R}v(y)^*\mathsf M(y)v(y)\,dy,
\quad
\mathsf M(y)=\int_0^T w_y(s)\begin{pmatrix}1&s\\s&s^2\end{pmatrix}\,ds.
\]
Lemma~\ref{rev:intersection} gives $m:=\int_0^Tw_y\,ds\ge\mu_T$.
For $\bar s=m^{-1}\int_0^Tsw_y(s)\,ds$, the bound $0\le w_y\le1$
implies that $\{|s-\bar s|\le r\}$ contains at most $2r$ units of
weighted mass. The layer-cake identity therefore yields
\[
\int_0^T(s-\bar s)^2w_y(s)\,ds
\ge\int_0^{m/2}2r(m-2r)\,dr=\frac{m^3}{12}.
\]
Thus $\det\mathsf M(y)\ge m^4/12$ and
$\operatorname{tr}\mathsf M(y)\le m(1+T^2)$, so
$\lambda_{\min}(\mathsf M(y))\ge\det\mathsf M(y)/
\operatorname{tr}\mathsf M(y)\ge\mu_T^3/[12(1+T^2)]$.
Integration proves the lower bound. The upper bound follows from
$|v_1+sv_2|^2\le(1+s^2)|v|^2$ and $w_y\le1$; density gives both
estimates on the completion.

For the exact formula, let $E_y=\{s\in(0,T):y-\lambda s\in\Omega\}$.
Lemma~\ref{rev:intersection} with the actual parameters
$\ell_\bullet,G^\bullet$ gives $|E_y|\ge\delta$. For any
$z\in\mathbb C^2$ put $p(s)=|z_1+sz_2|^2$. Among measurable subsets of
$[0,T]$ of measure at least $\delta$, its integral is minimized on
an interval $I$ of length $\delta$. Indeed, when $z_2\ne0$, the
sublevel sets of $p$ in $[0,T]$ are intervals, and their lengths vary
continuously from zero to $T$. Choose $I=\{p\le h\}$ of length $\delta$.
On $E_y\setminus I$ the integrand is at least $h$, and on $I\setminus E_y$
it is at most $h$. For $z_2=0$ the assertion is immediate.

For $I=[u,u+\delta]$ the moment matrix on $I$ has determinant
$\delta^4/12$ and trace $\delta[1+(u+\delta/2)^2+\delta^2/12]$.
At fixed determinant the smaller eigenvalue decreases as the trace
increases. The trace is largest at $u=T-\delta$. Hence
$\int_{E_y}p(s)\,ds\ge\lambda_{\min}(M_\delta)|z|^2$.
Taking $z=v(y)$ and integrating yields
$C_{\mathrm{obs}}(T,\Omega)\le1/\lambda_{\min}(M_\delta)$.

For the opposite bound choose gaps $(\alpha_n,\beta_n)$ of lengths
$g_n\to G^\bullet$ and $0<\varepsilon_n<g_n$ tending to zero.
If $\lambda>0$, take $J_n=(\beta_n-\varepsilon_n,\beta_n)$, and if
$\lambda<0$, take $J_n=(\alpha_n,\alpha_n+\varepsilon_n)$.
For $y\in J_n$, the path $y-\lambda s$ lies in the gap for
$0<s<(g_n-\varepsilon_n)/|\lambda|$. Thus $E_y\subset[T-d_n,T]$,
where $d_n=T-(g_n-\varepsilon_n)/|\lambda|\to\delta$.
Let $z_n$ be a unit eigenvector for $\lambda_{\min}(M_{d_n})$ and choose
$h_n\in C_c^\infty(J_n)$ with $\|h_n'\|_2=1$. The datum
$f_n=(z_{n,1}h_n',z_{n,2}h_n)$ has unit intrinsic norm and observed
energy at most $\lambda_{\min}(M_{d_n})$. Passing to the limit proves
equality, including when the supremum of the gap lengths is not attained.
Finally $\det M_\delta=\delta^4/12$ and
$\lambda_{\max}(M_\delta)/\delta\to1+T_\bullet^2$ as $T\downarrow
T_\bullet$, which gives the asserted limit.
\end{proof}

Here $C_{\mathrm{obs}}$ compares squared norms. By
\eqref{glob:correct-cost}, the corresponding optimal minimal-control
norm cost, for a unit final displacement, is $C_{\mathrm{obs}}^{1/2}$;
its blow-up rate is therefore $(T-T_\bullet)^{-3/2}$.

\section{Translating observation regions}
\label{mov:section}
The intrinsic spaces are unchanged for the two classes of translating
regions considered below.

\subsection{Slowly translating observation sets}

\begin{theorem}
\label{thm:slow-translations-iff}
Assume the Kalman condition and $\sigma(A)\subset\mathbb R\setminus\{0\}$, and
let $\Omega$ be as in \eqref{eq:intro-gap-parameters}. Fix $T>0$ and let
$h\in W^{1,\infty}(0,T)$ satisfy $\|h'\|_{L^\infty(0,T)}\le v_0<\lambda_\bullet$,
with $\lambda_\bullet$ as in \eqref{eq:intro-geometric-time}.
Use the Lipschitz representative of $h$ on $[0,T]$, and define
\[
\Omega_t:=\Omega+h(t),\quad
L_\lambda:=|\lambda T-h(T)+h(0)|,\quad
L_\bullet:=\min_{\lambda\in\sigma(A)}L_\lambda.
\]
Then the following assertions are equivalent:
\begin{enumerate}
\item $L_\bullet>G^\bullet $;
\item there exists $c>0$ such that the adjoint solution with any
terminal datum $q_T\in\mathcal Y^-_{K,J}$ satisfies
\[
c\|q_T\|_{\mathcal Y^-_{K,J}}^2
\leq\int_0^T\int_{\Omega_t}|B^\top q(t,x)|^2\,dx\,dt;
\]
\item the system
$\partial_t y+A\partial_x y=B\mathbf1_{\Omega_t}u$
is exactly controllable in $\mathcal Y^+_{K,J}$ in time $T$, with controls in
$L^2((0,T)\times\mathbb R;\mathbb C^m)$.
\end{enumerate}
The observation in the second assertion is interpreted by completion
from Schwartz terminal data. Its upper bound by
$C_T\|q_T\|_{\mathcal Y^-_{K,J}}^2$ always holds.

In particular,
$T>\frac{G^\bullet }{\lambda_\bullet-v_0}$
is sufficient. The observation constants are uniform for fixed
$A,B,T,v_0$ on families satisfying
$\ell_\bullet\geq\ell_0>0$ and $L_\bullet-G^\bullet \geq\eta>0$.
The same statements include $\Omega=\mathbb R$, with $G^\bullet =0$.
\end{theorem}

The bound on $h'$ makes each relative characteristic
$t\mapsto\lambda t-h(t)$ monotone. Its travelled distance is therefore
$|\lambda T-h(T)+h(0)|$, which explains the endpoint form of the criterion.

\begin{proof}
The whole-line case follows from Theorem~\ref{glob:global-theorem} and
HUM. Otherwise use the notation of Section~\ref{jo:section} and put
$\lambda^\bullet=\max_{\lambda\in\sigma(A)}|\lambda|$.

\emph{Occupation and the hull.}
The maps $\psi_\lambda(t)=\lambda t-h(t)$ are monotone and bi-Lipschitz,
with $|\lambda|-v_0\leq|\psi_\lambda'|\leq|\lambda|+v_0$ a.e., and
image length $L_\lambda$. If $L_\bullet>G^\bullet $, choose
$0<\varepsilon<\frac14\min\{\ell_\bullet,L_\bullet-G^\bullet \}$ and set
$a_\Omega(x)=\min\{1,\varepsilon^{-1}\operatorname{dist}(x,\mathbb R\setminus\Omega)\}$,
$b(s,x)=a_\Omega(x-h(T-s))$. The intervals on which $a_\Omega=1$
have lengths at least $\ell_\bullet-2\varepsilon$ and gaps at most
$G^\bullet +2\varepsilon$. Lemma~\ref{rev:intersection} and the change
of variable $z=\psi_\lambda(t)$ give
\[
 \int_0^T b(s,y-\lambda s)\,ds\geq
 \mu:=\frac{\min\{\ell_\bullet-2\varepsilon,
 L_\bullet-G^\bullet -2\varepsilon\}}{\lambda^\bullet +v_0}>0.
\]
The bounds $0\leq b\leq1$, $\operatorname{Lip}_x b\leq\varepsilon^{-1}$
and $\operatorname{Lip}_s b\leq v_0\varepsilon^{-1}$ make its full spatial
hull compact in $C_{\mathrm{loc}}([0,T]\times\mathbb R)$; every hull
weight $\beta$ retains this occupation bound.

Define $\Gamma_\beta u(s,x)=\sqrt{\beta(s,x)}M_s(D_x)u(x)$. The proofs of Lemmas~\ref{jo:principal} and~\ref{rev:uniqueness}
apply with $b(x)$ replaced by $\beta(s,x)$. Proposition~\ref{abs:concrete} allows the time-dependent
coefficients $c_{\lambda,\beta}=\sqrt\beta V_\lambda$ and kernels
$k_\beta(s,x,y)=\sqrt{\beta(s,x)}k_s(x-y)$; Lemma~\ref{aux:polynomial-weight}, applied to $s\mapsto\beta(s,y-\lambda
s)$, and Lemma~\ref{jo:high-layers} give the coercivity of $\sum_\lambda\int_0^T\beta(s,y-\lambda
s)V_\lambda(s)^*V_\lambda(s)\,ds$; and in the uniqueness argument the set $U_b$ is replaced by
the open set $U_\beta=\{(s,x):\beta(s,x)>0\}$, from which every $y$ is reached along each
characteristic family by the occupation bound. Lemma~\ref{aux:pullback} applies on $U_\beta$.
In the last step of Lemma~\ref{rev:uniqueness}, the polynomial $s\mapsto M_s(0)u_0$ vanishes on a nonempty open time interval, since $U_\beta$ is open and nonempty; hence it vanishes identically, and \eqref{jo:low-identities} gives $u_0=0$.
Theorem~\ref{rev:hull}, Lemma~\ref{rev:normalizer} and
$b(s,x)\le\mathbf 1_{\Omega+h(T-s)}(x)$ now give the lower bound; whole-line admissibility gives
the upper bound and extension by completion.

\emph{Necessity, including the critical time.}
For $Fv=\lambda v\neq0$ with $v\neq0$, Kalman gives $Cv\neq0$.
Terminal data $v\varphi$ have norm equivalent to $\|\varphi\|_2$ by
whole-line observability, while their observed energy is
$|Cv|^2\int m_\lambda(y)|\varphi(y)|^2\,dy$, where
\[
 m_\lambda(y)=\int_0^T
 \mathbf1_\Omega(y-\lambda(T-t)-h(t))\,dt
 =\int_{I_\lambda}\mathbf1_\Omega(y-z)\rho_\lambda(z)\,dz.
\]
Here $I_\lambda$ is the image of $(0,T)$ under $z_\lambda(t)=\lambda(T-t)+h(t)$,
and $\rho_\lambda=|(z_\lambda^{-1})'|$ there, extended by zero elsewhere. The length of $I_\lambda$ is
$L_\lambda$, and
$(|\lambda|+v_0)^{-1}\leq\rho_\lambda\leq(|\lambda|-v_0)^{-1}$.
Here the bounds hold on $I_\lambda$. The zero extension belongs to
$L^1(\mathbb R)$, so $\|m_\lambda(\cdot+r)-m_\lambda\|_\infty
\le\|\rho_\lambda(\cdot+r)-\rho_\lambda\|_1\to0$. If $L_\lambda<G^\bullet $,
placing $y-I_\lambda$ inside a larger gap makes $m_\lambda(y)=0$.
If $L_\lambda=G^\bullet $, then $y-I_\lambda$ is an interval of length
$G^\bullet$ with left endpoint $y-\sup I_\lambda$, so
$m_\lambda(y)\le(|\lambda|-v_0)^{-1}\mathfrak m_{G^\bullet}(y-\sup I_\lambda)$ with
$\mathfrak m_{G^\bullet}$ as in Lemma~\ref{aux:critical}. That lemma gives
$\inf m_\lambda=0$, whether or not $G^\bullet$ is attained.
In either case its essential infimum is zero, and normalized smooth
packets supported where it is small contradict observability.

\emph{Uniformity and controllability.}
Choose $\varepsilon<\frac14\min\{\ell_0,\eta\}$ uniformly and replace
$\ell_\bullet,L_\bullet-G^\bullet $ in $\mu$ by $\ell_0,\eta$.
All continuous weights with these Lipschitz and occupation bounds form
a compact translation-invariant family. The preceding argument and
Corollary~\ref{abs:compact-family} give uniform constants.
Also $L_\bullet \geq(\lambda_\bullet -v_0)T$, proving the sufficient time bound.
Section~\ref{glob:section} gives uniform boundedness of
$B\mathbf1_{\Omega_t}:L^2\to \mathcal Y^+_{K,J}$; strong continuity of translations
gives strong measurability. The bounded endpoint map has adjoint
$q_T\mapsto\mathbf1_{\Omega_t}B^\top q(t)$. The closed range theorem
and the reversible group therefore identify the lower observation
bound with exact controllability between arbitrary states.
\end{proof}

\subsection{Translations at constant speed}
For affine translations, a change of frame removes the speed restriction
and also permits stationary characteristics of the original system.

\begin{corollary}
\label{cor:affine-translations-iff}
Assume the Kalman condition and $\sigma(A)\subset\mathbb R$, and let $\Omega$
be as in \eqref{eq:intro-gap-parameters}; $\Omega=\mathbb R$ is allowed, with
$G^\bullet=0$.

For arbitrary $h_0,v\in\mathbb R$ and $T>0$, set
\[
\Omega_t:=\Omega+h_0+vt,\quad
d(v):=\min_{\lambda\in\sigma(A)}|\lambda-v|.
\]
Observability in $\mathcal Y^-_{K,J}$ on the moving set, and equivalently exact
controllability in $\mathcal Y^+_{K,J}$ in time $T$, hold if and only if
\[
 G^\bullet =0\quad\text{or}\quad T\,d(v)>G^\bullet .
\]
More explicitly, when $G^\bullet >0$ and $v\notin\sigma(A)$ the
necessary and sufficient condition is
\[
T>\frac{G^\bullet }{\min_{\lambda\in\sigma(A)}|\lambda-v|}.
\]
When $G^\bullet >0$ and $v\in\sigma(A)$, observability fails for
every positive time. When $G^\bullet =0$, it holds for every $T>0$
and every $v$, including resonant velocities.
\end{corollary}

\begin{proof}
Set $r(t,z)=q(t,z+h_0+vt)$. Then
$\partial_t r+(A^\top-vI)\partial_zr=0$ and
$r(T,z)=q_T(z+h_0+vT)$; changing variables identifies the moving
observation of $q$ with the stationary observation of $r$ on $\Omega$.
For $F=A^\top$, $C=B^\top$ and $F_v=F-vI$, the binomial formula gives
$\bigcap_{j=0}^k\ker(CF_v^j)=\bigcap_{j=0}^k\ker(CF^j)$ for every $k$.
Thus the Kalman filtration and its projectors are unchanged. The
generalized spectral projectors are also unchanged, and
$F_v-(\lambda-v)I=F-\lambda I$ on $E_\lambda$, so the \(C_\lambda N_\lambda^j\)
coincide. Both intrinsic weights, and hence both completed energy
spaces, are identical. Spatial translation is an isometry, giving
$\|r(T)\|_{\mathcal Y^-_{K,J}(A-vI,B)}=\|q_T\|_{\mathcal Y^-_{K,J}(A,B)}$.

If $G^\bullet =0$, the observation set equals $\mathbb R$ up to countably
many endpoints, and Theorem~\ref{glob:global-theorem} applies.
If $d(v)>0$, Corollary~\ref{rev:boundedgaps} for $A-vI$ proves
sufficiency, while Proposition~\ref{aux:necessity} proves failure for
$T d(v)\leq G^\bullet $, including equality when the supremal gap is
not attained. If $G^\bullet >0$ and $d(v)=0$, take an eigenvector
$Fe=ve$ and a nonzero smooth packet $r_T=e\varphi$ supported in a gap.
This solution is stationary in the moving frame and invisible, while
its whole-line observation $T|Ce|^2\|\varphi\|_2^2$ is positive.
Finally, the change of frame is also an isometry on controls and on
both endpoint state spaces, so HUM gives the stated controllability
equivalence.
\end{proof}

Taking $v=0$ recovers the stationary criterion; for
$|v|<\min_{\lambda\in\sigma(A)}|\lambda|$, this agrees with
Theorem~\ref{thm:slow-translations-iff}.
\section[Transport with nonlocal dispersion in arbitrary dimension]{Transport with nonlocal dispersion in $\mathbb R^d$}
\label{sec:higher-dim}

Let $d\ge1$, $V\in\mathbb R^d$ and $J\in L^1(\mathbb R^d;\mathbb C)$.
For $u^0\in L^2(\mathbb R^d)$, let $u$ solve
\begin{equation}\label{eq:nonlocal}
  \partial_tu+V\cdot\nabla u=J*u-u\quad\text{in }(0,T)\times\mathbb R^d,
  \quad u(0)=u^0 .
\end{equation}
For $J\ge0$ with $\int J=1$, the operator $J*u-u$ is the classical
nonlocal dispersal operator; no sign
or normalization is needed below. For a weight $a:\mathbb R^d\to[0,1]$,
put
\begin{equation}\label{eq:nonlocal-occupation}
  \mu_a:=\inf_{z\in\mathbb R^d}\int_0^Ta(z+tV)\,dt .
\end{equation}
Balls $B_\varepsilon(z)$ are Euclidean in this subsection.

\begin{theorem}\label{thm:nonlocal}
Let $L\ge0$ and let $a:\mathbb R^d\to[0,1]$ satisfy $\operatorname{Lip}(a)\le L$.
\begin{enumerate}
\item[(i)] If $\mu_a\ge\mu>0$, then
\begin{equation}\label{eq:nonlocal-obs}
  \|u^0\|_{L^2(\mathbb R^d)}^2
  \le C\int_0^T\!\!\int_{\mathbb R^d}a(x)\,|u(t,x)|^2\,dx\,dt
  \quad(u^0\in L^2(\mathbb R^d)),
\end{equation}
with $C$ depending only on $d,T,V,J,L$ and $\mu$.
\item[(ii)] If $\mu_a=0$, then \eqref{eq:nonlocal-obs} fails for every $C$.
\end{enumerate}
\end{theorem}

\begin{proof}
\emph{Step 1. Representation.} The operator $\mathcal Jf:=J*f-f$ is
bounded on $L^2$ and commutes with translations, and
\begin{equation}\label{eq:nonlocal-semigroup}
  e^{t\mathcal J}f=e^{-t}f+K_t*f,\quad
  K_t:=e^{-t}\sum_{k\ge1}\frac{t^k}{k!}\,J^{*k}.
\end{equation}
The series converges in $L^1$, uniformly for $0\le t\le T$, and
\begin{equation}\label{eq:nonlocal-kernel}
  |K_t|\le F_T:=\sum_{r\ge1}\frac{T^r}{r!}|J|^{*r},
  \qquad \|F_T\|_1\le e^{T\|J\|_1}-1 .
\end{equation}
Since $w(t,x):=u(t,x+tV)$ solves $\partial_tw=\mathcal Jw$, we have
$u(t,x)=(e^{t\mathcal J}u^0)(x-tV)$. Thus, for $b:\mathbb R^d\to[0,1]$,
\begin{equation}\label{eq:nonlocal-observation}
  \Gamma_bu^0(t,x):=b(x)^{1/2}u(t,x)
  =e^{-t}b(x)^{1/2}u^0(x-tV)
  +\int_{\mathbb R^d}b(x)^{1/2}K_t(x-tV-y)\,u^0(y)\,dy .
\end{equation}

\emph{Step 2. Verification of (H1)--(H5).}
Let $\mathfrak A_{L,\mu}$ consist of the functions
$b:\mathbb R^d\to[0,1]$ satisfying $\operatorname{Lip}(b)\le L$ and
$\mu_b\ge\mu$. With locally uniform convergence and the translation
 action $\theta_h b=b(\cdot+h)$, this family is compact and metrizable
by Arzel\`a--Ascoli. The occupation bound passes to the limit on
each compact segment $\{z+tV:0\le t\le T\}$. Formula
\eqref{eq:nonlocal-observation} gives covariance. Local uniform
convergence of $b$ implies convergence of the observed cell energies
by dominated convergence.

The transport term has finite range between cells. For the convolution
term, put $D=1+T|V|_2$. Young's inequality bounds the block from
$Q_j$ to $Q_k$, at each time, by
$\int_{k-j+[-D,D]^d}F_T$. These integrals are summable in $k-j$ by
finite overlap, which proves (H3) after integration in time.
The principal multiplication operator $D_b$ has coefficient
$d_b(z)=\int_0^T e^{-2t}b(z+tV)\,dt\in[e^{-2T}\mu,T]$, which
depends continuously on $b$ locally uniformly. This gives (H4).

To verify (H5), choose $J_m\in C_c^\infty$ with $J_m\to J$ in $L^1$.
For any $\alpha>0$, let $M_m=\|e^{\alpha|\cdot|_2}J_m\|_1$ and
$L_m=\|e^{\alpha|\cdot|_2}J_m\|_\infty$.
The triangle inequality and Young's inequality give
$\|e^{\alpha|\cdot|_2}|J_m|^{*r}\|_\infty\le L_m M_m^{r-1}$.
Thus $|K^{(m)}_t(x)|\le C_{m,T}e^{-\alpha|x|_2}$ on $[0,T]$.
Proposition~\ref{abs:concrete} applies to these kernels with one
velocity $-V$ and $c_{1,b}(t,x)=e^{-t}\sqrt{b(x)}$. In particular,
the bounded-window restrictions of
$(\Gamma_b^{(m)})^*\Gamma_b^{(m)}-D_b$ are collectively compact in $b$.
The convolution generators converge in operator norm on $L^2$, hence
their exponentials converge in norm uniformly on $[0,T]$. It follows
that $\Gamma_b^{(m)}\to\Gamma_b$ and their Gramians converge in norm,
uniformly in $b$. Since $D_b$ is independent of $m$, the restricted
remainders for $J$ are uniform norm limits of collectively compact
families and are collectively compact by
Lemma~\ref{lem:cc-tools}\textup{(b)}.
This proves (H5), and Lemma~\ref{lem:gramR} gives (R).

\emph{Step 3. Uniqueness on $\mathfrak A_{L,\mu}$.} Let
$b\in\mathfrak A_{L,\mu}$ and $w\in L^2_{\mathrm{uloc}}$ with
$\Gamma_bw=0$. By \eqref{eq:nonlocal-kernel}, the series
\eqref{eq:cellwise-extension} defining the extension of $\Gamma_b$ can be summed
under the integral, so $\Gamma_bw(t,x)=b(x)^{1/2}g(t)(x-tV)$ with
$g(t):=e^{-t}w+K_t*w$. Young's inequality on the cells gives
$\|J*f\|_{L^2_{\mathrm{uloc}}}\le2^d\|J\|_{L^1}\|f\|_{L^2_{\mathrm{uloc}}}$,
so $g(t)=e^{-t}\sum_{k\ge0}\frac{t^k}{k!}J^{*k}*w$ extends to an entire
function of $t$ with values in $L^2_{\mathrm{uloc}}$.

Fix $y_0\in\mathbb R^d$. Since $\int_0^Tb(y_0+tV)\,dt\ge\mu>0$ and $b$ is
continuous, there are $t_0\in(0,T)$ and $\varepsilon>0$ such that
$(t_0-\varepsilon,t_0+\varepsilon)\subset(0,T)$ and $b(y+tV)>0$ whenever
$|t-t_0|<\varepsilon$ and $y\in B_\varepsilon(y_0)$. The change of
variables $x=y+tV$ at fixed $t$ shows that $g(t)(y)=0$ for almost every
$(t,y)\in(t_0-\varepsilon,t_0+\varepsilon)\times B_\varepsilon(y_0)$. By
Fubini's theorem, the entire function
$t\mapsto\mathbf 1_{B_\varepsilon(y_0)}g(t)\in L^2(B_\varepsilon(y_0))$
vanishes for almost every $t\in(t_0-\varepsilon,t_0+\varepsilon)$, hence
identically. At $t=0$ this gives $w=0$ on $B_\varepsilon(y_0)$. As $y_0$
is arbitrary, $w=0$.

\emph{Step 4. Conclusion of (i).} Corollary~\ref{abs:compact-family},
applied to $\{\Gamma_b\}_{b\in\mathfrak A_{L,\mu}}$, gives $c_*>0$,
depending only on $d,T,V,J,L,\mu$, such that
$\|\Gamma_bu^0\|_{\mathcal Z}\ge c_*\|u^0\|_{L^2}$ for all
$b\in\mathfrak A_{L,\mu}$ and $u^0\in L^2$. Since $a\in\mathfrak A_{L,\mu}$
and $\|\Gamma_au^0\|_{\mathcal Z}^2$ is the integral in
\eqref{eq:nonlocal-obs}, this proves (i) with $C=c_*^{-2}$.

\emph{Step 5. Necessity.} Let $\mu_a=0$. The function
$m(z):=\int_0^Ta(z+tV)\,dt$ is $TL$-Lipschitz. Given $k\ge1$, choose $z_k$
with $m(z_k)\le1/k$ and put $\varepsilon_k:=1/(k(TL+1))$, so that
$m\le2/k$ on $B_{\varepsilon_k}(z_k)$. Fix $\chi\in C_c^\infty(B_1(0))$
with $\|\chi\|_{L^2}=1$ and a unit vector $\eta\in\mathbb R^d$, and set
\[
  u^0_{k,N}(x):=\varepsilon_k^{-d/2}\chi\big((x-z_k)/\varepsilon_k\big)\,e^{iN\eta\cdot x}
  \qquad(N\ge1).
\]
Then $\|u^0_{k,N}\|_{L^2}=1$, and \eqref{eq:nonlocal-semigroup} together
with the change of variables $x=y+tV$ gives
\[
  \int_0^T\!\!\int_{\mathbb R^d}a|u|^2\,dx\,dt
  \le2\int m(y)\,|u^0_{k,N}(y)|^2\,dy
  +2T\sup_{0\le t\le T}\|K_t*u^0_{k,N}\|_{L^2}^2
  \le\frac4k+2T\sup_{0\le t\le T}\|K_t*u^0_{k,N}\|_{L^2}^2 .
\]
The Fourier multiplier of $f\mapsto K_t*f$ is
$e^{-t}\big(e^{t\widehat J_0(\xi)}-1\big)$, where
$\widehat J_0(\xi):=\int e^{-ix\cdot\xi}J(x)\,dx$ tends to zero as
$|\xi|\to\infty$ by the Riemann--Lebesgue lemma. This multiplier is
bounded by $e^{T\|J\|_{L^1}}+1$ and tends to zero as $|\xi|\to\infty$,
uniformly in $t\in[0,T]$. Since the Fourier transform of $u^0_{k,N}$ is
the translate by $N\eta$ of a fixed $L^2$ function, dominated convergence
shows that the last supremum tends to zero as $N\to\infty$, for fixed
$k$. Choosing $N=N_k$ large, we obtain unit data whose observation is at
most $(4+2T)/k$, so \eqref{eq:nonlocal-obs} fails.

\end{proof}
\begin{corollary}\label{cor:nonlocal-slabs}
Let $V\neq0$, $e:=V/|V|_2$ and $\Omega\in\mathcal G_{\ell,G}$, and consider the slab
\[
S_\Omega:=\{x\in\mathbb R^d:\ x\cdot e\in\Omega\}.
\]
For $\Omega=\bigcup_{j\in\mathbb Z}(a_j,b_j)$, let
$G^\bullet:=\sup_{j\in\mathbb Z}(a_{j+1}-b_j)\le G$
denote the supremum of its gap lengths, with $G^\bullet=0$ if $\Omega=\mathbb R$.
\begin{enumerate}
\item[(i)] For every $T>G^\bullet/|V|_2$ there is $C$, depending only on
$d,T,V,J,\ell,G^\bullet$, such that
\[
\|u^0\|_{L^2}^2\le C\int_0^T\!\!\int_{S_\Omega}|u|^2\,dx\,dt
\qquad (u^0\in L^2(\mathbb R^d)).
\]
For $T>G/|V|_2$, the constant can be chosen uniformly for all
$\Omega\in\mathcal G_{\ell,G}$.
\item[(ii)] If $G^\bullet>0$, this inequality fails for every
$0<T\le G^\bullet/|V|_2$.
\item[(iii)] The system
$\partial_ty+V\cdot\nabla y=J*y-y+\mathbf 1_{S_\Omega}h$ is exactly
controllable between arbitrary states of $L^2(\mathbb R^d)$ with
$h\in L^2((0,T)\times S_\Omega)$ if and only if $T>G^\bullet/|V|_2$.
\end{enumerate}
\end{corollary}
\begin{proof}
(i) Since $G^\bullet$ is the supremum of the gap lengths of $\Omega$, we have
$\Omega\in\mathcal G_{\ell,G^\bullet}$, and $T>G^\bullet/|V|_2$.
Lemma~\ref{kal:compactweights}, applied with $\Lambda=\{-|V|_2\}$ and the gap
bound $G^\bullet$, therefore provides the weight $a_\Omega$ of
\eqref{kal:minorant}. Put $a(x):=a_\Omega(x\cdot e)$. Then
$0\le a\le\mathbf 1_{S_\Omega}$,
$\operatorname{Lip}(a)\le\operatorname{Lip}(a_\Omega)$, and
\[
\int_0^T a(z+tV)\,dt=m_{-|V|_2,a_\Omega}(z\cdot e)\ge m_0>0
\qquad(z\in\mathbb R^d),
\]
where $m_0$ and $\operatorname{Lip}(a_\Omega)$ depend only on
$T,|V|_2,\ell,G^\bullet$. Theorem~\ref{thm:nonlocal}(i), applied with
$L=\operatorname{Lip}(a_\Omega)$ and $\mu=m_0$, and the bound
$a\le\mathbf 1_{S_\Omega}$ give the asserted inequality, with $C$ depending
only on $d,T,V,J,\ell,G^\bullet$. If $T>G/|V|_2$, the same lemma applied with
the common bound $G$ yields constants that are uniform over
$\mathcal G_{\ell,G}$.

(ii) Let
\[
m(z):=\int_0^T\mathbf 1_{S_\Omega}(z+tV)\,dt
=\frac{\bigl|(z\cdot e,\,z\cdot e+T|V|_2)\cap\Omega\bigr|}{|V|_2}.
\]
This function is $(2/|V|_2)$-Lipschitz, and its infimum is zero. Indeed, if
$T|V|_2<G^\bullet$, an interval of length $T|V|_2$ fits in a gap of
$\Omega$; if $T|V|_2=G^\bullet$, then
$m(z)=|V|_2^{-1}\mathfrak m_{G^\bullet}(z\cdot e)$, and Lemma~\ref{aux:critical}
gives $\inf m=0$, whether or not $G^\bullet$ is attained. Step~5 of the proof of
Theorem~\ref{thm:nonlocal}, with $a$ replaced by $\mathbf 1_{S_\Omega}$,
uses only these two properties of $m$, and yields unit data with vanishing
observations.

(iii) The operator $-V\cdot\nabla+\mathcal J$ generates a group on
$L^2(\mathbb R^d)$. Its adjoint is $V\cdot\nabla+\widetilde{\mathcal J}$,
where $\widetilde{\mathcal J}f:=\widetilde J*f-f$ and
$\widetilde J(x):=\overline{J(-x)}$. With control space
$L^2(S_\Omega)$, the control operator is extension by zero to
$L^2(\mathbb R^d)$, and its adjoint is restriction to $S_\Omega$.
After the time reversal
$t\mapsto T-t$, the adjoint evolution is therefore \eqref{eq:nonlocal} with
$(V,J)$ replaced by $(-V,\widetilde J)$, and
$\|\widetilde J\|_{L^1}=\|J\|_{L^1}$. Since $-V/|V|_2=-e$, the slab is
$S_\Omega=\{x:\ x\cdot(-e)\in-\Omega\}$, and $-\Omega\in\mathcal G_{\ell,G}$
has the same gap supremum $G^\bullet$. Parts (i) and (ii), applied to this
adjoint system, therefore characterize its observability by
$T>G^\bullet/|V|_2$. By Hilbert space duality
\cite{Lions1988,TucsnakWeiss2009}, surjectivity of the endpoint map is
equivalent to this observability inequality.
\end{proof}
\begin{example}[No transport drift]\label{rem:nonlocal}
If $V=0$, then $\mu_a=T\inf a$. Theorem~\ref{thm:nonlocal} therefore
gives weighted observability exactly when $\inf a>0$.
For instance, the Lipschitz weight $a(x)=\min\{1,|x|_2\}$ does not
give observability, even when $J$ is nonzero: nonlocal dispersion
alone does not provide a uniform observation estimate near a zero of
the weight.
\end{example}

\section{Anisotropic heat flux beyond the Gramian assumptions}
\label{sec:anisotropic-heat}
\label{new:sec-heat-flux}

Let $\mathsf A\in\mathbb R^{d\times d}$ be symmetric positive definite, with
eigenvalues $0<\alpha_{\min}\le\dots\le\alpha_{\max}$, and let
$u(t)=P_tf:=e^{t\nabla\cdot(\mathsf A\nabla)}f$ on $\mathbb R^d$. In an
anisotropic conductor, the energy identity
\begin{equation}\label{new:eq-aniso-energy}
  \|f\|_2^2=\|P_Tf\|_2^2+2\int_0^T\langle\mathsf A\nabla u(t),\nabla u(t)\rangle\,dt
\end{equation}
involves the quadratic form $\langle\mathsf A\nabla u,\nabla u\rangle$. Heat-flux
sensors, however, measure Fourier's flux $\mathbf q=-\mathsf A\nabla u$. We
observe this flux on $\mathbb R^d$ and the temperature on a measurable set
$\Omega$. The resulting observability statement is elementary. The point of
this section is structural: the measured outputs satisfy (R), but when
$\mathsf A$ has two distinct eigenvalues their Gramian admits no
decomposition of the form (H4)--(H5).

\begin{theorem}\label{new:thm-aniso-flux}
Fix $T>0$. There is a constant $C$ such that
\begin{equation}\label{new:eq-aniso-obs}
  \|f\|_2^2\le C\int_0^T\Big(\|\mathsf A\nabla u(t)\|_2^2
  +\|\mathbf 1_\Omega u(t)\|_2^2\Big)\,dt
  \qquad(f\in L^2(\mathbb R^d))
\end{equation}
if and only if $\Omega$ is thick, that is,
$\inf_{h\in\mathbb R^d}|\Omega\cap(h+[0,L]^d)|>0$ for some $L>0$. On the
weak-star hull $\mathcal W$ of $\mathbf 1_\Omega$, the observation family
\[
  \Gamma_bf=\big(-\mathsf A\nabla P_tf,\ b^{1/2}P_tf\big),\qquad b\in\mathcal W,
\]
satisfies (H1)--(H3) and (R). If $\mathsf A$ has two distinct eigenvalues, then
for no $b\in\mathcal W$ does $\Gamma_b^*\Gamma_b$ admit a decomposition
satisfying (H4)--(H5).
\end{theorem}

\begin{proof}
Let $\mathcal W$ be the weak-star closure of
$\{\mathbf1_\Omega(\cdot+h):h\in\mathbb R^d\}$ in $L^\infty$.
Its elements satisfy $0\le b\le1$. The hull is compact and metrizable
because $L^1(\mathbb R^d)$ is separable, and the action
$\theta_hb=b(\cdot+h)$ is jointly continuous by translation continuity
in $L^1$. Covariance follows from translation invariance of the heat flow.

\emph{Step 1: energy normalization.} Put
$\Gamma_b^0f=(\mathsf A^{1/2}\nabla P_tf,\ b^{1/2}P_tf)$. By
\eqref{new:eq-aniso-energy},
\[
(\Gamma_b^0)^*\Gamma_b^0=\tfrac12\,I+K_b^0,\qquad
K_b^0=-\tfrac12P_{2T}+\int_0^TP_t\,b\,P_t\,dt .
\]
The heat kernel is
$p_t(x)=(4\pi t)^{-d/2}(\det\mathsf A)^{-1/2}\exp(-x\cdot\mathsf A^{-1}x/4t)$,
and $\|\nabla P_t\|_{2\to2}\le(2e\,\alpha_{\min}t)^{-1/2}$. Hence
$\sup_{0\le b\le1}\|K^0_b\|_{L^2\to H^1}\le C_T$, and Rellich's theorem gives
(H4)--(H5) with $d_b=\tfrac12$. The Gaussian bounds for $p_t$ and $\nabla p_t$ bound
the blocks between cells at distance $r\ge2$ by $C_Te^{-c(r-1)^2/T}$, and the
energy identity bounds the remaining blocks. This proves (H3). By Lemma~\ref{lem:gramR},
$\Gamma^0_b$ satisfies (R).

\emph{Step 2: the measured flux.} We have $\Gamma_b=M\Gamma_b^0$ with
$M=\mathrm{diag}(-\mathsf A^{1/2},1)$. Lemma~\ref{new:lem-output}, applied with
$L=\mathrm{diag}(-\mathsf A^{-1/2},1)$, gives (H1), (H3) and (R) for $\Gamma_b$.
The flux channel does not depend on $b$. On a cell $Q$, the energy of the
temperature channel is the pairing of $b$ with
$\mathbf 1_Q\int_0^T|P_tf|^2\,dt\in L^1$, which is weak-star continuous. This
proves (H2); see Remark~\ref{new:rem-H2}(ii).

\emph{Step 3: uniqueness and thickness.} Let $\Gamma_bw=0$ with
$w\in L^2_{\rm uloc}$. Since $\mathsf A$ is invertible, $\nabla P_tw=0$, so
$P_tw$ is constant in space, and hence also in time by the heat equation. Since
$P_tw\to w$ in $\mathcal S'$, the datum $w$ is a constant $c$. Then
$b^{1/2}c=0$ forces $c=0$ unless $b=0$.
If $\Omega$ is thick, choose $L,\eta>0$ such that
$|\Omega\cap(h+[0,L]^d)|\ge\eta$ for every $h$.
Passing to weak-star limits preserves
$\int_{[0,L]^d}b\ge\eta$, so $0\notin\mathcal W$.
Conversely, if thickness fails, choose $h_n$ such that
$|\Omega\cap(h_n+[-n,n]^d)|<1/n$.
Then $\mathbf1_\Omega(\cdot+h_n)\to0$ in $L^1_{\rm loc}$;
the uniform bound by one also gives weak-star convergence in $L^\infty$.
Thus $0\in\mathcal W$. Theorem~\ref{rev:hull} now gives the equivalence.

\emph{Step 4: no Gramian decomposition.} The flux part of
$\Gamma_b^*\Gamma_b$ is the Fourier multiplier
\[
  m_T(\xi)=\frac{\xi\cdot\mathsf A^2\xi}{2\,\xi\cdot\mathsf A\xi}
  \Big(1-e^{-2T\xi\cdot\mathsf A\xi}\Big).
\]
Let $\mathsf Ae_i=\alpha_ie_i$, $i=1,2$, with $|e_i|_2=1$ and
$\alpha_1\ne\alpha_2$. Fix a unit vector $\varphi\in C_c^\infty(Q)$ in $L^2$,
where $Q$ is a bounded cube, and put $f_{n,i}=\varphi\,e^{in\,e_i\cdot x}$. These
sequences are weakly null. By Plancherel and dominated convergence, their flux
energies tend to $\alpha_i/2$. The temperature Gramian maps $L^2$ into $H^1$
uniformly. Its compression to $Q$ is therefore compact, so its
contribution tends to zero. Suppose that
$\Gamma_b^*\Gamma_b=D_b+K_b$ satisfies (H4)--(H5). Since $|f_{n,i}|=|\varphi|$,
both limits would equal $\int d_b|\varphi|^2$, because
$\mathbf 1_QK_b\mathbf 1_Q$ is compact. This gives $\alpha_1=\alpha_2$, a
contradiction.
\end{proof}

\begin{remark}\label{new:rem-isotropic}
If $\mathsf A=\alpha I$, then
$m_T(\xi)=\tfrac\alpha2\big(1-e^{-2T\alpha|\xi|_2^2}\big)$, and (H4)--(H5) hold
directly. The obstruction is created by the anisotropy of the medium, not by an
ad hoc choice of output weights. Measuring the plain gradient $\nabla u$
instead of the flux leads to the same conclusion, with limits $1/(2\alpha_i)$
and $M=\mathrm{diag}(\mathsf A^{-1/2},1)$.
\end{remark}

\appendix
\renewcommand{\thetheorem}{\Alph{section}.\arabic{theorem}}
\section{Polynomial coercivity on time weights}
\label{aux:section-polynomial}

\begin{lemma}\label{aux:polynomial-weight}
Let $r\geq1$, $T>0$, and let $a:(0,T)\to[0,1]$ be measurable, with
$m=\int_0^T a(s)\,ds>0$. For vectors $v_0,\ldots,v_{r-1}$ in a
finite-dimensional Hilbert space, set $P(s)=\sum_{j=0}^{r-1}s^jv_j$.
There is $c_r>0$, depending only on $r$, such that
\begin{equation}\label{aux:polynomial-bound}
 \int_0^T a(s)|P(s)|^2\,ds
 \geq c_r\frac{m^{2r-1}}{(1+T)^{2r-2}}
\sum_{j=0}^{r-1}|v_j|^2.
\end{equation}
Factorials in the coefficients may be absorbed in $c_r$.
\end{lemma}
\begin{proof}
Write $I=\int a|P|^2$. If $I=0$, then $P$ vanishes on a set of positive
Lebesgue measure and hence is identically zero, component by component.
Assume $I>0$ and let
$E=\{s\in(0,T):|P(s)|^2\leq2I/m\}.$
Since the weighted mass of its complement is at most $m/2$, one has
$\int_Ea\geq m/2$, and therefore $|E|\geq m/2$. Put $d=m/(4r)$. A greedy selection produces $r$ points $s_1,\ldots,s_r$ in $E$ at pairwise distances at least $d$: before the last selection, the union of the forbidden intervals has measure at most
$2(r-1)d<m/2$. The vector-valued Lagrange formula reads
\[
 P(s)=\sum_{i=1}^r P(s_i)\prod_{h\neq i}\frac{s-s_h}{s_i-s_h}.
\]
Every denominator has modulus at least $d^{r-1}$, while every coefficient
of a numerator is bounded by $C_r(1+T)^{r-1}$. Consequently
\[
 \sum_j|v_j|^2
 \leq C_r(1+T)^{2r-2}d^{-2r+2}\max_i|P(s_i)|^2
 \leq C_r(1+T)^{2r-2}m^{-2r+1}I.
\]
Rearrangement proves the assertion.
\end{proof}

\section{Pullback uniqueness along a family of characteristics}
\label{aux:section-pullback}

\begin{lemma}\label{aux:pullback}
Let $U\subset(0,T)\times\mathbb R$ be open, $\lambda\in\mathbb R$, and
$p_\lambda(s,x)=x+\lambda s$. If $p_\lambda(U)=\mathbb R$ and
$g\in\mathcal D'(\mathbb R)$ satisfies $g(x+\lambda s)=0$ on $U$,
then $g=0$.
\end{lemma}
\begin{proof}
Given $z_0\in\mathbb R$, choose $(s_0,x_0)\in U$ with
$x_0+\lambda s_0=z_0$. In the coordinates $(s,z)=(s,x+\lambda s)$,
the image of $U$ contains a rectangle $J\times V$ around $(s_0,z_0)$.
Choose $\rho\in C_c^\infty(J)$ with $\int\rho=1$. For
$\phi\in C_c^\infty(V)$ the function $\eta(s,x)=\rho(s)\phi(x+\lambda s)$
is supported in $U$. Testing there and changing variables yields
\[
 0=\langle g(x+\lambda s),\eta\rangle
   =\left\langle g,\phi\int_J\rho(s)\,ds\right\rangle=\langle g,\phi\rangle.
\]
Thus $g$ vanishes near each $z_0$, proving the claim componentwise.
\end{proof}

\section{Obstructions and the sharp geometric threshold}
\label{aux:section-obstructions}

\begin{proposition}\label{aux:necessity}
If the Kalman rank condition fails, no observation inequality for a
nondegenerate norm on Schwartz data can hold, even with observation
on the whole line. If $A$ has a zero eigenvalue and $\Omega$ has a
nonempty open gap, the same conclusion holds for localized observation
at every $T>0$.

For a real-spectrum Kalman pair with nonzero speeds and a geometry
\eqref{eq:intro-gap-parameters} with $G^\bullet>0$, observability in
$\mathcal Y^-_{K,J}$ fails whenever $0<T\le T_\bullet$, with $T_\bullet$ as in
\eqref{eq:intro-geometric-time}.
\end{proposition}

The critical case $T=T_\bullet$ rests on the following elementary fact,
which is also used in Theorem~\ref{thm:slow-translations-iff} and
Corollary~\ref{cor:nonlocal-slabs}.

\begin{lemma}[Critical length]
\label{aux:critical}
Let $\Omega$ be as in \eqref{eq:intro-gap-parameters} with $G^\bullet>0$,
and let $0<L\le G^\bullet$. The function
$\mathfrak m_L(y):=|(y,y+L)\cap\Omega|$ is $2$-Lipschitz and
$\inf_{y\in\mathbb R}\mathfrak m_L(y)=0$, whether or not $G^\bullet$ is attained.
\end{lemma}

\begin{proof}
The symmetric difference of $(y,y+L)$ and $(z,z+L)$ has measure at most
$2|y-z|$. If $(b_j,a_{j+1})$ is a gap of length $g_j$, then
$\mathfrak m_L(b_j)\le(L-g_j)_+$. Since $\sup_jg_j=G^\bullet\ge L$, the infimum
of these bounds is zero.
\end{proof}

\begin{proof}[Proof of Proposition~\ref{aux:necessity}]
Write $s=T-t$.

If Kalman fails, the Cayley--Hamilton theorem gives a nonzero vector
$v$ such that $B^\top (A^\top)^jv=0$ for every $j\geq0$.
Choose a nonzero scalar Schwartz function $h$ with compactly supported
Fourier transform. Then
\[
B^\top e^{is\xi A^\top }v\,\widehat h(\xi)=0
\quad (s\in[0,T],\ \xi\in\mathbb R),
\]
so the nonzero terminal datum $vh$ is invisible on the whole line.
If instead $A^\top v=0$ for some $v\neq0$, choose $h\in C_c^\infty(I)\setminus\{0\}$, where $I$ is an open gap of $\Omega$.
The solution with terminal datum $vh$ is stationary and invisible.

For the last assertion, assume the stated Kalman and spectral
hypotheses. Choose a real eigenvector $v\neq0$ such that $A^\top v=\lambda v$ and
$|\lambda|=\lambda_\bullet$.
Kalman implies $B^\top v\neq0$, and
$e^{s A^\top \partial_x}(vh)(x)=v\,h(x+\lambda s).$
Let $C_T^{\mathrm{glob}}>0$ be a constant in the whole-line upper
observation estimate. It follows that
\[
T|B^\top v|^2\|h\|_{L^2}^2
=
\int_0^T
\|B^\top e^{sA^\top \partial_x}(vh)\|_{L^2}^2\,ds
\leq
C_T^{\mathrm{glob}}
\|vh\|_{\mathcal Y_{K,J}^{-}}^2.
\]
Thus scalar packets with unit $L^2$ norm have intrinsic norms bounded
away from zero, without any frequency modulation.

For $y\in\mathbb R$, define
\[
J_\lambda(y):=\{y-\lambda s:0<s<T\},
\quad
m_\lambda(y):=
\int_0^T\mathbf1_\Omega(y-\lambda s)\,ds
=
\frac{|J_\lambda(y)\cap\Omega|}{|\lambda|}.
\]
Changing variables gives
\[
\int_0^T\!\!\int_\Omega
|B^\top e^{s A^\top \partial_x}(vh)(x)|^2\,dx\,ds
=
|B^\top v|^2\int_{\mathbb R}
m_\lambda(y)|h(y)|^2\,dy.
\]

If $0<T<T_\bullet$, there is an open gap of length greater than
$|\lambda|T$. The terminal points $y$ for which $J_\lambda(y)$
lies in that gap contain a nonempty open interval.
A nonzero smooth packet supported in this interval has zero
observation and positive intrinsic norm.

Now let $T=T_\bullet$, so that $|J_\lambda(y)|=G^\bullet $.
The left endpoint of $J_\lambda(y)$ is $y_-:=\min\{y,y-\lambda T\}$, so
$m_\lambda(y)=|\lambda|^{-1}\mathfrak m_{G^\bullet}(y_-)$ with $\mathfrak m_{G^\bullet}$
as in Lemma~\ref{aux:critical}. That lemma gives points $y_n$ with
$m_\lambda(y_n)\to0$, whether or not the supremal gap length is attained,
and
$|m_\lambda(y)-m_\lambda(z)|\leq\frac{2|y-z|}{|\lambda|}.$
Choose $\varepsilon_n\downarrow0$ and real functions
\[
h_n\in C_c^\infty(y_n-\varepsilon_n,y_n+\varepsilon_n),
\quad
\|h_n\|_{L^2}=1.
\]
On their supports,
\[
m_\lambda(y)
\leq
\delta_n:=
m_\lambda(y_n)+\frac{2\varepsilon_n}{|\lambda|}
\longrightarrow0.
\]
Consequently, the real Schwartz data $f_n:=vh_n$ satisfy
\[
\|f_n\|_{\mathcal Y_{K,J}^{-}}^2
\geq\frac{T|B^\top v|^2}{C_T^{\mathrm{glob}}}>0,
\quad\int_0^T\!\!\int_\Omega
|B^\top e^{s A^\top \partial_x}f_n|^2\,dx\,ds
\leq
|B^\top v|^2\delta_n\longrightarrow0.
\]
Dividing by their intrinsic norms gives unit-norm counterexamples to observability at $T=T_\bullet$, whether or not the supremal gap
length is attained.
\end{proof}
\section*{Acknowledgments}
The author would like to express his sincere gratitude to Enrique Zuazua for his valuable discussions, insightful suggestions, and continuous encouragement throughout the development of this work.

The author gratefully acknowledges the hospitality of the Chair for Dynamics, Control, Machine Learning and Numerics at Friedrich-Alexander-Universität Erlangen–Nürnberg, where this research was carried out during a research stay. The author also acknowledges the financial support of the Alexander von Humboldt Foundation through an Alexander von Humboldt Research Fellowship.

\end{document}